\documentclass[10pt]{amsart}

\usepackage{geometry}        
\usepackage{graphicx}
\usepackage{amssymb}
\usepackage{epstopdf}
\usepackage{amsmath}
\usepackage{color}
\usepackage{hyperref}
\usepackage{pdfsync}
\usepackage{tikz}
\usepackage{tikz-cd}
\usepackage{datetime2}
\usepackage[all]{xy}
\xyoption{matrix}
\xyoption{arrow}

\usepackage{amsfonts}
\usepackage{amsthm}
\usepackage{mathdots}
\usepackage{enumerate}
\usepackage[latin1]{inputenc}
\usepackage{graphicx}
\usepackage{latexsym}
\usepackage{accents}
\usepackage{mathrsfs}

\usetikzlibrary{arrows,decorations.pathmorphing,backgrounds,positioning,fit,petri}

\tikzset{help lines/.style={step=#1cm,very thin, color=gray},
help lines/.default=.5} 
\tikzset{thick grid/.style={step=#1cm,thick, color=gray},
thick grid/.default=1} 

\newtheorem{thm}{Theorem}[subsection]
\newtheorem{lem}[thm]{Lemma}
\newtheorem{cor}[thm]{Corollary}
\newtheorem{prop}[thm]{Proposition}

\theoremstyle{definition}
\newtheorem{defn}[thm]{Definition}
\newtheorem{eg}[thm]{Example}

\newtheorem{rem}[thm]{Remark}
\newtheorem{notation}[thm]{Notation}

\numberwithin{equation}{section}

\DeclareMathOperator{\Hom}{Hom}%
\DeclareMathOperator{\Ext}{Ext}%
\DeclareMathOperator{\End}{End}%
\newcommand{\kk}{\ensuremath{\Bbbk}}

\newcommand{\commentout}[1]{}

\newcommand{\cA}{\ensuremath{{\mathcal{A}}}}

\newcommand{\cB}{\ensuremath{{\mathcal{B}}}}

\newcommand{\cF}{\ensuremath{{\mathcal{F}}}}

\newcommand{\mmod}{{\rm mod}\hspace{0.3pt}}
\newcommand{\rad}{{\rm rad}\hspace{0.3pt}}

\newcommand{\undim}{\underline{\rm dim}\hspace{0.5pt}}
\newcommand{\ndp}{{\rm dp}}
\newcommand{\dt}{{\accentset{\hspace{1pt}\mbox{\large\bfseries .}}{}}}

\newcommand{\cdt}{\dt\hspace{4pt}}

\newcommand{\pdt}{{\hspace{.8pt}\dt\hspace{1.5pt}}}
\newcommand{\id}{{\rm id}\hspace{0.3pt}}
\def\Z{\hbox{$\mathbb{Z}$}}
\def\Da{\hbox{${\mathit\Delta}$}}
\def\Sa{\hbox{${\mathit\Sigma}$}}
\def\Oa{\hbox{${\mathit\Omega}$}}
\def\La{\hbox{${\mathit\Lambda}$}}
\def\Ga{\hbox{${\mathit\Gamma}$}}
\def\GaA{\hbox{$\mathit\Gamma_{\hspace{-1.2pt}\mathcal A}$}}
\def\GaH{\hbox{$\mathit\Gamma_{\hspace{.8pt}{\rm mod}H}$}}
\def\DaH{\hbox{$\mathit\Delta_{\hspace{-.2pt}H}$}}
\def\SaH{\hbox{$\mathit\Sigma_{\hspace{-.3pt}H}$}}
\def\QH{\hbox{$Q_{\hspace{-.8pt}H}$}}
\def\OQH{\hbox{$\hspace{2.5pt}\overline{\hspace{-1.5pt}Q}_{\hspace{-.9pt}H\hspace{-1pt}}$}}
\def\ZQH{\hbox{$\mathbb{Z}\hspace{.8pt}Q^{\rm op}_{\hspace{-1pt}H}$}}
\def\QoH{\hbox{$Q\hspace{-.8pt}_H^{\hspace{.8pt}\rm op}$}}

\title{\sc Representation theory of hereditary artin algebras of finite representation type}
\author{Shiping Liu}
\address{Department of Mathematics, Sherbrooke University, Sherbrooke, Quebec, Canada}
\email{shiping.liu@usherbrooke.ca}
\author{Gordana Todorov}
\address{Department of Mathematics, Northeastern University, Boston, MA 02115, USA}
\email{g.todorov@northeastern.edu}



\begin{document}

\begin{abstract}

Let $H$ be a hereditary artin algebra of finite representation type. 
We first determine all hammocks in the Auslander-Reiten quiver $\GaH$ of $\mmod H$, the category  of finitely generated left $H$-modules. This enables us to obtain an effective method to construct $\GaH$ by simply viewing the ext-quiver of $H$. As easy applications, we compute the numbers of non-isomorphic indecomposable objects in $\mmod H$ and the associated cluster category $\mathscr{C}_H$, as well as the nilpotencies of the radicals of $\mmod H\hspace{-.4pt},$ $\hspace{-.5pt} D^{\hspace{.5pt}b\hspace{-.6pt}}(\hspace{-.5pt}\mmod H\hspace{-.5pt})$ and $\mathscr{C}_H$.

\vspace{-20pt}

\end{abstract}

\maketitle

\section*{Introduction}

Since Gabriel's pivotal work on representations of quivers; see \cite{Gab}, the representation theory of hereditary artin algebras has been extensively studied over the past fifty years; see, for example, \cite{DRi, DR, PAu, Ri1, Zac}. The representation type of such an algebra is finite precisely when its ext-quiver is of Dynkin type; see \cite{ARS, DR, Gab}. In this case, the information on the module category is encoded in its Auslander-Reiten quiver. In the linearly oriented $\mathbb{A}_n$-case, the Auslander-Reiten quiver has a wing shape; see \cite[(2.6)]{LiY} and \cite[(3.3)]{Ri1}. In all other cases, the Auslander-Reiten quiver has only been roughly described; see \cite[(VIII.1.15)]{ARS} and \cite[(1.13)]{LiY}.

The main objective of this paper is to provide an effective method to construct the Auslander-Reiten quiver simply by viewing the ext-quiver of the algebra. This enables us to compute the number of non-isomorphic indecomposable objects in the module category and the associated cluster category, and also the nilpotency of the radicals of the module category, its bounded derived category and the associated cluster category. More details are outlined below.

In this introduction, $H$ denotes a hereditary artin algebra with a Dynkin extension quiver $\QH$; see (\ref{def_Ext_quiver}). Let $\mmod H$ be the category of finitely generated left $H$-modules and $\GaH$ be its Auslander-Reiten quiver; see (1.5). It is well-known that every $\tau$-orbit in $\GaH$ is finite and contains a projective module and an injective module; see \cite[(1.8)]{PAu}, and also \cite[(1.9)]{DR}. Moreover, $\GaH$ embeds as a convex valued translation quiver in the repetitive quiver $\mathbb{\Z}\DaH$, where $\DaH$ is the full valued subquiver of $\GaH$ generated by the projective modules; see \cite[(1.13)]{LiY}, and also \cite[(VIII.1.15)]{ARS}. Thus, in order to describe the precise shape of $\GaH$, it suffices to determine the injective module and the number of modules in the $\tau$-orbit of any projective module.

For giving a combinatorial criterion for a finite translation quiver to be the Auslander-Reiten quiver of a finite dimensional algebra over an algebraically closed field, Brenner introduced the notions of hammocks and hammock functions; see \cite{Bre}. We adapt and extend these notions to the artin setting; see (\ref{hammock_defn}) and (\ref{hamm_func}). After a thorough study of the valued translation quivers of tree type in Section 2, we describe all the hammocks in $\GaH$; see (\ref{hamm_neg_section}). This allows us to determine  which injective module lies in the $\tau$-orbit of any given projective module; 
see, for example, (\ref{An}), (\ref{Dn}) and (\ref{E6}).

To study representations of a valued graph, Dlab and Ringel defined a Coxeter transformations as a product of all reflections in the corresponding rational vector space; see \cite[Page 8]{DRi}. Alternatively, for any hereditary artin algebra, Auslander and Platzeck defined the Coxeter transformation as an special automorphism of the Grothendieck group of the module category; see (\ref{Cox_trans}) and \cite[Section 2]{PAu}. Since $\QH$ is of Dynkin type, the Coxeter transformation $C_H$ for $H$ is of finite order $|C_H|$; see \cite[(4.1)]{PAu}, called the {\it Coxeter order}. The Coxeter order has been computed for each Dynkin diagram with a particular orientation; see \cite[Pages 289, 290]{ARS}. We show that $|C_H|$ is independent of the orientation of $\QH$, so it is explicitly given no matter how $\QH$ is oriented; see (\ref{co_table}). Further, we shall pair the $\tau$-orbits of projective modules in $\GaH$ in such a way that the sum of the numbers of modules in the paired $\tau$-orbits is equal to $|C_H|$; see (\ref{order_c}), while the difference equals the difference of the numbers of arrows in two reduced walks in $\QH$ between the corresponding vertices; see (\ref{pi-index}). So, the number of modules in any $\tau$-orbit is expressed in terms of $|C_H|$ and the number of arrows in a specific reduced walk in $\QH$. This ultimately 
completes the picture of $\GaH$ in terms of $\QH$; see (\ref{ARQ_Dyn}).

Since $H$ is of finite representation type, it is interesting to compute the number of non-isomorphic indecomposable $H$-modules. Since the center $\kk$ of $H$ is a field; see \cite[(3.1)]{PAu}, $\mmod H$ is equivalent to the representation category of a $\kk$-linear species of $\QH$; see \cite[Theorem C]{DRi}. Hence, this number coincides with the number of positive roots of $\QH$; see \cite{DR, Gab}, which can be found in \cite{Cox}.
However, we directly obtain this number from our results, that is, half of $|C_H|$ multiplied by the number of non-isomorphic simple $H$-modules; see (\ref{rep_nb}).

A well-known result of Auslander says that an artin algebra $\La$ is representation-finite 
if and only if the radical of $\mmod \La$ is nilpotent; see \cite[(V.7.7)]{ARS} and \cite[(1.1)]{Sko}. In this case, it is interesting to compute the nilpotency of $\rad(\mmod \La)$, namely the least integer $r$ for which $\rad^r\hspace{-1pt}(\mmod \La)=0$. Using preprojective partitions, Zacharia proved that $|C_H|-2$ is the maximal length of chains of irreducible maps in $\mmod \La$ with a non-zero composite; see \cite[(4.11)]{Zac}, and consequently, $\rad(\mmod H)$ is of nilpotency $|C_H|-1$. Alternatively, using some nice properties of $\GaH;$ see (\ref{prep_cpt_path}), we show that the radicals of $\mmod H$ and $D^{\hspace{.5pt}b\hspace{-.6pt}}(\mmod H)$ have the same nilpotency $|C_H|-1$; see (\ref{Mod_nil}) and (\ref{Der_nil}).

In order to categorify Fomin and Zelevinsky's cluster algebra associated with a finite acyclic unvalued quiver \cite{FoZ}, Buan, Marsh, Reineke, Reiten and Todorov introduced the associated cluster category, which is an orbit category of the bounded derived category of finite dimensional representations of the quiver; see \cite{BMRRT}. The same construction yields the cluster category $\mathscr{C}_H$ associated with $H$. We show that the radical of $\mathscr{C}_H$ is also of nilpotency $|C_H|-1$, and the number of non-isomorphic indecomposable objects in $\mathscr{C}_H$ is half of $|C_H|+2$ multiplied by the number of non-isomorphic simple $H$-modules; see (\ref{clus_nil}), which coincides with the number of cluster variables of the cluster algebra associated with $\QH$; see \cite[(5.9.1)]{FWZ}. Finally, we want to express our gratitude to Platzeck and Zacharia for some helpful discussions.


\section{Preliminaries}

The objective of this section is to lay the foundation for this paper. Besides fixing some terminology and notation, we shall collect and also prove some results in the general context, which are needed later. Throughout this paper, $R$ denotes a commutative artinian ring. All algebras are $R$-algebras, and all categories are additive $R$-categories whose morphisms are composed from right to left.

\subsection{\sc Quivers} 
We start with laying the combinatorial foundation. Let $Q=(Q_0, Q_1)$ be a quiver, where $Q_0$ is the set of vertices and $Q_1$ is the set of arrows from a vertex, called the {\it start point}, to another vertex, called the {\it end point}. An arrow $\alpha$ from $x$ to $y$ is usually represented graphically as $x\to y$, and we write $x=s(\alpha)$ and $y=e(\alpha)$. Given $x\in Q_0$, write $x^+$ for the set of vertices $y$ such that there exists an arrow $x\to y$, and $x^-$ for the set of vertices $z$ such that there exists an arrow $z\to x$. We call $x$ a {\it source} or {\it sink } if $x^-$ or $x^+$ is empty, respectively. With each vertex $x\in Q_0$, one associates a {\it trivial path} $\varepsilon_x$ of length 0 with $s(\varepsilon_x)=e(\varepsilon_x)=x$.
A {\it path} of length $r>0$ is a sequence $\eta=\alpha_r \cdots \alpha_1$, where $\alpha_i\in Q_1$ such that $e(\alpha_i)=s(\alpha_{i+1})$ for $1\le i<r$; and in this case, write $s(\eta)=s(\alpha_1)$ and $e(\eta)=e(\alpha_r)$. In the sequel, the length of a path $\eta$ will be written as $l(\eta)$.
A non-trivial path $\eta$ is called an {\it oriented cycle} if $s(\eta)=e(\eta)$, and an oriented cycle of length $2$ is called a $2$-{\it cycle}. Two paths $\eta, \zeta$ are called {\it parallel} if $s(\eta)=s(\zeta)$ and $e(\eta)=e(\zeta)$. A subquiver $Q'$ of $Q$ is called {\it full} if it contains all arrows $x\to y$ in $Q$ with $x, y\in Q'$, and {\it convex} if it contains all paths $x\rightsquigarrow y$ in $Q$ with $x, y\in Q'$.

For each arrow $\alpha: x\to y$ in $Q$, one introduces its {\it formal inverse} $\alpha^{-1}$ with $s(\alpha^{-1})=e(\alpha)$ and $e(\alpha^{-1})=s(\alpha)$. The formal inverses of arrows will be called  {\it inverse arrows.}
A {\it trivial walk} of length 0 is a trivial path. A \emph{walk} of length $t (\hspace{.5pt}>0)$ is a sequence $w=c_t c_{t-1}\cdots c_1$, where $c_i$ is an arrow or an inverse arrow such that $e(c_i)=s(c_{i+1})$ for $1\le i<t$. In this case, we write $s(w)=s(c_1)$ and $e(w)=e(c_t)$, and call $w$ a {\it walk from $s(w)$ to $e(w)$}. Note that $w^{-1}=c_1^{-1} \cdots c_{t-1}^{-1}c_t^{-1}$ is a walk from $y$ to $x$, and the number of arrows in $w^{-1}$ is the number of inverse arrows in $w$, while the number of inverse arrows in $w^{-1}$ is the number of arrows in $w$. A walk in $w$ in $Q$ is called \emph{reduced} if $w$ is trivial or $w=c_t c_{t-1}\dots c_1$ such that $c_{i+1}\neq c_i^{-1}$ for all $1\le i<t$. One says that $Q$ is a {\it tree} if, for any $x, y\in Q_0$, there exists at most one reduced walk from $x$ to $y$. The following notation is important for later purposes.

\begin{notation}\label{No_arrows}

Let $Q$ be a connected tree with $x, y\in Q_0$. In case $x=y$, put
$\mathfrak{a}^+(x, y)=\mathfrak{a}^-(x, y)=0$. Otherwise, there exists a unique reduced $w=c_t\cdots c_1$ from $x$ to $y$, where each $c_i$ is an arrow or an inverse arrow in $Q$. Then we denote by $\mathfrak{a}^+(x, y)$ the number of indices $i$ such that $c_i$ is an arrow; and by $\mathfrak{a}^-(x, y)$ the number of indices $j$ such that $c_j$ is an inverse arrow.

\end{notation}

We define the {\it opposite quiver} $Q^{\rm op}$ of $Q$ by $(Q^{\rm op})_0=Q_0$ and $(Q^{\rm op})_1=\{\alpha^{\rm o}: y\to x \mid \alpha: x\to y \in Q_1\}.$ The trivial paths in $Q$ are identified with the trivial paths in $Q^{\rm op}$. And a non-trivial walk $w$ in $Q$ from $x$ to $y$ induces a non-trivial walk $w^{\rm o}$ in $Q^{\rm op}$ from $x$ to $y$ so that the number of arrows in $w^{\rm o}$ is the number of inverse arrows in $w$, and the number of inverse arrows in $w^{\rm o}$ is the number of arrows in $w$.

\subsection{\sc Valued quivers} Let $\mathbb Z^+$ denote the set of positive integers. A {\it valued quiver} is a pair $(Q, v)$, where $Q=(Q_0, Q_1)$ is a quiver without multiple arrows and $v$ is the {\it valuation}, that is a map $$v: Q_1\to \mathbb{Z}^+ \times \mathbb Z^+: \alpha \mapsto v(\alpha):=(v_\alpha, v'_\alpha);$$ see, for example, \cite[Page 69]{ARS}, and compare \cite[Page 287]{HPR}. The valuation $v(\alpha)$ of an arrow $\alpha$ is called {\it trivial} if $v(\alpha)=(1,1).$ In a graphic representation of a valued quiver, we will omit the trivial valuations. In case $Q$ has no loop,  we will write $v(\alpha)=(v_{xy}, v_{xy}')$ for an arrow $\alpha: x\to y$. And in case $Q$ has no $2$-cycle, we will write $v(\alpha)=(v_{xy}, v_{yx})$ for an arrow $\alpha: x\to y$.
%
%
%
A {\it full valued subquiver} of a valued quiver $(Q, v)$ is a full subquiver of $Q$ with the valuation obtained by restricting $v$. In order to define the opposite valued quiver, we write
$(a,b)^\circ=(b,a)$ for any pair $(a,b)\in \Z^+\times \Z^+$.

\begin{defn}\label{opp_vq}

Let $(Q, v)$ be a valued quiver. Its {\it opposite valued quiver}  $(Q^{\rm op}, v^{\rm o})$ is defined by

\begin{enumerate}[$(1)$]

\item $Q^{\rm op}$ is the opposite quiver of $Q;$

\item $v^{\rm o}$ is defined by $v^\circ(\alpha^\circ)=v(\alpha)^\circ,$ for any arrow $\alpha$ in $Q$.

\end{enumerate}

\end{defn}


\subsection{\sc Valued graphs} We  denote by $\mathbb N$ the set of non-negative integers. The following definition is slightly modified from the one given in \cite[Page 241]{ARS}; compare \cite[Page 1]{DR}.

\begin{defn} A {\it valued graph} is a triplet $(\Da_0, \Da_1, v)$, where

\begin{enumerate}[$(1)$]

\item $\Da_0$ is a set of vertices$\,;$

\item $\Da_1$ is a set of edges between vertices, containing no loop or multiple edge$\,;$

\item $v$ is the {\it valuation}, that is a function
$$v: \Da_0\times \Da_0\to \mathbb N: (x, y)\mapsto v(x,y):=v_{xy}$$ such that $v_{xy} \ne 0$ if and only if $v_{yx}\ne 0$ if and only if there exists an edge between $x$ and $y$.

\end{enumerate} \end{defn}

In the sequel, a valued graph $(\Da_0, \Da_1, v)$ will be written as $(\Da, v)$, where $\Da=(\Da_0, \Da_1);$ or simply as $\Da$. If $e$ is an edge between $x$ and $y$ such that $v_{xy}=v_{yx}=1$, then we say that $e$ is {\it trivially valued} with {\it trivial valuation} $(1, 1)$. One says that $\Da$ is {\it simply laced} if every edge in $\Da$ is trivially valued. In this paper, we shall only consider valued graphs with a planar representation defined as follows.

\begin{notation}\label{v_notation}

Let $(\Da, v)$ be a planar valued graph such that every vertex is incident to at most one non-trivially valued edge. Consider an edge $e$ between vertices $x$ and $y$. If $e$ is trivially valued, then it is represented by a blank line between $x$ and $y$. Otherwise, $e$ is represented by a horizontal line $\xymatrixcolsep{12pt}\xymatrix{i\ar@{-}[r] & j\hspace{-2pt}}$ labeled with the pair $(v_{xy}, v_{yx})$ or $\hspace{-1.5pt}\xymatrixcolsep{12pt}\xymatrix{\hspace{-2pt} y \ar@{-}[r] & x\hspace{-2pt}}$ labeled with the pair $(v_{yx}, v_{xy})$. In either case, the pair is called the {\it valuation} of $e$ in the respective representation.

\end{notation}

\begin{rem}

In contrast to the valuation of an arrow in a valued quiver, the valuation of a non-trivially valued edge in a valued graph depends on how the edge is written, namely, which vertex is on the left and which is on the right.

\end{rem}

Now, we introduce valued graph isomorphisms, which we cannot find in any existing literature.

\begin{defn}

Let $(\Da, v)$ and $(\Da', v')$ be valued graphs. We define a {\it valued graph isomorphism} $\varphi: (\Da, v) \to (\Da', v')$ to be a bijection $\varphi: \Da_0 \to \Da_0'$ such that $v'_{\varphi(x) \varphi(y)}=v_{xy},$ for all $x,y\in \Da_0$.

\end{defn}

\begin{rem} Let $\varphi: (\Da, v) \to (\Da', v')$ be a valued graph isomorphism. Then, for $x,y\in \Da_0$, there exists an edge $e$ between $x$ and $y$ in $\Da$ if and only if \vspace{-2pt} there exists an edge $e'$ between $\varphi(x)$ and $\varphi(y)$ in $\Da'$. And if $e$ is written as $\hspace{-2pt}\xymatrixcolsep{14pt}\xymatrix{x\ar@{-}[r] & y}\hspace{-2pt}$ with valuation $(a,a')$, then $e'$ can be written as $\hspace{-3pt}\xymatrixcolsep{14pt}\xymatrix{\varphi(x)\ar@{-}[r] & \varphi(y)}\vspace{-4pt}\hspace{-3pt}$ with valuation $(a,a')$ or as $\hspace{-3pt}\xymatrixcolsep{14pt}\xymatrix{\varphi(y)\ar@{-}[r] & \varphi(x)}\vspace{-3pt}\hspace{-3pt}$ with valuation $(a', a)$.

\end{rem}

\vspace{2pt}


\begin{defn}

Let $(\Da, v)$ be a valued graph. The {\it weight} of $x\in \Da_0$ is defined by
${\rm w}(x)\hspace{-.5pt}:\hspace{.5pt}=\sum_{y\in \mathit{\Delta}_0} \hspace{-2pt} v_{xy}$.

\end{defn}

The following statement follows immediately from the definition a valued graph isomorphism.

\begin{lem}\label{iso-wt}

Let $\varphi : (\Da, v) \to (\Da', v')$ be a valued graph isomorphism. Given any vertex $x$ in $\Da,$ we have ${\rm w}(x)={\rm w}(\varphi(x))$.

\end{lem}

For later reference, we introduce the following definition; see \cite[Page 242]{ARS}.

\begin{defn} \label{Dyn_diag}

A valued graph is called a {\it Dynkin diagram} if it is isomorphic to one of the following {\it canonical Dynkin diagrams}$\hspace{.5pt}:$

\begin{enumerate}[$(1)$]

\vspace{4pt}

\item[] $\xymatrixcolsep{18pt}\xymatrix{\mathbb{A}_n: & 1 \ar@{-}[r] & 2 \ar@{-}[r] & \cdots \ar@{-}[r] &  n,}$ where $n\ge 1$.

\vspace{2pt}

\item[] $\xymatrixcolsep{18pt}\xymatrix{\mathbb{B}_n: & 1\ar@{-}[r]^{(1,2)} & 2 \ar@{-}[r] & 3 \ar@{-}[r] & \cdots \ar@{-}[r] & n,}$ where $n\ge 2$.

\vspace{2pt}

\item[] $\xymatrixcolsep{18pt}\xymatrix{\mathbb{C}_n: & 1\ar@{-}[r]^{(2,1)} & 2 \ar@{-}[r] & 3 \ar@{-}[r] & \cdots \ar@{-}[r] & n,}$ where $n\ge 3$.

\vspace{3pt}

\item[] $\xymatrixrowsep{14pt}\xymatrixcolsep{18pt}
\xymatrix{
&&2\ar@{-}[d]&& \\
\mathbb{D}_n: & 1\ar@{-}[r] & 3 \ar@{-}[r] &4 \ar@{-}[r]&\cdots \ar@{-}[r] & n, \mbox{ where } n\ge 4.}$

\vspace{3pt}

\item[] $\xymatrixrowsep{14pt}\xymatrixcolsep{18pt}
\xymatrix{&&&4\ar@{-}[d]&&\\
\mathbb{E}_n: &
1\ar@{-}[r] & 2 \ar@{-}[r]&3 \ar@{-}[r]&5 \ar@{-}[r] & 6 \ar@{-}[r] & \cdots \ar@{-}[r] & n, \mbox{ where } n=6,7,8.}$

\vspace{3pt}

\item[] $\xymatrixcolsep{18pt}\xymatrix{\mathbb{F}_4: & 1 \ar@{-}[r] & 2 \ar@{-}[r]^{(1,2)} & 3 \ar@{-}[r] &  4.}$

\vspace{3pt}

\item[] $\xymatrixcolsep{18pt}\xymatrix{\mathbb{G}_2: & 1 \ar@{-}[r]^{(1,3)} & 2.}$

\end{enumerate}

\end{defn}

\vspace{2pt}

It is important for us to relate the valued quivers without $2$-cycles to the valued graphs.

\begin{defn}\label{under_vq}

Let $(Q, v)$ be a valued quiver without $2$-cycles. The {\it underlying valued graph} $(\hspace{1.5pt}\overline{\hspace{-1.6pt}Q\hspace{-.5pt}}, \, \overline{\hspace{-.6pt}v})$ of $(Q, v)$ is a valued graph defined in the following way.

\begin{enumerate}[$(1)$]

\vspace{-.5pt}

\item The graph $\hspace{1.5pt}\overline{\hspace{-1.6pt}Q\hspace{-.5pt}}$ is obtained by forgetting the orientation of the arrows in $Q.$

\vspace{.0pt}

\item Given vertices $x, y$ in $\hspace{1.5pt}\overline{\hspace{-1.6pt}Q\hspace{-.5pt}}$, we define the valuation $\overline{\hspace{-.6pt}v}_{xy}$ by
$$\overline{\hspace{-.6pt}v}_{xy}=
\left\{ \hspace{-5pt} \begin{array}{ll}
0, & \mbox{ if $Q$ has no arrow $x \to y$ or $y\to x$;} \vspace{1pt} \\
v_{xy}, & \mbox{ if $Q$ contains an arrow from $x$ to $y$ 
with valuation $(v_{xy}, v_{yx});$}\\
v_{xy}, & \mbox{ if $Q$ contains an arrow from $y$ to $x$ 
with valuation $(v_{yx}, v_{xy})$.}

\end{array} \right.$$

\end{enumerate} \end{defn}

\begin{rem} \label{under_vq_val} Let $Q$ be a valued quiver without $2$-cycles. An arrow $x\to y$  with valuation $(a,b)$ in $Q$ induces an edge $\xymatrixcolsep{14pt}\xymatrix{x\ar@{-}[r] & y}\vspace{-1pt}$ with valuation $(a,b)$, or equivalently, $\xymatrixcolsep{14pt}\xymatrix{y\ar@{-}[r] & x}\vspace{-1pt}$ with valuation $(b,a)$ in $\hspace{2.5pt}\overline{\hspace{-1.6pt}Q\hspace{-.5pt}}\hspace{1pt}.$




\end{rem}

\vspace{-8pt}

\begin{eg} If $Q$ is $\xymatrixcolsep{16pt}\xymatrix{1\ar[r]^{(1,3)} & 2 & \ar[l] 3,}$ then $\hspace{2.5pt}\overline{\hspace{-1.8pt}Q\hspace{-.5pt}}\hspace{2pt}$ is
$\xymatrixcolsep{16pt}\xymatrix{1\ar@{-}[r]^-{(1,3)} & 2 \ar@{-}[r] & 3,}$ or equivalently,
$\xymatrixcolsep{16pt}\xymatrix{3 \ar@{-}[r] & 2 \ar@{-}[r]^-{(3,1)} & 1.}$





\end{eg}

The following statement is evident.

\begin{lem}

Let $(Q, v)$ be a valued quiver without $2$-cycles. Then $(Q, v)$ and $(Q^{\rm op}, v^\circ)$ have the same underlying valued graph.

\end{lem}



In the sequel, a valued quiver will be called a ({\it canonical}\hspace{.5pt}) {\it Dynkin quiver} if it admits no $2$-cycle and its underlying valued graph is a (canonical) Dynkin diagram.

\subsection{\sc Valued translation quivers} A {\it valued translation quiver} is a triplet $(\Ga, v, \tau)$, where $(\Ga, v)$ is a valued quiver and $\tau$ is a bijection, called the {\it translation}, from a subset $\Ga_0'$ of $\Ga_0$ to another $\Ga_0''$ such, for any $x\in \Ga_0'$ and any arrow $y\to x$ with valuation $(v_{yx}, v'_{yx})$, that there exists an arrow $\tau x\to y$ with valuation $(v_{yx}', v_{yx})$. In this case, we shall write $\tau^-$ for the inverse of $\tau$. The $\tau$-{\it orbit} of a vertex $x$ is the set of vertices $\tau^sx$ with $s\in \Z$.
A path $x\rightsquigarrow y$ in $\Ga$ is called \emph{sectional} if it does not contain any subpath of the form $\tau a \to b \to a$; and in this case, we call $x$ a {\it sectional predecessor} of $y$, and $y$ a {\it sectional successor} of $x$. Note that a trivial path in $\Ga$ is sectional. A sectional path is called {\it strictly sectional} if it meets any $\tau$-orbit in $\Ga$ at most once. The following easy statement follows immediately from the definition of a valued translation quiver, which will be used frequently.

\begin{lem}\label{pa_valuation}

Let $(\Ga, v, \tau)$ be a valued translation quiver. Let $x\to y$ be an arrow  in $\Ga$ with valuation $(a, a')$. If $\tau^sx$ and $\tau^sy$ are defined, then $\Ga$ contains an arrow $\tau^sx\to \tau^sy$ with valuation $(a, a')$.

\end{lem}



\subsection{\sc Krull-Schmidt categories} Let $\mathcal A$ be a Hom-finite Krull-Schmidt $R$-category. We write ${\rm rad} \mathcal{A}$ for the Jaconson radical of $\cA$, and $\rad^s \! \cA$ for the $s$-th power of $\rad \cA$ for all $s\ge 0$. A subcategory $\mathcal{B}$ of $\cA$ is called {\it convex} provided that any sequence of morphisms $X_0 \rightarrow X_1 \rightarrow \cdots \rightarrow X_{r-1} \rightarrow  X_r$
between indecomposable objects in $\mathcal A$ lies entirely in $\cB$ whenever $X_0, X_r\in \mathcal{B}.$

\vspace{2pt}

For our later purpose, we briefly recall the Auslander-Reiten theory in this general setting; see \cite{Bau} and \cite{Liu}. Given indecomposable objects $X, Y$ in $\mathcal A$, we put $k_X={\rm End}_\mathcal{A}(X)/{\rm rad}({\rm End}_\mathcal{A}(X))$, that is a division $R$-algebra; and
${\rm irr}(X, Y)={\rm rad}_\mathcal{A}(X, Y)/{\rm rad}^2_\mathcal{A}(X, Y)$, that is a $k_Y$-$k_X$-bimodule. It is well-known that a morphism $f: X\to Y$ is irreducible if and only if $f\in {\rm rad}_\mathcal{A}(X, Y) \backslash {\rm rad}^2_\mathcal{A}(X, Y)$. Write $d_{XY}={\rm dim}\, _{k_Y}{\rm irr}(X, Y)$ and $d_{XY}'={\rm dim} \, {\rm irr}(X, Y) _{k_X}$, which are related by Bautista to minimal left almost morphism and minimal right almost split morphisms as follows; see \cite[(3.3), (3.4)]{Bau}.

\begin{prop}\label{Irr_mor}

Let $\cA$ be a Hom-finite Krull-Schmidt $R$-category, and let $f:X\to Y$ be an irreducible morphism between indecomposable objects in $\cA$.

\begin{enumerate}[$(1)$]

\item If $\cA$ has a minimal left almost split morphism $g: X \to M,$ then $d_{XY}$ is the multiplicity of $Y$ as a direct summand of $M$.

\item If $\cA$ has a minimal right almost split morphism $h: N \to Y,$ then $d_{XY}'$ is the multiplicity of $X$ as a direct summand of $N$.

\end{enumerate} \end{prop}

The following definition unifies the notions of almost split sequences in abelian categories and almost split triangles in triangulated categories; see \cite[(1.3)]{Liu}, and compare \cite[(2.7)]{INP}.

\begin{defn}\label{ass_def}
Let $\cA$ be a Hom-finite Krull-Schmidt $R$-category. A sequence \hspace{-6pt} $\xymatrixcolsep{18pt}\xymatrix{X\ar[r]^-f & Y\ar[r]^-g& Z}$ \hspace{-6pt} of morphisms in $\cA$ with $Y\ne 0$ is called {\it almost split} if $f$ is minimal left almost split and a pseudo kernel of $g,$ while $g$ is minimal right almost split and a pseudo cokernel of $f.$
\end{defn}

\vspace{-2pt}

The {\it Auslander-Reiten quiver} 
of $\cA$ is a valued translation quiver $\GaA$ defined as follows; \cite[(2.1)]{Liu}. The vertex set of $\GaA$ is a complete set of representatives of isomorphism classes of indecomposable objects in $\cA$. Given vertices $X, Y$ in $\GaA$, there exists an arrow $X\to Y$ in $\GaA$ if and only if there exists an irreducible morphism $f: X\to Y$ in $\cA$; and in this case, the valuation of $X\to Y$ is $(d_{XY}, d'_{X,Y})$. The translation $\tau_{_{\hspace{-1pt}\mathcal A}}$, \vspace{-2pt} called the {\it Auslander-Reiten translation},
is defined in such a way that $X=\tau_{_{\hspace{-1pt}\mathcal A}}Z$ if and only if $\mathcal A$ has an almost split sequence \hspace{-3pt} $\xymatrixcolsep{18pt}\xymatrix{X\ar[r] & Y\ar[r] & Z.}$

\vspace{2pt}

To conclude this subsection, we recall the construction of an orbit category of $\cA$; compare \cite[(2.1)]{Ass}.

\begin{defn}\label{orbit_cat}

Let $\cA$ be a Hom-finite Krull-Schmidt $R$-category with an action by a group $G$. The {\it orbit category} $\cA/G$ is defined as follows.

\begin{enumerate}[$(1)$]

\item  The objects of $\cA/G$ are those of $\cA$.

\item Given objects $X, Y$, one has $\Hom_{\mathcal{A}/G}(X, Y)=\textstyle\oplus_{g\in G}\Hom_{\mathcal{A}}(X, g \!\cdot\! Y).$

\item Given morphisms $u=(u_g)_{g\in G}: X\to Y$ and $v=(v_g)_{g\in G}: Y\to Z$ in $\cA/G$, where $u_g: X\to g \!\cdot\! Y$ and $v_g: Y\to g \!\cdot\! Z$ are morphisms in $\cA$, one has
$vu=(w_g)_{g\in G},$ where $w_g=\sum_{fh=g} (f \!\cdot\! v_h) u_f$.

\end{enumerate}\end{defn}

The action of a group $G$ on $\cA$ is called {\it free} provided that $g\cdot X\not\cong X$ for any non-identity $g\in G$ and any indecomposable $X$ in $\cA;$ and {\it locally bounded} provided that
$\Hom_{\mathcal{A}}(X, g\cdot Y)\ne 0$ for all but finitely many $g\in G$ and all objects $X, Y$ in $\cA$; see \cite[(1.3)]{BaL}.

\begin{prop}\label{orbit_cat_prop}

Let $\cA$ be a Hom-finite Krull-Schmidt $R$-category with a free and locally bounded action by a group $G$. Then $\cA/G$ is a Hom-finite Krull-Schmidt $R$-category, whose indecomposable objects are those of $\cA$. Moreover, if $X, Y\in \cA$ are indecomposable, then

\begin{enumerate}[$(1)$]

\vspace{-.5pt}

\item $\rad^s_{\hspace{-1pt}\mathcal{A}/G}(X, Y)= \textstyle\oplus_{g\in G}\,\rad^s_{\hspace{-1pt}\mathcal{A}}(X, g \hspace{-1pt} \cdot\! Y),$ for $s\ge 0;$

\vspace{1pt}

\item $X\cong Y$ in $\cA/G$ if and only if $X\cong g\cdot Y$ in $\cA$ for some $g\in G.$

\end{enumerate} \end{prop}

\vspace{-5pt}

\noindent{\it Proof.} First of all, a direct sum in $\cA$ is a direct sum in $\cA/G$. And since the action of $G$ is locally bounded, $\cA/G$ is a Hom-finite additive $R$-category. Let $X, Y\in \cA$ be indecomposable. Consider a morphism $u=(u_g)_{g\in G}: X\to Y$ in $\cA/G$, where $u_g\in \Hom_\mathcal{A}(X, g\!\cdot\! Y)$. Suppose that $u_f\not\in \rad_\mathcal{A}(X, f\!\cdot\! Y)$ for some $f\in G$. Then, $u$ has an inverse $(u^{-1}_g)_{g\in G}: Y\to X$ in $\cA/G$, where $u^{-1}_g=0$ for all $ g \, (\ne f^{-1}) \hspace{.4pt}\in G$, and $u^{-1}_{f^{-1}}=f^{-1}\!\cdot u_f^{-1}\hspace{-2pt}: Y \to f^{-1} \!\cdot \hspace{-1pt} X$. Suppose now that $u_g\in \rad_\mathcal{A}(X, g\!\cdot\! Y)$ for all $g\in G$. Given any $v=(v_g)_{g\in G}: Y\to X$ in $\cA/G$, we have $vu=((vu)_g)_{g\in G}$, where
$(vu)_g=\sum_{fh=g} (f \cdot v_h\hspace{-.8pt})\hspace{.4pt} u_f \in \rad\hspace{-1pt}_\mathcal{A}(X, g\cdot X)$ for all $g\in G$. Let $e$ be the identity of $G$. Then, $1_{\hspace{-1pt}X}-vu=(w_g)_{g\in G}$, where $w_e=1_{\hspace{-1pt}X}-(vu)_e$, and $w_g=-(vu)_g$ for all $g \, (\ne e) \in G$. Since $w_e$ is invertible, as seen previously, so is $1_X-vu$. Hence, $u\in \rad_{\mathcal{A}/G}(X, Y)$. This shows that $\rad_{\hspace{-1pt}\mathcal{A}/G}(X, Y)= \textstyle\oplus_{g\in G}\,\rad_{\hspace{-1pt}\mathcal{A}}(X, g \!\cdot\! Y).$

\vspace{.5pt}

Since $G$ acts freely on $\cA$, $\rad_{\mathcal{A}/G}(X,X)=\rad_\mathcal{A}(X, X) \oplus (\oplus_{g \, (\ne e) \in G} \Hom_\mathcal{A}(X, g\cdot \hspace{-1pt} X)).$ So, $\End_{\mathcal{A}/G}(X)$ is local. Since a non-zero object in $\cA$ is non-zero in $\cA/G$, an object is indecomposable in $\cA$ if and only if it is indecomposable in $\cA/G$, and consequently, $\cA/G$ is Krull-Schmidt. Now, by induction, we easily see that Statement (1) for all $s\ge 0$. Finally, $X\cong Y$ in $\cA/G$ if and only if $\Hom_{\hspace{-1pt}\mathcal{A}/G}(X, Y)\ne \rad_{\hspace{-1pt}\mathcal{A}/G}(X, Y),$ if and only if $\Hom_\mathcal{A}(X, g\cdot Y)\ne \rad_{\hspace{-1pt}\mathcal{A}}(X, g \!\cdot\! Y)$ for some $g\in G$, if and only if $X\cong g\cdot Y$ in $\cA$ for some $g\in G.$ The proof of the proposition is completed.

\subsection{\sc Module category} Let $\La$ be an artin algebra. We denote by $\mmod \La$ the category of finitely generated left $\La$-modules, and by $\rad(\mmod \La)$ the Jacobson radical of $\mmod \La$. The maps in $\rad(\mmod \La)$ are called {\it radical maps}. 
Given a module $M$ in $\mmod \La$, we shall write $\rad M$, ${\rm top}M$ and ${\rm soc}M$ for the radical, the top and the socle of $M$ respectively.
The Auslander-Reiten quiver $\Ga_{\hspace{-.8pt}{\rm mod}\mathit\Lambda}$ of $\mmod \La$ carries the essential information of the finite powers of $\rad(\mmod \La)$. The Auslander-Reiten translations $\tau_{\hspace{-1pt}_{\mathit\Lambda}}$ and $\tau^-_{{\hspace{-1pt}_\mathit{\Lambda}}}$ of $\Ga_{\hspace{-.8pt}{\rm mod}\mathit\Lambda}$ are given by $D{\rm Tr}$ and ${\rm Tr}D$, respectively; see \cite[Page 22]{ARS}. Besides the Auslander-Reiten quiver, the ext-quiver of $\La$; see \cite[Page 69]{ARS} also plays an important role in our study.

\begin{defn} \label{def_Ext_quiver} Let $\La$ be an artin algebra, and let $S_1,...\,,S_n$ be the non-isomorphisc simple modules in $\mmod \La$. The {\it ext-quiver} of $\La$ is a valued quiver $Q_{\hspace{-1pt}\mathit\Lambda}$ defined in the following way.
\begin{enumerate}[$(1)$]

\item The vertex set is $(Q_{\hspace{-1pt}\mathit\Lambda})_0=\{1,...\,,n\}.$

\item The arrow set is $(Q_{\hspace{-1pt}\mathit\Lambda})_1=\{i\xrightarrow{}j \ |\  \Ext^1_{\hspace{-1pt}\mathit\Lambda}(S_i,S_j)\neq 0\}.$

\item The valuation for an arrow $i\xrightarrow{}j$ in $Q_{\hspace{-1pt}\mathit\Lambda}$ is the pair $(d_{ij}, d_{ij}')$, \vspace{-1pt} where $d_{ij}=\dim_{\hspace{.4pt} \End_{\mathit\Lambda}(S_j)}\!\Ext^1_{\hspace{-1pt}\mathit\Lambda}(S_i,S_j)$ and $d_{ij}'= \dim\Ext^1_{\hspace{-1pt}\mathit\Lambda}(S_i,S_j)_{\End_{\hspace{-1pt}\mathit\Lambda}(S_i)}$.

\end{enumerate}
\end{defn}




In the study of the Auslander-Reiten quiver $\Ga_{\hspace{-.8pt}{\rm mod}\mathit\Lambda}$, sectional paths play an important role; see \cite{BSm, Liu3, Liu4, IgT}. The following statement tells us 
when we can extend a sectional path  in $\Ga_{\hspace{-.8pt}{\rm mod}\mathit\Lambda}$.

\begin{prop}\label{p_sectional_path}
Let $\La$ be an artin algebra, and let  $X_0\to X_1\to \cdots \to X_{s-1}\to X_s$  be a non-trivial sectional path in the Auslande-Reiten quiver $\Ga_{\hspace{-.8pt}{\rm mod}\mathit\Lambda}$.

\begin{enumerate}[$(1)$]

\vspace{-1pt}

\item Suppose that $X_0$ is projective. If $Y\to X_s$ is an arrow in $\Ga_{\hspace{-.8pt}{\rm mod}\mathit\Lambda}$ with $Y\ne X_{s-1}$, then $Y$ is not injective, and consequently, $\Ga_{\hspace{-.8pt}{\rm mod}\mathit\Lambda}$ contains a sectional path $X_0\to 
\cdots \to X_{s-1}\to X_s \to \tau^-_{\hspace{-1pt}_{\mathit\Lambda}\hspace{-1pt}}Y.$

\vspace{.5pt}

\item Suppose that $X_s$ is injective. If $X_0\to Y$ is an arrow in $\Ga_{\hspace{-.8pt}{\rm mod}\mathit\Lambda}$ with $Y\ne X_1$, then $Y$ is not projective, and consequently, $\Ga_{\hspace{-.8pt}{\rm mod}\mathit\Lambda}$ contains a sectional path $\tau_{\hspace{-1pt}_{\mathit\Lambda}\hspace{-.5pt}}Y \to X_0\to X_1\to \cdots 
\to X_s.$

\end{enumerate}\end{prop}

\vspace{-4pt}

\noindent{\it Proof.} We shall only prove Statement (1). Let $Y\to X_s$ be an arrow in $\Ga_{\hspace{-.8pt}{\rm mod}\mathit\Lambda}$ with $Y\ne X_{s-1}$. Suppose that $Y$ is injective. Then, we have an irreducible epimorphism $f_s: Y\to X_s$. 
Since $Y\ne X_{s-1}$, there exists an irreducible map $(f_s, g): Y \oplus X_{s-1}\to X_s$; see \cite[(3.2)]{Bau}. So, $\mmod \La$ has an almost split sequence \vspace{-8pt}
$$\xymatrixcolsep{35pt}\xymatrix{0\ar[r] & \tau_{\hspace{-1pt}_{\mathit\Lambda}\hspace{-.5pt}} X_s\ar[r]^-{(h, u , f_{s-1})} & Y \oplus Z \oplus X_{s-1} \ar[r]^-{\begin{pmatrix} f_s \\ w \\ g\end{pmatrix}} & X_s\ar[r] & 0.}$$

Since $f_s$ is an irreducible epimorphism, so is $f_{s-1}: \tau_{\hspace{-1pt}_{\mathit\Lambda}\hspace{-.5pt}} X_s\to X_{s-1}$.
If $s>1$, then $\tau_{\hspace{-1pt}_{\mathit\Lambda}\hspace{-.5pt}} X_s\ne X_{s-2}$, and we similarly obtain an irreducible epimorphism $f_{s-2}: \tau_{\hspace{-1pt}_{\mathit\Lambda}\hspace{-.5pt}} X_{s-1}\to X_{s-2}$. Continuing this process, we obtain an irreducible epimorphism $f_1: \tau_{\hspace{-1pt}_{\mathit\Lambda}\hspace{-.5pt}} X_1\to X_0$, absurd. The proof of the proposition is completed.

\vspace{2pt}

It is well-known that $A$ is representation-finite if and only if $\rad(\mmod \La)$ is nilpotent; see
\cite[(V.7.7)]{ARS} and \cite[(1.1)]{Sko}. To compute the nilpotency of $\rad(\mmod \La)$, one introduces the notion of depth for maps 
in terms of the radical series of $\mmod \La$. 

\begin{defn}\label{dep_def}

Let $\La$ be an artin algebra. The \emph{depth} of a map $f: M\to N$ in $\mmod \La$ is defined by
$$\ndp(f):=\sup\{ s \, \in \mathbb{N} \ |\  f\in \rad^s\hspace{-1pt}(M,N)\}.$$


\end{defn}

\begin{rem} \label{nilp_dpt} If $M, N\in \mmod \La$ are indecomposable, then a map $f: M\to N$ is irreducible if and only if $\ndp(f)=1$.

\end{rem}


\vspace{1pt}

The following statement relates maps of finite depth in $\mmod \La$ to paths in $\Ga_{\mmod \mathit\Lambda}$. 

\begin{lem}\label{t_irreducibles}

Let $\La$ be an artin algebra. Consider a radical map $f: M\to N$ between indecomposable modules in $\mmod \La$. If $\ndp(f)=t$, then there exists a chain of irreducible maps between indecomposable modules $M\xrightarrow{f_1} M_1 \to \cdots \to M_{t-1}\xrightarrow{f_t} N$ in $\mmod \La$ such that $\ndp(f_t\cdots f_1)=t$.

\end{lem}

\vspace{-4pt}

\noindent {\it Proof.} \vspace{.5pt} Assume that ${\rm dp}(f)=t>0$. Since $f\in \rad^t\hspace{-1pt}(M, N)$, we may write $f= \sum_{i=1}^p f_{it}\cdots f_{i1},$ where the $f_{ij}$ are radical maps between indecomposable modules in $\mmod \La$. Since $f\notin \rad^{t+1}(M, N)$, there exists some $1\le s\le p$ such that $f_{st}\cdots f_{s1}\notin \rad^{t+1}\hspace{-1pt}(M, N)$. Thus,
$f_{s1}, \ldots, f_{st}$ are irreducible such that $\ndp(f_{st}\cdots f_{s1})=t$. The proof of the lemma is completed.

\subsection{\sc Hammocks} Let $\La$ be an artin algebra with $S$ a simple module in $\mmod \La$. For any module $M$ in $\mmod A$, we write $\ell_S(M)$ for the multiplicity of $S$ as a composition factor of $M$. For convenience of reference, we state the following well-known  statement; see, for example, \cite[Page 45]{ARS}.

\begin{lem}\label{m_com_fac}

Let $\La$ be an artin algebra. Consider a simple module in $\mmod \La$ with  projective cover $P$ and injective envelope $I$. If $M\in \mmod \La$, then $\ell_S(M)$ is equal to the length of the right ${\rm End}_{\mathit\Lambda}(P)$-module $\Hom_{\mathit\Lambda}(P, M)$, as well as the length of the left ${\rm End}_{\mathit\Lambda}(I)$-module $\Hom_{\mathit\Lambda}(M, I)$.

\end{lem}

Recall that $\La$ is of {\it directed representation type} if $\Ga_{\hspace{-.8pt}{\rm mod}\mathit\Lambda}$ is finite and contains no oriented cycle. The following definition is due to Brenner; see \cite{Bre}.

\begin{defn}\label{hammock_defn}

Let $\La$ be an artin algebra of directed representation type. Consider a simple module $S$ in the Auslander-Reiten quiver $\Ga_{\hspace{-.8pt}{\rm mod}\mathit\Lambda}$.

\begin{enumerate}[$(1)$]

\item The {\it hammock} associated with $S$ is the full valued translation subquiver $\mathcal{H}_S$ of $\Ga_{\hspace{-.8pt}{\rm mod}\mathit\Lambda}$ generated by the modules $M$ of which $S$ is a composition factor.

\item The {\it canonical hammock function} associated with $S$ is the function
$h_S: \mathcal{H}_S\to \mathbb{N}: M\mapsto \ell_S(M).$

\end{enumerate}

\end{defn}

The following statement explains in particular the terminology of a hammock.

\begin{prop}\label{hammock_prop}

Let $\La$ be an artin algebra of directed representation type. Consider a simple module $S$ with projective cover $P$ and injective envelope $I$, where $S, P, I\in \Ga_{\hspace{-.8pt}{\rm mod}\mathit\Lambda}$.

\begin{enumerate}[$(1)$]

\vspace{-1pt}

\item Every module in the hammock $\mathcal{H}_S$ is a successor of $P$ and a predecessor of $I$.

\item If $M$ is a sectional successor of $P$ or a sectional predecessor of $I$ in $\Ga_{\hspace{-.8pt}{\rm mod}\mathit\Lambda}$, then $M\in \mathcal{H}_S.$

\end{enumerate} \end{prop}

\vspace{-4pt}

\noindent{\it Proof.} It is clear that $P, I\in \mathcal{H}_S$. Consider a module $M\in \mathcal{H}_S$. By Lemma \ref{m_com_fac}, we have non-zero maps $f: P\to M$ and $g: M\to I$.
Since $\La$ is representation-finite, $\Ga_{\hspace{-.8pt}{\rm mod}\mathit\Lambda}$ contains paths $P\rightsquigarrow M$ and $M\rightsquigarrow
I$; see \cite[(V.7.8)]{ARS}. We claim that $\mathcal{H}_S$ contains a path $P\rightsquigarrow M$. Since $\GaH$ is finite without oriented cycles, we have a maximal length $l_M$ of paths $P\rightsquigarrow M$ in $\GaH$. We may assume that $l_M>0$. Then, $f$ factors through the minimal right almost split map for $M$. Thus,  $\mathcal{H}_S$ contains an arrow $N\to M$. Since $l_N<l_M$, by the induction hypothesis, $\mathcal{H}_S$ contains a path $P\rightsquigarrow N$, and hence, a path $P \rightsquigarrow M$. Our claim holds. Dually, $\mathcal{H}_S$ contains a path $M \rightsquigarrow I$. This establishes Statement (1). Next, if $M$ is sectional successor of $P$ or a sectional predecessor of $I$ in $\Ga_{\hspace{-.8pt}{\rm mod}\mathit\Lambda}$, then $\Hom_{\mathit\Lambda}(P, M)\ne 0$ or $\Hom_{\mathit\Lambda}(M, I)\ne 0$; see \cite[(13.4)]{IgT} and also \cite[(VII.2.4)]{ARS}. So,
$M\in \mathcal{H}_S$; see (\ref{m_com_fac}). The proof of the proposition is completed.

\section{Valued translation quivers with sections}

The notion of sections is useful in describing Auslander-Reiten components without oriented cycles for artin algebras; see \cite{Liu4, Liu2}. The objective of this section is to further investigate valued translation quivers with sections, by studying the repetitive quiver $\Z\Da$ of a valued quiver $\Da$ without oriented cycles. The sectional paths in $\Z\Da$ are precisely described in terms of the reduced walks in $\Da$; see (\ref{reduced_sectional}) and (\ref{sec_paths}). In case $\Da$ is a tree, a sectional path in $\Z\Da$ is the only path between its end-points, and any two parallel paths have the same length; see (\ref{uni_sec_path}). Finally, we shall consider additive functions on stable valued translation quivers with sections. These results will be applied to preprojective and preinjective components of hereditary artin algebras in the next two sections.

\subsection{\sc Repetitive quiver} A valued translation quiver $(\Ga, v, \tau)$ is called {\it stable} if $\tau x$ and $\tau^{-}x$ are defined for all $x\in\Ga_0$. In her classification of self-injective algebras of finite representation type, Riedtmann introduced a canonical way to construct a stable translation quiver by ``knitting" repetitively a tree; see \cite{Rie}. This has been generalized to ``knit" any valued quivers without oriented cycles.

\begin{defn} \label{repetivie_quivers}
Let $(\Da, v)$ be a valued quiver without oriented cycles. The {\it repetitive quiver} $\mathbb Z\Da$ of $\Da$ is a stable valued translation quiver defined as follows:

\begin{enumerate}[$(1)$]

\item The vertex set is
$(\mathbb Z\Da)_0=\{(s,x)\ |\ s\in \mathbb Z; \ x\in \Da_0 \}.$

\vspace{.5pt}

\item The arrow set is $(\mathbb Z\Da)_1= \{(s,x)\xrightarrow{(s,\alpha)} (s,y) \mbox{ and } (s,y)\xrightarrow{(s,\alpha)^{*}} (s+1,x) \mid s\in \mathbb{Z};\, \alpha: x\to y\in \Da_1\}.$

\vspace{1.5pt}

\item The translation
$\tau: (\mathbb Z\Da)_0 \to (\mathbb Z\Da)_0$ is defined by $\tau(s,x)=(s-1,x)$, for all $(s, x)\in \mathbb Z\times \Da_0$.

\vspace{1.5pt}

\item The valuation $\tilde{v}$ is such that $\tilde{v}((s, \alpha))=v(\alpha)$ and $\tilde{v}((s, \alpha)^*)=v(\alpha)^\circ$, for all $\alpha\in\Da_1$ and $s\in \Z$.

\end{enumerate}\end{defn}

\begin{rem}

Since $\Da$ contains no oriented cycles, it is evident that $\mathbb{Z} \Da$ contains no oriented cycles.

\end{rem}


We denote by $\mathcal{P}(\mathbb Z\Da)$ the set of paths in $\Z\Da$, and by $\mathcal{W}(\Da)$ the set of walks in $\Da$. 

\begin{prop} \label{covering}

Let $\Z\Da$ be  the repetitive quiver of a valued quiver $\Da$ without oriented cycles. Setting $\pi(s,\alpha)=\alpha$ and $\pi(s,\alpha)^{*}=\alpha^{-1},$ we obtain a length-preserving map $\pi: \mathcal{P}(\mathbb Z\Da) \to \mathcal{W}(\Da)$ such that $\pi(\eta \cdot \xi)= \pi(\eta)\cdot \pi(\xi)$, for all $\eta, \xi\in \mathcal{P}(\mathbb Z\Da).$

\end{prop}

\vspace{-4pt}

\noindent{\it Proof.} Given a trivial path $\varepsilon_{(s, x)}$ with $(s, x) \in (\mathbb{Z}\Da)_0$, we set $\pi(\varepsilon_{(s, x)})=x$. Consider a path $\eta=\delta \gamma$, where $\gamma, \delta$ are arrows in $\mathbb{Z}\Da$. Assume that $\gamma=(s, \alpha)$ for some $s\in \mathbb Z$ and $\alpha: x\to y$ in $\Da_1$. Then $\delta=(s, \alpha)^{*}$ or $\delta=(s, \beta)$ for some $\beta: y\to z$ in $\Da_1$. This yields a walk $\pi(\delta) \cdot \pi(\gamma)=\alpha^{-1} \alpha$  in the first case and
a path $\pi(\delta) \cdot \pi(\gamma)=\beta \alpha$ in the second case. In any case, we set $\pi(\eta)=\pi(\delta) \cdot \pi(\gamma) \in \mathcal{W}(\Da)$. Similarly in case $\gamma=(s, \alpha)^*$, we set $\pi(\eta)=\pi(\delta) \cdot \pi(\gamma)\in \mathcal{W}(\Da).$ By induction, $\pi$ extends to a map $\pi: \mathcal{P}(\Z\Da)\to \mathcal{W}(\Da)$. It is easy to see that $\pi(\eta \cdot \xi)= \pi(\eta)\cdot \pi(\xi)$, for all $\eta, \xi\in \mathcal{P}(\mathbb Z\Da)$. The proof of the proposition is completed.

\vspace{3pt}

In the sequel, the map $\pi: \mathcal{P}(\mathbb Z\Da) \to \mathcal{W}(\Da)$ defined above will be called the {\it covering map,} which can be used to determine sectional paths in $\Z\Da$ as shown below.

\begin{lem} \label{reduced_sectional}

Let $\Z\Da$ be the repetitive quiver of a valued quiver $\Da$ without oriented cycles. Then a path $\eta$ in $\mathbb Z\Da$ is sectional if and only if $\pi(\eta)$ is a reduced walk in $\Da$, where $\pi$ is the covering map.

\end{lem}

\noindent {\it Proof.} Let $\eta\in \mathcal{P}(\Z\Da)$. If $\eta$ is not sectional, then it contains a subpath $(s-1, \alpha)^{*}  (s, \alpha)$ or $(s,\alpha)(s, \alpha)^{*}$, for some $s\in \mathbb Z$ and $\alpha \in \Da_1$. This yields a subwalk $\alpha^{-1} \alpha$ or $\alpha \alpha^{-1}$ of $\pi(\eta)$. So, $\pi(\eta)$ is not reduced. If $\pi(\eta)$ is not reduced, then we similarly show that $\eta$ is not sectional. The proof of the lemma is completed.

\vspace{3pt}

The following statement says that $\pi$ is indeed a covering map, which we do not rigorously define.

\begin{lem} \label{eta_eta}

Let $\Z\Da$ be the repetitive quiver of a valued quiver $\Da$ without oriented cycles, with paths $\eta, \xi$ such that $\pi(\eta)\hspace{-1pt}=\hspace{-1pt}\pi(\xi)$, where $\pi$ is the covering map. If $s(\eta)\hspace{-1pt}=\hspace{-1pt}s(\xi)$ or $e(\eta)\hspace{-1pt}=\hspace{-1pt}e(\xi)$, then $\eta=\xi$.

\end{lem}

\vspace{-4pt}

\noindent{\it Proof.} Consider only the case where $s(\eta)=s(\xi)=(t,x)=:\tilde x$, for some $t\in \mathbb Z$ and $x\in \Da_0$. Since $\pi$ preserves the length, it is easy to see that the statement holds if $l(\eta)=0$ or $1$. Suppose that $l(\eta)\ge 2$. Write $\eta=\eta' \tilde\gamma$ and $\xi=\xi' \tilde\delta$, for some arrows $\tilde\gamma: \tilde x\to \tilde y$ and $\tilde\delta: \tilde x\to \tilde z$, and some paths $\eta', \xi'$ such that $s(\eta')=\tilde y$ and $s(\xi')=\tilde z$. Then,
$\pi(\eta')\cdot \pi(\tilde\gamma)= \pi(\eta)=\pi(\xi)=\pi(\xi')\cdot \pi(\tilde\delta)$. Since $\pi(\tilde\gamma)$ and $\pi(\tilde\delta)$ are arrows or inverse arrows in $\Da$, we have $\pi(\tilde\gamma)=\pi(\tilde\delta)$ and $\pi(\eta') = \pi(\xi')$. As seen above, $\tilde\gamma=\tilde\delta$. Then $s(\eta')=s(\xi')$. By the induction hypothesis, $\eta' = \xi'$. Hence, $\eta=\xi$. The proof of the lemma is completed.

\vspace{3pt}

We precisely describe the sectional paths in $\Z \Da$ in terms of reduced walks in $\Da$ as follows.

\begin{lem} \label{sec_paths}

Let $\Da$ be a valued quiver without oriented cycles, and let $w=c_t\cdots c_2c_1$ be a reduced walk in $\Da$, where $c_i$ is an arrow or an inverse arrow from $x_{i-1}$ to $x_i$, such that the edge $\hspace{-2pt}\xymatrixcolsep{13pt}\xymatrix{
x_{i-1} \ar@{-}[r] & x_i}\hspace{-3pt}$ \vspace{-2pt} in $\hspace{2pt}\overline{\hspace{-2pt}\Da}$ has valuation $(v_i, v_i')$, for $i=1, \ldots, t$. Then, $\mathbb Z\Da$ contains a sectional path \vspace{-2pt}
$$\xymatrixcolsep{22pt}\xymatrix{\hspace{-2pt}  (0, x_0) \ar[r] & (r_1, x_1) \ar[r] & \cdots \ar[r] & (r_t, x_t),} \vspace{-2pt}$$
where $r_i$ is the number of inverse arrows in the subwalk of $w$ from $x_0$ to $x_i$, and $(r_{i-1}, x_{i-1})\to (r_i, x_i)$ has valuation $(v_i, v_i')$, for $i=1, \ldots, t.$

\end{lem}

\vspace{-4pt}

\noindent {\it Proof.} Observe that $r_0=0$. We may assume that $t>0$ and that $\mathbb Z\Da$ contains a desired sectional path $(x_0, r_0)\to \cdots \to (r_{t-1}, x_{t-1})$.
Suppose first that $\Da$ has an arrow $x_{t-1}\to x_t$. Then, its valuation is $(v_t, v_t')$; see (\ref{under_vq_val}) and $r_t=r_{t-1}$. By definition, $\Z \Da$ has an arrow $(r_{t-1}, x_{t-1})\to (r_t, x_t)$ with valuation $(v_t, v_t')$. Suppose now that $\Da$ has an arrow $x_{t-1}\leftarrow x_t$. Then, its valuation is $(v_t', v_t)$; see (\ref{under_vq_val}) and $r_t=r_{t-1}+1$. Thus, $\Z\Da$ has an arrow $(r_{t-1}, x_{t-1})\to (r_t, x_t)$ with valuation $(v_t, v_t')$. Since $x_t\ne x_{t-2}$, we obtain a desired sectional path as stated in the lemma. The proof of the lemma is completed.

\vspace{1pt}

The following easy statement is useful for inductive arguments.

\begin{lem}\label{non-sectional}

Let $(\Ga, v, \tau)$ be a stable valued translation quiver. If $\Ga$ contains a non-sectional path
$\zeta: x\rightsquigarrow y$, then it contains a path $\xi:  x \rightsquigarrow \tau y$ with $l(\xi)=l(\zeta)-2$.

\end{lem}

Now, we state some nice properties of repetitive quivers of valued trees.

\begin{prop}\label{uni_sec_path}

Let $\mathbb Z\Da$ be the repetitive quiver of a valued tree $\Da$ with translation $\tau$.

\begin{enumerate}[$(1)$]

\item If $\eta: \tilde x\rightsquigarrow \tilde y$ is a sectional path in $\mathbb Z\Da$, then it is the only path from $\tilde x$ to $\tilde y$. 

\item Any two parallel paths $\eta: \tilde x\rightsquigarrow \tilde y$ and $\zeta: \tilde x\rightsquigarrow \tilde y$ in $\mathbb Z\Da$ have the same length.

\end{enumerate} \end{prop}

\vspace{-4pt}

\noindent{\it Proof.}  Consider the covering map $\pi: \mathcal{P}(\Z\Da)\to \mathcal{W}(\Da)$. Let $\eta: \tilde x\rightsquigarrow \tilde y$ and $\zeta: \tilde x\rightsquigarrow \tilde y$ be paths in $\mathbb Z\Da$.

(1) If $\ell(\eta)=0$, then clearly $\xi=\eta$. Assume that $\ell(\eta)>0$. Write $\eta= \tilde \alpha \eta'$, where $\tilde\alpha: \tilde z\to \tilde y$ is an arrow and $\eta': \tilde x\rightsquigarrow \tilde z$ is a sectional path with $l(\eta')=l(\eta)-1.$ Consider the arrow $\tilde \beta: \tau \tilde y\to \tilde z$ in $\Z \Da$. If $\xi$ is not sectional then, by Lemma \ref{non-sectional}, $\mathbb{Z}\Da$ contains a path $\xi': \tilde x\rightsquigarrow \tau \tilde y$, and a path $\tilde\beta \xi': \tilde x\rightsquigarrow \tau \tilde y\to \tilde z$. By the induction hypothesis, $\eta'=\tilde \beta \xi'$, and consequently, $\eta=\tilde \alpha \tilde \beta \xi'$, a contradiction to $\eta$ being sectional. Thus, $\xi$ is sectional. By Lemma \ref{reduced_sectional}, both $\pi(\eta)$ and $\pi(\xi)$ are reduced walks in $\Da$ from $\pi(\tilde x)$ to $\pi(\tilde y)$. Since $\Da$ is tree, $\pi(\eta) = \pi(\xi)$, and by Proposition \ref{eta_eta}, $\eta=\xi$.

(2) By Statement (1), we may assume that $\eta$ and $\eta'$ are non-sectional. By Lemma \ref{non-sectional}, there exist paths $\xi:  \tilde x \rightsquigarrow \tau \tilde y$ and $\xi':  \tilde x \rightsquigarrow \tau \tilde y$, where $\tau$ is the translation of $\mathbb Z\Da$, such that $l(\xi)=l(\eta)-2$ and $l(\xi')=l(\eta')-2$. By the induction hypothesis, $l(\xi)=l(\xi')$. Therefore, $l(\eta)=l(\eta')$. The proof of the lemma is completed.

\subsection{\sc Sections} In the existing literature, sections are defined only for connected valued translation quivers; see \cite[(2.1)]{Liu2}. We shall drop this restriction in order to include the Auslander-Reiten quiver of the bounded derived category of a simple artin algebra; see (\ref{Der_ARQ}). Throughout this subsection, $\Ga$ denotes a valued translation quiver with translation $\tau$.

\begin{defn}\label{Section}

Let $\Ga$ be a valued translation quiver with translation $\tau$. A valued subquiver $\Da$ of $\Ga$ is called a {\it section} if the following conditions are satisfied:

\begin{enumerate}[$(1)$]

\item $\Da$ is connected and convex in $\Ga$.

\item $\Da$ contains no oriented cycle.

\item $\Da$ meets every $\tau$-orbit in $\Ga$ exactly once.

\end{enumerate}\end{defn}

The following statement says that $\mathbb{Z}\mathbb{A}_1$ is essentially the only valued translation quiver which is not connected and contains sections.

\begin{lem}\label{conn_sec}

Let $\Ga$ be a valued translation quiver with translation $\tau$, containing a non-trivial section $\Da$. Then,
$x^+$ and $x^-\hspace{-3pt}$ are non-empty for all $x\in \Ga_0$. And consequently, $\Ga$ is connected.

\end{lem}

\vspace{-4pt}

\noindent{\it Proof.} Let $x\in \Ga_0$. Then $x=\tau^s\hspace{-.5pt}y$ for some $y\in \Da_0$ and $s\in \Z$. Since $\Da$ is non-trivial and connected, we may assume that it contains an arrow $y\to z$. Being stable, $\Ga$ contains arrows $\tau^s\hspace{-.5pt}y \to \tau^s\hspace{-.5pt}z$ and $\tau^{s+1}\hspace{-.5pt}z \to \tau^s\hspace{-.5pt}y$. Thus,  $x^+$ and $x^-\hspace{-3pt}$ are non-empty. The proof of the lemma is completed.

\vspace{3pt}

Let $\Da$ be a section in $\Ga$. By definition, every vertex in $\Ga$ is uniquely written as $\tau^{-r}x$ with $r\in \Z$ and $x\in \Da_0.$ We shall say that $\Da$ is a {\it leftmost section} if all vertices in $\Ga$ are of the form $\tau^{-r}x$ with $r\ge 0$ and $x\in \Da_0;$ and a {\it rightmost section} if all vertices in $\Ga$ are of the form $\tau^{r}x$ with $r\ge 0$ and $x\in \Da_0$.

\vspace{-1pt}

\begin{prop}\label{sec_property}

Let $\Ga$ be a valued translation quiver with translation $\tau$, containing a section $\Da$.

\begin{enumerate}[$(1)$]

\vspace{-1pt}

\item If $x\to y$ is an arrow in $\Ga$, then $y\in \Da$ implies $x$ or $\tau^-x\in \Da;$ and $x\in \Da$ implies $y$ or $\tau y\in \Da$.

\item There exists a valued translation quiver embedding $\phi: \Ga\to \Z\Da$,
sending $\tau^sx$ to $(-s,x)$.

\item If $\Da$ is leftmost or rightmost, then it is the unique leftmost or rightmost section in $\Ga$, respectively.

\end{enumerate} \end{prop}

\vspace{-4pt}

\noindent {\it Proof.}  Statements (1) and (2) are quoted from \cite[(2.2), (2.3)]{Liu2}, and Statement (2) follows immediately from the definitions. 
The proof of the proposition is completed.



\vspace{3pt}

In view of Proposition \ref{sec_property}(2), we will see that valued translation quivers with sections inherits many nice properties of repetitive quivers. We will also need the following notion.

\begin{defn}\label{S-section}

Let $\Ga$ be a valued translation quiver with a vertex $x$. A section $\Da$ in $\Ga$ is called an {\it $x$-source section} if every vertex in $\Da$ is a successor of $x$ in $\Da;$ and an {\it $x$-sink section} if every vertex in $\Da$ is a predecessor of $x$ in $\Da.$

\end{defn}

Given a vertex $x$ in $\Ga$, we denote by
${\rm Suc}(x)$ and ${\rm Pred}(x)$ the full valued translation subquivers of $\Ga$ generated by the successors of $x$ and by the predecessors of $x$, respectively.

\begin{lem}\label{s_sec_unicity}

Let $\Ga$ be a valued translation quiver with a vertex $x$.

\begin{enumerate}[$(1)$]

\vspace{-1.5pt}

\item If $\Da$ is an $x$-source section in $\Ga$, then it is unique in $\Ga$ and a leftmost section in ${\rm Suc}(x)$.

\item If $\Sa$ is an $x$-sink section in $\Ga$, then it is unique in $\Ga$ and a rightmost section in ${\rm Pred}(x)$.

\end{enumerate}  \end{lem}

\vspace{-4pt}

\noindent{\it Proof.} We shall only prove Statement (1). Let $\Da$ be an $x$-source section in $\Ga$. Then, $\Da\subseteq {\rm Suc}(x)$. Thus, $\Da$ is a section in ${\rm Suc}(x)$. Suppose that $\Ga$ contains a path $x\rightsquigarrow y$. Then, $y=\tau^r\hspace{-1pt}z$ for some $r\in \Z$ and $z\in \Da_0$. If $r<0$, then $\Ga$ contains a path $x\rightsquigarrow y \rightsquigarrow z$. Since $\Da$ is convex in $\Ga$, both $y$ and $\tau^ry$ belong to $\Da$, a contradiction. Thus, $r\ge 0$. That is, $\Da$ is a leftmost section in ${\rm Suc}(x)$, which is unique by Lemma \ref{sec_property}. Thus, $\Da$ is the unique $x$-source section in $\Ga$. The proof of the lemma is completed.

\vspace{3pt}

The following statement is important for our investigation in the next section.

\begin{prop}\label{ZD_sec}

Let $\Z\Da$ be the repetitive quiver of a connected valued tree $\Da$ with a vertex $(r,x)$. Then the sectional successors of $(r,x)$ in $\Z\Da$ generate an $(r,x)$-source section  in $\Z\Da;$ and the sectional predecessors of $(r,x)$ in $\Z\Da$ generate an $(r,x)$-sink section in $\Z\Da$.

\end{prop}

\vspace{-4pt}

\noindent{\it Proof.} We shall only prove the first part of the statement. Let $\Oa$ be the full valued subquiver of $\Z\Da$ generated by the sectional successors of $(r,x)$ in $\Z\Da$. It is easy to see that every vertex in $\Oa$ is a successor of $(r, x)$ in $\Oa$. In particular, $\Oa$ is connected and contains no oriented cycle. Now, suppose that $\Z\Da$ contains a path $(s,y) \rightsquigarrow (s', y') \rightsquigarrow (t,z)$, where $(s, y), (t, z)\in \Oa_0$. Then, $\Z\Da$ contains sectional paths $\zeta: (r,x) \rightsquigarrow (s, y)$ and $\xi: (r,x) \rightsquigarrow (t, z)$. By Lemma \ref{uni_sec_path},
$(r,x) \rightsquigarrow (s, y) \rightsquigarrow (s', y') \rightsquigarrow (t,z)$ coincides with $\xi$. Thus, $(r,x) \rightsquigarrow (s, y) \rightsquigarrow (s', y')$ is sectional, and $(s', y')\in \Oa$. So, $\Oa$ is convex in $\Z\Da$.

Consider finally any $y\in \Da_0$. Since $\Da$ is connected, by Lemma \ref{sec_paths}, $\Z\Da$ has a sectional path $(r,x) \rightsquigarrow (s_0, y)$ with $s_0\in \Z$. So, $(s_0, y)\in \Oa$. Suppose that $(t_0, y)\in \Z\Da$ for another $t_0\in \Z$, say $s_0\le t_0$. Then, $\Z\Da$ contains a path $(s_0, y) \rightsquigarrow (t_0,y)$ and a sectional path $(r,x) \rightsquigarrow (t_0, y)$. By Lemma \ref{uni_sec_path}, the path $(r,x) \rightsquigarrow (s_0, y) \rightsquigarrow (t_0,y)$ is sectional. So, $s_0=t_0.$ Therefore, $\Oa$ meets every $\tau$-orbit in $\Z\Da$ exactly once. Thus, $\Oa$ is a $(r,x)$-source section in $\Z\Da$. The proof of the proposition is completed.

\subsection{\sc Additive functions} In the study of Auslander-Reiten components of artin algebras, a powerful tool is additive functions on translation quivers; see, for example, \cite{Bre, HPR, Rin}. We shall need it to describe hammocks for hereditary  artin algebras of finite representation type. 

\begin{defn}

Let $(\Ga, \tau, v)$ be a valued translation quiver. A function $f: \Ga_0\to \Z$ is called {\it additive} on $\Ga$ provided, for any $x\in \Ga_0$ with $\tau x$ defined, that $f(\tau x)+f(x)=\sum_{y\in x^-} v'_{yx} f(y),\vspace{-1pt}$ where the sum is zero in case $x^-$ is empty.

\end{defn}

\begin{prop}\label{stab_af}

Let $(\Ga, \tau, v)$ be a stable valued translation quiver with a finite section $\Da$. Given $r_x\in \Z$ with $x\in \Da_0$, there exists a unique additive function $f$ on $\Ga$ such that $f(x)=r_x$ for all $x\in \Da_0$.

\end{prop}

\vspace{-4pt}

\noindent{\it Proof.} It suffices to establish the existence of $f$. We start with setting $f(x)=r_x\in \Z$, for all $x\in \Da_0.$ Given $s\in \Z$,  the vertices $\tau^sx$ with $x\in \Da_0$ generate a section $\tau^s\hspace{-2pt} \Da$ in $\Ga$. We shall first define $f$ on $\tau \Da$. Being finite without oriented cycles, $\Da$ contains a sink vertex $x_1$. By Lemma \ref{sec_property}(1), $x_1^-\subseteq \Da_0$. Put $f(\tau x_1)=\sum_{y\in x_1^-} v_{x_1\hspace{-.6pt},y}'f(y) - f(x_1)$. \vspace{-2pt} Clearly, the vertices in $\{\tau x_1\}\cup (\Da_0\backslash\{x_1\})$ generate a section $\Da^{(1)}$, on which $f$ is defined. Suppose that $f$ is defined on a section $\Da^{(i)}$ generated by the vertices in $\{\tau x_1, \ldots, \tau x_i\} \cup (\Da_0\backslash\{x_1, \ldots, x_i\})$ for some $i\ge 1$. If
$\Da_0=\{x_1, \ldots, x_i\}$, then $\Da^{(i)}=\tau \Da$. Otherwise, $\Da$ is non-trivial and there exists $x_{i+1}\in \Da_0\backslash \{x_1, \ldots, x_i\}$ such that $\Da^{(i)}$ contains no arrow $x_{i+1}\to y$ with $y\in \Da\backslash \{x_1, \ldots, x_i\}$. Suppose that $\Da^{(i)}$ contains an arrow $x_{i+1}\to \tau x_j$ for some $1\le j\le i$. Since $x_j^-\ne \emptyset$; see (\ref{conn_sec}), $\Ga$ contains a path $x_{i+1} \to \tau x_j \rightsquigarrow x_j$, and hence, $\tau x_j\in \Da$, absurd. Thus, $x_{i+1}$ is a sink in $\Da^{(i)}$. By Lemma \ref{sec_property}(1), $x_{i+1}^-\subseteq \Da^{(i)}$. \vspace{-2pt}  Define $f(\tau x_{i+1})=\sum_{y\in x_{i+1}^-} v_{x_{i+1}\hspace{-.6pt},y}'f(y) - f(x_{i+1}).$ Then, $f$ is defined on the section $\Da^{(i+1)}$ generated by the vertices in $\{\tau x_1, \ldots, \tau x_{i+1}\} \cup (\Da_0\backslash\{x_1, \ldots,  x_{i+1}\})$. By induction, we may define $f$ on $\tau \Da$. Repeating this process, we may define $f$ on $\tau^s\hspace{-2pt}\Da$ for all $s> 0$. Dually, considering sources in sections, we define $f$ on $\tau^{-s}\hspace{-2pt}\Da$ for all $s>0$. This yields a desired additive function $f$ on $\Ga$. The proof of the proposition is completed.

\section{Hereditary artin algebras}

The main objective of this section is to study the preprojective and preinjective components of a connected hereditary artin algebra. Although they have already been well studied; see \cite[Section VIII.1]{ARS} and \cite[Pages 80 - 81]{Ri1}, we shall further study them in terms of the ext-quiver of the algebra. In case the ext-quiver is a tree, a sectional path in the preprojective component is the unique path between the end-points, two parallel paths have the same length,
and every projective module is the source of a source section, which is explicitly described in terms of reduced walks in the ext-quiver;
see (\ref{pic_sec_path}) and (\ref{prep_cpt_path}). Finally, we shall show that the Coxeter transformation is compatible with the derived Auslander-Reiten translation; see (\ref{Cox_tau}).

Throughout this section, $H$ stands for a connected hereditary artin algebra with ext-quiver $\QH$, and $\mmod H$ for the category of finitely generated left $H$-modules. We simply write $\tau$ for the Auslander-Reiten translation $\tau_{\hspace{-.6pt}_H}$ of the Auslander-Reiten quiver $\GaH$. And we associate a simple module $S_i$, a projective module $P_i$ and an injective module $I_i$ in $\GaH$  with each $i\in (\QH)_0$ such that ${\rm top}P_i\cong S_i\cong {\rm soc} I_i$.

\subsection{\sc Preprojective and preinjective components}
The following result is probably well-known. We sketch a proof for the first statement and refer to \cite[Page 267]{ARS} for a proof for the second.

\begin{lem}\label{pi_endo}

Let $H$ be a hereditary artin algebra with ext-quiver $\QH$. Consider the simple modules $S_i$, the projective modules $P_i$ and the injective modules $I_i$ in $\GaH$ associated with $i\in \QH$.

\begin{enumerate}[$(1)$]

\vspace{-1.5pt}

\item If $i$ is a vertex in $\QH$, then $\End_H(P_i)\cong \End_H(I_i) \cong \End_H(S_i)$ and $\ell_{S_i}(P_i)=\ell_{S_i}(I_i)=1$.

\vspace{.5pt}

\item There exists an arrow $i\to j$ in $\QH$ if and only if there exists an irreducible map $f: P_j\to P_i$, or equivalently, an irreducible map $g: I_j\to I_i$, in $\mmod H$.

\end{enumerate}\end{lem}

\vspace{-4pt}

\noindent{\it Proof.} Fix $i\in (\QH)_0$. Since  $I_i$ is not a direct summand of the injective module $I_i/S_i$, we have $\Hom_H(S_i, I_i/S_i)=0.$ Thus, we deduce from Lemma \ref{m_com_fac} that $\ell_{S_i}(I_i)=1$. Since $S_i\cong {\rm soc} I_i$, every non-zero map in $\End_H(S_i)$ induces a non-zero map in $\End_H(I_i)$. Consider a non-zero map $f: I_i\to I_i$. Since ${\rm Im}(f)$ is injective, $f$ is invertible. Since $\Hom_H(S_i, I_i/S_i)=0$, we see that $f$ is induced from a non-zero map $g: S_i\to S_i$.  Thus, $\End_H(I_i) \cong \End_H(S_i)$. 
The proof of the lemma is completed.

\vspace{3pt}

Since $\QH$ is connected, by Lemma \ref{pi_endo}(2), the projective modules in $\Ga_{{\rm mod}H}$ lie in the same connected component $\mathscr{P}_H$, called the {\it preprojective component}; and the injective modules lie in the same connected component  $\mathscr{I}_H$, called the {\it preinjective component}. The following statement is reformulated from Proposition 1.15 in \cite[Chapter VIII]{ARS} and its dual statement. Since our valuation for $\GaH$ is different from that given in \cite[Page 225]{ARS}, we include a detailed proof for the preinjective component.

\begin{thm}\label{pisec}

Let $H$ be a connected hereditary artin $R$-algebra with
ext-quiver $\QH$. Consider the projective modules $P_i$ and the injective modules $I_i$ in $\GaH$ associated with $i\in \QH$.

\begin{enumerate}[$(1)$]

\vspace{-1pt}

\item The projective modules $P_i$ with $i\in \QH$ generate a leftmost section $\DaH$ in the preprojective component $\mathscr{P}_H$. Moreover, there exists a valued quiver isomorphism $\QoH\to \DaH,$ sending $i$ to $P_i$.

\vspace{.5pt}

\item The injective modules $I_i$ with $i\in \QH$ generate a rightmost section $\SaH$  in the preinjective component $\mathscr{I}_H$. Moreover, there exists a valued quiver isomorphism $\QoH\to \SaH,$ sending $i$ to $I_i$.


\end{enumerate}\end{thm}

\vspace{-3pt}

\noindent{\it Proof.} We shall only prove Statement (2). We write $\ell_R(V)$ for the $R$-length of an $R$-module $V,$ and put $k_M={\rm End}_H(M)/\rad({\rm End}_H(M))$ for $M\in \Ga_{{\rm mod}H}$. Let $\SaH$ be the full valued subquiver of $\mathscr{I}_H$ generated by the $I_i$ with $i\in \QH$. Set $S_i={\rm soc}I_i$ for $i\in \QH$. By Lemma \ref{pi_endo}(2), $\QH$ contains an arrow $i\to j$ if and only if $\SaH$ contains an arrow $I_j\to I_i$; and in this case, the valuation for $i\to j$ is $(d_{ij}, d_{ij}')$, where $d_{ij}={\rm dim}_{{\rm End}_H(S_j)}\Ext_H^1(S_i, S_j)$ and $d_{ij}'={\rm dim} \Ext_H^1(S_i, S_j)_{{\rm End}_H(S_i)}.$ Thus, \vspace{-1pt} $$d_{ij}' \cdot \ell_R(\End_H(S_i))=\ell_R(\Ext_H^1(S_i, S_j))=d_{ij} \cdot \ell_R(\End_H(S_j)).
\vspace{-2pt}$$

On the other hand, we see from Proposition \ref{Irr_mor} that $I_j\to I_i$ has valuation $(d_{I_{\hspace{-.8pt}j\hspace{-.6pt}}, I_i}\hspace{-.5pt}, d'_{I_{\hspace{-.8pt}j\hspace{-.6pt}}, I_i})$, where $d_{I_{\hspace{-.8pt}j\hspace{-.6pt}}, I_i}={\rm dim}\, _{k_{I_i}}{\rm irr}(I_j, I_i)$ and $d_{I_{\hspace{-.8pt}j\hspace{-.6pt}}, I_i}'={\rm dim} \, {\rm irr}(I_j, I_i) _{k_{I_j}}\hspace{-2pt}.$ Since $k_{I_j}\cong \End_H(S_j)$ and $k_{I_i}\cong \End_H(S_i)$ by Lemma \ref{pi_endo}(1), $d_{I_j, I_i}' \cdot \ell_R(\End_H(S_j))=\ell_R({\rm Irr}(I_j, I_i))=d_{I_j, I_i} \cdot \ell_R(\End_H(S_i)).$

Now since $I_j / S_j$ is injective, by Proposition 1.15(b) in \cite[(III.1.15)]{ARS}, $d_{ij}'$ is the multiplicity of $I_i$ as a direct summand of $I_j / S_j$. And since the canonical projection $\varphi\hspace{-1pt} : \hspace{-1pt} I_j\to I_j/S_j$ is minimal left almost split, $d_{I_{\hspace{-.8pt}j\hspace{-.6pt}}, I_i}$ is the multiplicity of $I_i$ as a direct summand of $I_j / S_j$; see ( \ref{Irr_mor}). So, $d_{ij}'=d_{I_{\hspace{-.8pt}j\hspace{-.6pt}}, I_i}$. Hence,\vspace{-3pt}
$$\hspace{30pt} d_{I_{\hspace{-.8pt}j\hspace{-.6pt}}, I_i}' \cdot \ell_R(\End_H(S_j))= d_{I_{\hspace{-.8pt}j\hspace{-.6pt}}, I_i} \cdot \ell_R(\End_H(S_i))=
d_{ij}' \cdot \ell_R(\End_H(S_i)) = d_{ij} \cdot \ell_R(\End_H(S_j)).
\vspace{-3pt}$$

Thus, $d_{I_{\hspace{-.8pt}j\hspace{-.6pt}}, I_i}'=d_{ij}$. This proves the second part of Statement (2). In particular, $\Sa_{\hspace{-.8pt}H}$ is connected. Being injective, the $I_i$ lie in pairwise distinct $\tau$-orbits in $\GaH$. Note that the irreducible maps between the modules $I_i$ are all epimorphisms. Thus, since $\End_H(I_i)$ is divisible; see (\ref{pi_endo}), $\SaH$ has no oriented cycle. Hence, $\SaH$ is a section in $\mathscr{I}_H$, which is clearly rightmost. The proof of the theorem is completed.






\vspace{3pt}

In the sequel, we shall call the underlying valued graph $\OQH$ of $\QH$ the {\it ext-graph} of $H$. Moreover, $\DaH$ will be called the {\it projective section} in $\mathscr{P}_H$, and $\SaH$ called the {\it injective section} in $\mathscr{I}_H$. The following well-known statement is a consequence of Theorem \ref{pisec} and well-known; see \cite[(V.7.5), (VIII.1.9)]{ARS}.

\begin{cor}\label{rad_maps}

Let $H$ be a connected hereditary artin algebra. Consider modules $M, N$ in $\GaH$ with $M$ in the preinjective component $\mathscr{I}_H$ or $N$ in the preprojective component $\mathscr{P}_H$. Then, every non-zero radical map $f: M\to N$ is a sum of composites of irreducible maps between modules in $\GaH$.
\end{cor}


We shall need the following easy statement.

\begin{lem}\label{arrow_translate}

Let $H$ be a hereditary artin algebra. Consider an arrow $M\to N$ in $\GaH$.

\begin{enumerate}[$(1)$]

\vspace{-2pt}

\item If $\GaH$ contains $\tau^{r}\hspace{-1pt}M$ for some $r\ge 0$, then it contains an arrow $\tau^{r}\hspace{-1pt}M\to \tau^{r}\hspace{-1pt}N.$

\item If $\GaH$ contains $\tau^{-r}\hspace{-1pt}N$ for some $r\ge 0$, then it contains an arrow $\tau^{-r}\hspace{-1pt}M\to \tau^{-r}N$.

\end{enumerate} \end{lem}

\vspace{-4pt}

\noindent{\it Proof.} We shall only prove Statement (1). Assume that
 $\tau^{r}\hspace{-1pt}M \hspace{-1pt}\in \hspace{-1pt} \GaH$ for some $r>0$, but $\tau^{r}\hspace{-1pt}N \hspace{-1pt}\not\in \hspace{-1pt} \GaH$. Then, $\tau^{s}\hspace{-1pt}N$ is projective for some $0\le s<r$. Since $\tau^{s}M\in \GaH$, by Lemma \ref{pa_valuation}, $\GaH$ contains an arrow $\tau^{s}\hspace{-1pt}M\to \tau^{s}\hspace{-1pt}N$. Since $H$ is hereditary, $\tau^{s}\hspace{-1pt}M$ is projective, and hence, $\tau^{r}\hspace{-1pt}M\not\in \GaH$, absurd. Thus, $\tau^{r}\hspace{-1pt}N\in \GaH$. So, $\GaH$ contains an arrow $\tau^{r}\hspace{-1pt}M \to \tau^{r}\hspace{-1pt}N$. The proof of the lemma is completed.

\vspace{3pt}

The following statement describes in terms of $\QH$ the sectional paths in $\GaH$, which start with a projective module or end with an injective module. 

\begin{lem}\label{pic_sec_path}

Let $H$ be a hereditary artin algebra with ext-quiver $\QH$. Consider the projective modules $P_i$ and the injective modules $I_i$ in $\GaH$ associated with $i\in \QH$. And let $w=c_t\cdots c_2c_1$ be a non-trivial reduced walk in $\QH$, where $c_j$ is an arrow or inverse arrow from $i_{j-1}$ to $i_j$, for $j=1, \ldots, t$.

\begin{enumerate}[$(1)$]

\vspace{-1pt}

\item The preprojective component $\mathscr{P}_H$ of $\GaH$ contains a sectional path \vspace{-4pt}
$$\xymatrixcolsep{16pt}\xymatrix{P_{i_0} \ar[r] & \tau^{-r_{1}} \! P_{i_1} \ar[r] & \cdots \ar[r] & \tau^{-r_{t-1}}\!P_{i_{t-1}} \ar[r] & \tau^{-r_{t}}\!P_{i_t},}
\vspace{-4pt}$$ where $r_j$ is the number of arrows in the subwalk of $w$ from $i_0$ to $i_j$, for $j=1, \ldots, t$.

\vspace{1.5pt}

\item The preinjective component $\mathscr{I}_H$ of $\GaH$ contains a sectional path
\vspace{-3pt} $$\xymatrixcolsep{16pt}\xymatrix{\tau^{s_{0}}\! I_{i_0} \ar[r] & \tau^{s_{1}} \! I_{i_1} \ar[r] &  \cdots \ar[r] & \tau^{s_{t-1}}\!I_{i_{t-1}} \ar[r] & I_{i_t},}  \vspace{-4pt} $$ where $s_j$ is the number of arrows in the subwalk of $w$ from $i_j$ to $i_t$, for $j=0, 1, \ldots, t-1$.

\vspace{-2pt}

\end{enumerate} \end{lem}

\noindent{\it Proof.} We shall only prove Statement (1). Since $w$ is reduced,
the $P_{i_j}$ with $0\le j\le t$ lie in pairwise distinct $\tau$-orbits in $\mathscr{P}_H$. \vspace{-2pt} Set $r_0=0$. Assume that $\mathscr{P}_H$ contains a desired sectional path $\xymatrixcolsep{18pt}\xymatrix{\hspace{-3pt} P_{i_0} \ar[r]
& \cdots \ar[r] & \tau^{-r_{t-1}}\!P_{i_{t-1}}.}$ \vspace{-2pt} Assume first that $\QH$ contains an arrow $i_{t-1}\to i_t$. Then $r_t=r_{t-1}+1$ and by Theorem \ref{pisec}(1), $\mathscr{P}_H$ contains an arrow $P_{j_t} \to P_{j_{t-1}}.$ By
Lemma \ref{arrow_translate}(2), $\mathscr{P}_H$ contains an arrow
$\tau^{-r_{t-1}} \! P_{i_t} \!\to\! \tau^{-r_{t-1}} \! P_{i_{t-1}}$.
Since $r_{t-1}+1 = r_t,$ \vspace{-2.5pt} it follows from Lemma \ref{p_sectional_path}(1) that $\mathscr{P}_H$ contains a sectional path
$\xymatrixcolsep{16pt}\xymatrix{\hspace{-2pt} P_{i_0} \ar[r] & \cdots \ar[r] & \tau^{-r_{t-1}}\!P_{i_{t-1}} \ar[r] & \tau^{-r_t}\!P_{i_t}.}$

\vspace{-2pt}

Assume now that $\QH$ contains an arrow $i_{t-1} \leftarrow i_t$. Then $r_t=r_{t-1}$ and $\mathscr{P}_H$ contains an arrow $P_{i_{t-1}}\to P_{i_t}$. We claim that $\mathscr{P}_H$ contains an arrow $\tau^{-r_{t-1}} \! P_{i_{t-1}} \!\to\! \tau^{-r_t}\! P_{i_t}$.
This is the case if $r_{t-1}=0$. Otherwise, $\mathscr{P}_H$ has an arrow $P_{i_t} \to \tau^-\hspace{-1.2pt} P_{i_{t-1}}$.
Since $\tau^{1-r_{t-1}}(\tau^-P_{i_{t-1}}) = \tau^{-r_{t-1}} \! P_{i_{t-1}}$, by Lemma
\ref{arrow_translate}(1), $\mathscr{P}_H$ contains an arrow $\tau^{1-r_{t-1}} \! P_{i_t}\to \tau^{-r_{t-1}} \! P_{i_{t-1}}$. By Lemma \ref{p_sectional_path}(1), $\tau^{1-r_{t-1}} \! P_{i_t}$ is not injective. Since $r_t=r_{t-1}$, there exists an arrow $\tau^{-r_{t-1}} \! P_{i_{t-1}} \!\to\! \tau^{-r_t}\! P_{i_t}$ in $\mathscr{P}_H$. \vspace{-2.5pt} This establishes our claim. Hence, $\mathscr{P}_H$ contains a sectional path \hspace{-5pt} $\xymatrixcolsep{16pt}\xymatrix{P_{i_0} \ar[r] & \cdots \ar[r] & \tau^{-r_{t-1}}\!P_{i_{t-1}} \ar[r] & \tau^{-r_t}\!P_{i_t}.}$\hspace{-2pt} The proof of the lemma is completed.

\vspace{2pt}

The following statement is crucial for us to study hammocks in the next section.

\begin{prop} \label{special section}

Let $H$ be a connected hereditary artin algebra. 

\begin{enumerate}[$(1)$]

\vspace{-1pt}

\item The preprojective component $\mathscr{P}_H$ contains a $P$-source section for every projective module $P$ in $\mathscr{P}_H$.

\vspace{.5pt}

\item  The preinjective component $\mathscr{I}_H$ contains an $I$-sink section for every injective module $I$ in $\mathscr{I}_H$.

\end{enumerate}\end{prop}

\vspace{-4pt}

\noindent{\it Proof.} We shall only prove Statement (1). Consider a module $P$ in the projective section $\DaH$ in $\mathscr{P}_H$. Let $\Sa$ be the full valued subquiver of $\mathscr{P}_H$ generated by the successors $M$ of $P$ in $\mathscr{P}_H$ such that $\tau M$ is not a succesor of $P$. Then, $\Sa$ is contains no oriented cycle and meets any $\tau$-orbit in $\mathscr{P}_H$ at most once. Given $P'\in \DaH$, since the ext-quiver of $H$ is connected, we deduce from Lemma \ref{pic_sec_path}(1) that there exists a minimal $s\ge 0$ such that $\tau^{-s} \hspace{-1.5pt}P'$ is a successor of $P$ in $\mathscr{P}_H$. Then $\tau^{-s} \hspace{-1.5pt} P'\in \Sa$. Hence, $\Sa$ meets every $\tau$-orbit in $\mathscr{P}_H$ exactly once. Finally, suppose that $\mathscr{P}_H$ contains a path $\eta: M\rightsquigarrow L \rightsquigarrow N$ with $M, N\in \Sa$. Then, $\mathscr{P}_H$ contains a path $P\rightsquigarrow M \rightsquigarrow L$. Assume that $L\not\in \Sa$. Then, $\mathscr{P}_H$ contains a path $P\rightsquigarrow \tau L$. 
Since $H$ is hereditary, $L\rightsquigarrow N$ contains no projective module. Hence, $\mathscr{P}_H$ contains a path $P\rightsquigarrow \tau L \rightsquigarrow \tau N$, a contradiction. Therefore, $\Sa$ is convex in $\mathscr{P}_H$. In particular, every module in $\Sa$ is a succesor of $P$ in $\Sa$. Hence, $\Sa$ is a $P$-source section in $\mathscr{P}_H$.
The proof of the proposition is completed.

\subsection{\sc Tree type} We shall say that $H$ is of {\it tree type} or {\it Dynkin type} if $\QH$ is a tree or a Dynkin quiver, respectively. These hereditary algebras have the following important properties.

\begin{prop}\label{prep_cpt_path}

Let $H$ be a connected hereditary artin algebra of tree type. Consider the preprojective component $\mathscr{P}_H$ and the preinjective component $\mathscr{I}_H$ of $\GaH$.

\begin{enumerate}[$(1)$]

\vspace{-1.5pt}

\item If $\eta, \zeta$ are parallel paths in $\mathscr{P}_H$ or $\mathscr{I}_H$, then $l(\eta)=l(\zeta);$ and $\eta=\zeta$ in case $\eta$ is sectional.

\item The $P$-source section in $\mathscr{P}_H$ with $P$ projective is generated by the sectional successors of $P$ in $\mathscr{P}_H$.

\item The $I$-sink section in $\mathscr{I}_H$ with $I$ injective is generated by the sectional predecessors of $I$ in $\mathscr{I}_H$.
\end{enumerate}

\end{prop}

\vspace{-4pt}

\noindent {\it Proof.} (1) Let $\DaH$ be the projective section in $\mathscr{P}_H$. By Proposition \ref{sec_property}(2), $\mathscr{P}_H$ embeds in $\mathbb Z\DaH$ as a full valued translation subquiver. In particular, a (sectional) path in $\mathscr{P}_H$ gives rise to a (sectional) path in $\mathbb Z\DaH$. Since $\DaH$ is a tree; see (\ref{pisec}), we see that Statement (1) follows from Proposition \ref{uni_sec_path}.

(2) Given $P\in \DaH$, by Proposition \ref{special section}, $\mathscr{P}_H$ contains a $P$-source section $\Sa$. Clearly, every module in $\Sa$ is a sectional successor $P$. Suppose that $\mathscr{P}_H$ contains a sectional path $\eta: P\rightsquigarrow M$. By Lemma \ref{s_sec_unicity}, $M=\tau^{-s}\hspace{-1pt}N$ for some $N\in \Sa$. Then, $\mathscr{P}_H$ contains a path $\zeta: P\rightsquigarrow N \rightsquigarrow M$. By Statement (1), $\zeta=\eta$. So, $\zeta$ is sectional, and hence, $M=N\in \Sa$. The proof of the proposition is completed.

\vspace{3pt}

The following statement will enable us to define extended hammocks in the next section.

\begin{lem} \label{sp_mul}

Let $H$ be a connected hereditary artin algebra of tree type. Consider a sectional path $\hspace{-2pt}\xymatrixcolsep{18pt}\xymatrix{P=M_0\ar[r] & M_1 \ar[r] & \cdots \ar[r] & M_{t-1}\ar[r] & M_t}\hspace{-2pt}\vspace{-1pt}$ in $\GaH$, where $P$ is a projective module with top $S$, and $M_{i-1}\to M_i$ with $1\le i\le t$ has valuation $(d_i, d_i')$. Then $\ell_S(M_i)=d_i' \cdots d_1',$ for $i=1, \ldots, t.$

\end{lem}

\vspace{-5pt}

\noindent{\it Proof.} By Lemma \ref{pi_endo}(1), $\ell_S(M_0)=1.$ Suppose that $t>0$. 
Consider a minimal right almost split map $f: L\to M_t$. Since $P\not\cong M_t$,
the ${\rm End}_H(P)$-linear map $f_*: \Hom_H(P, L) \to \Hom_H(P, M_t)$ is an epimorphism. We claim that it is an isomorphism. If $M_t$ is projective, then $f$ is a monomorphism, and so is $f_*$. Otherwise, we have an exact sequence \vspace{-6pt}
$$\xymatrix{0\ar[r] & \Hom_H(P, \tau M_t) \ar[r] &\Hom_H(P, L) \ar[r]^-{f_*} & \Hom_H(P, M_t) \ar[r] & 0.}
\vspace{-5pt}$$

Since the path stated in the lemma is sectional, by Lemma \ref{prep_cpt_path}(1), $\mathscr{P}_H$ contains no path from $P$ to $\tau M_t$. Hence, $\Hom_H(P, \tau M_t)=0$; see (\ref{rad_maps}). This establishes our claim. Therefore,  $\ell_S(M_t)=\ell_S(L)$; see (\ref{m_com_fac}). On the other hand, by Proposition \ref{Irr_mor}, $L=M_{t-1}^{d_t'} \oplus L_1 \cdots \oplus L_r$, where $L_1 \ldots, L_r$ with $r\ge 0$ are modules in $\mathscr{P}_H$ different from $M_{t-1}$. By Lemma \ref{prep_cpt_path}(1), $\mathscr{P}_H$ contains no path from $P$ to $L_j$, and hence,
$\Hom_H(P, L_j)=0$ and $\ell_S(L_j)=0$, for $j=1, \ldots, r$. By the induction hypothesis, we see that $\ell_S(M_t)=\ell_S(L)=d_t' \ell_S(M_{t-1})=d_t' \cdots d_1'.$ The proof of the lemma is completed.

\subsection{\sc Coxeter transformation} Let us recall Auslander and Platzeck's Coxeter transformation of the Grothendieck group $K_0(\mmod H)$. Consider the simple modules $S_i$, the projective modules $P_i$ and the injective mdoules in $\GaH$ associated with $i\in (\QH)_0=\{1, \ldots, n\}.$ For any module $M$ in $\mmod H$, one defines its {\it dimension vector} to be $\undim M: = (\ell_{S_1}(M), \ldots, \ell_{S_n}(M))\in K_0(\mmod H).$ Since $\{\undim P_1, \ldots, \undim P_n\}$ and $\{\undim I_1, \ldots, \undim I_n\}$ are bases for $K_0(\mmod H)$; see \cite[(2.1)]{PAu}, one may introduce the following definition; see \cite[Section 2]{PAu}, and compare \cite[Page 8]{DRi}.

\begin{defn}\label{Cox_trans} Let $H$ be a hereditary artin algebra. The {\it Coxeter transformation} of the Grothendieck group $K_0(\mmod H)$ is the unique automorphism
$$C_H: K_0(\mmod H)\to K_0(\mmod H) \mbox{ defined by }  C_H(\undim P_i)=-\hspace{.8pt}\undim I_i, \mbox{ for } i=1, \ldots, n.$$

\end{defn}

An important property of the Coxeter transformation $C_H$ of $K_0(\mmod H)$ is its compatibility with the Auslander-Reiten translation $\tau$ of $\GaH$; see \cite[(2.2)]{PAu}, and also \cite[(VIII.2.2)]{ARS}. We shall extend this to the bounded derived category $D^{\hspace{.4pt}b\hspace{-.5pt}}(\mmod H)$ of $\mmod H$, which is a Hom-finite Krull-Scmidt $R$-category.
It is well-known; see \cite[(7.3)]{BLP}, \cite[(3.6)]{Hap} and \cite[(I.3.3)]{RVDB} that $D^{\hspace{.4pt}b\hspace{-.5pt}}(\mmod H)$ has almost split triangles as defined in \cite[(3.1)]{Hap}. Note that a sequence of morphisms $\hspace{-3pt}\xymatrixcolsep{17pt}\xymatrix{L^\pdt \ar[r] & M^\pdt \ar[r] & N^\pdt}\hspace{-2pt}$ \vspace{-1.5pt} in $D^{\hspace{.4pt}b\hspace{-.5pt}}(\mmod H)$ with $M^\pdt\ne 0$ is an almost split sequence if and only if it embeds in an almost split triangle $\hspace{-3pt}\xymatrixcolsep{17pt}\xymatrix{L^\pdt \ar[r] & M^\pdt \ar[r] & N^\pdt\ar[r] & L^\pdt[1]}\vspace{-2pt}\hspace{-2pt}$; see \cite[(6.1)]{Liu}. Thus, the Auslander-Reiten quiver $\Ga_{\hspace{-1pt}D^{\hspace{.4pt}b\hspace{-.5pt}}(\mmod H)}$ of $D^{\hspace{.4pt}b\hspace{-.5pt}}(\mmod H)$ as defined in Subsection 1.5 coincides with that defined by Happel in the algebraically closed setting; see \cite[(3.7)]{Hap}.

Given $M\in \mmod H$ and $s\in \mathbb Z$, we write $M[s]$ for the stalk complex concentrated in degree $-s$ where the component is $M$. Then, the indecomposable objects in $D^{\hspace{.4pt}b\hspace{-.5pt}}(\mmod H)$ are the stalk complexes $M[s]$, where $s\in \Z$ and $M$ is an indecomposable module in $\mmod H$; see \cite[(3.1)]{Len}. Thus, we may choose the vertices in $\Ga_{\hspace{-1pt}D^{\hspace{.4pt}b\hspace{-.4pt}}(\mmod H)}$ to be the stalk complexes $M[s]$ with $M\in \GaH$ and $s\in \Z$. We write $\tau_{\hspace{-.5pt}_D}$ for the Auslander-Reiten translation of $\Ga_{\hspace{-1pt}D^{\hspace{.4pt}b\hspace{-.5pt}}(\mmod H)}$. On the other hand, given a complex $M^{\hspace{.5pt}\dt}\hspace{2pt}$ in $D^{\hspace{.4pt}b\hspace{-.5pt}}(\mmod H)$, Happel defined its {\it dimension vector} by
$$\undim M^{\hspace{.8pt}\dt}:=\textstyle\sum_{s\in \mathbb{Z}} (-1)^s \undim M^s\in K_0(\mmod H),$$ which is invariant on isomorphism classes in $D^{\hspace{.4pt}b\hspace{-.5pt}}(\mmod H);$
see \cite[(2.2)]{Hap}.
As shown below, the Coxeter transformation of $K_0(\mmod H)$ is compatible with the Auslander-Reiten translation of $\Ga_{\hspace{-1pt}D^{\hspace{.4pt}b\hspace{-.5pt}}(\mmod H)}$.

\begin{prop} \label{Cox_tau}

Let $H$ be a hereditary artin algebra. Consider the Coxeter transformation $C_{\hspace{-.5pt}H}$ of $K_0(\mmod H)$ and the Auslander-Reiten translation $\tau_{\hspace{-.5pt}_D}$ of $\Ga_{\hspace{-1pt}D^{b\hspace{-.5pt}}(\mmod H)}$.
Then, $\undim\hspace{.5pt}\tau_{\hspace{-.5pt}_D}^{\hspace{1pt}t\hspace{-1.5pt}}(M[s])=C_{\hspace{-.5pt}H}^{\hspace{1pt}t\hspace{-1.5pt}}(\undim M[s]),$ for all $M\in \GaH$ and $s, t\in \mathbb{Z}$.

\end{prop}

 \vspace{-4pt}

\noindent{\it Proof.} We shall only show that $\undim\hspace{.5pt}\tau_{\hspace{-.5pt}_D}\hspace{-1pt}(M[s]) = C_{\hspace{-.5pt}H} \hspace{.5pt} (\undim M[s])$, for $M\in \GaH$ and $s\in \Z$. Suppose first that $M$ is not projective. Then, $\underline\dim \hspace{.4pt} \tau M = C_{\hspace{-.5pt}H}(\underline\dim M)$; see \cite[(2.2)]{PAu}, and also \cite[(VIII.2.2)]{ARS}. Since $\tau_{\hspace{-.5pt}_D}$ commutes with the shift functor $[1]$, we have $\tau_{\hspace{-.5pt}_D} \hspace{-.5pt}(\hspace{-.5pt}M[s])=(\tau_{\hspace{-.5pt}_D} M)[s]=(\tau M)[s]$; see \cite[(7.2)]{BLP}, and also \cite{Hap}. By definition, $\undim N[s]=(-1)^s \undim N$ for any $N\in \mmod H$. Thus,
\vspace{-2pt}  $$\undim \tau_{\hspace{-.5pt}_D}(M[s]) \hspace{-1pt} = \hspace{-1pt}
 \undim (\tau M)[s] \hspace{-1pt} = \hspace{-1pt} (-1)^s \undim \tau M \hspace{-1.5pt} = \hspace{-1.5pt} (-1)^s C_{\hspace{-.5pt}H}\hspace{.5pt}(\undim  M) \hspace{-1pt} = \hspace{-1pt} C_{\hspace{-.5pt}H}\hspace{.5pt}((-1)^s\undim M) \hspace{-1pt} = \hspace{-1pt} C_{\hspace{-.5pt}H}(\undim M[s]). \vspace{-2pt} $$

Suppose now that $M=P$, the projective cover of a simple module $S$ in $\GaH$. By definition, $C_{\hspace{-.5pt}H}\hspace{.5pt}(\undim  P)=-\hspace{1pt} \undim I$, where $I$ is the injective envelope of $S$. Moreover, $\tau_{\hspace{-.5pt}_D}(P[0])= I[-1];$ see \cite[(7.2)]{BLP}, and also \cite{Hap}. This yields \vspace{-2pt}
$$\undim  \tau_{\hspace{-.5pt}_D}\hspace{-.5pt}(P[s]) \hspace{-1pt} = \hspace{-1pt} \undim (\tau_{\hspace{-.5pt}_D}\hspace{-.5pt}P[0])[s] \hspace{-1pt} = \hspace{-1pt} \undim I[s \hspace{-1pt} - \hspace{-1pt} 1]  \hspace{-1pt} = \hspace{-1pt} (-1)^s C_{\hspace{-.5pt}H}\hspace{.5pt}(\undim P) \hspace{-1pt} = \hspace{-1pt} C_{\hspace{-.5pt}H}\hspace{.5pt}((-1)^s \undim  P)
\hspace{-1pt} = \hspace{-1pt} C_{\hspace{-.5pt}H}\hspace{.5pt}(\underline\dim \hspace{.4pt} P[s]).$$ The proof of the proposition is completed.

\section{Main Results}

The objective of this section is to study the module category of a hereditary artin algebra of finite representation type with a connection to the derived category and the associated cluster category. We first determine all the hammocks in the Auslander-Reiten quiver; see (\ref{hamm_neg_section}). This leads to a description of the precise shape of the Auslander-Reiten quiver in terms of the ext-quiver of the algebra; see (\ref{ARQ_Dyn}). As applications, we obtain the number of non-isomorphic indecomposable objects in the module category and the associated cluster category; see (\ref{rep_nb}) and (\ref{clus_nil}). Moreover, the radicals of the module category, the bounded derived category and the associated cluster category all have the same nilpotency; see (\ref{Mod_nil}), (\ref{Der_nil}), and (\ref{clus_nil}).

It is well-known that a hereditary artin algebra is connected of finite representation type if and only if its ext-graph is a Dynkin diagram; see \cite[(VIII.5.4)]{ARS}. In the algebraically closed case, this is equivalent to the ext-graph being $\mathbb A_n (n\ge 1)$, $\mathbb D_n (n\ge 4)$ or $\mathbb E_n (n=6,7,8)$; see \cite[(VIII.5.5)]{ARS}.

Throughout this section, $H$ stands for a hereditary artin algebbra with a Dynkin ext-quiver $\QH$, and $\mmod H$ for the category of finitely generated left $H$-modules. Let $\GaH$ and $\Ga_{\hspace{-.8pt}D^{\hspace{.4pt}b\hspace{-.5pt}}(\mmod H)}$ be the Auslander-Reiten quivers of $\mmod H$ and $D^{\hspace{.4pt}b\hspace{-.5pt}}(\mmod H)$, respectively. We shall simply write $\tau$ for the Auslander-Reiten translation $\tau_{\hspace{-.6pt}_H}$ of $\GaH$, and $\tau_{\hspace{-.5pt}_D}$ for the Auslander-Reiten translation of $\Ga_{\hspace{-.8pt}D^{\hspace{.4pt}b\hspace{-.5pt}}(\mmod H)}$. With each vertex $i\in \QH$, we associate a simple module $S_i$, a projective module $P_i$ and an injective module $I_i$ in $\GaH$  such that ${\rm top}P_i\cong S_i\cong {\rm soc} I_i$.

\subsection{\sc The canonical embedding} Since $H$ is of finite representation type, it is well-known that $\GaH$ coincides with its preprojective component $\mathscr{P}_H$, and its preinjective component $\mathscr{I}_H$; see, for example, \cite[(VIII.3.13)]{ARS}. The following statement provides a rough description of the shape of $\GaH$.

\begin{prop}\label{H_alg_Dynkin}

Let $H$ be a hereditary artin algebra with a Dynkin ext-quiver $\QH$. Consider the projective modules $P_i$ in $\GaH$ associated with vertices $i\in \QH$. Then there exists a valued translation quiver embedding $\varphi: \GaH \to \ZQH: \tau^{-r}\hspace{-1.5pt}P_i\mapsto (r, i),$ with a convex image in $\ZQH.$ Moreover, $\varphi$ maps the $P_i$-source section in $\GaH$ onto the $(0,i)$-source section in $\ZQH$, for every vertex $i\in \QH$.

\end{prop}

\vspace{-4pt}

\noindent{\it Proof.} Let $\DaH$ be the projective section in $\GaH$. Then, there exists valued translation quiver embedding $\psi: \GaH \to \Z \DaH: \tau^{-r}\hspace{-1.5pt} P_i\mapsto (r, P_i),$ whose image is convex in $\Z \DaH;$ see \cite[(1.13)]{LiY}. Since $\DaH\cong Q^{\rm op}_{\hspace{-.5pt}H}$; see (\ref{pisec}), we have a valued quiver isomorphism $\theta:  \Z \DaH\to \Z Q^{\rm op}_{\hspace{-.5pt}H}: (r,P_i) \mapsto (r, i)$. This yields an embedding $\varphi=\theta\circ \psi: \GaH \to \Z Q^{\rm op}_{\hspace{-.6pt}H}: \tau^{-r}\hspace{-1.5pt}P_i\mapsto (r, i),$ whose image is convex in $\ZQH.$


Fix $i\in (\QH)_0$. The $P_i$-source section $\Da_i$ in $\GaH$ is generated by the sectional successors of $P_i$ in $\GaH$; see (\ref{prep_cpt_path}), and the $(0,i)$-source section $\Oa_i$ in $\ZQH$ is generated by the sectional successors of $(0,i)$ in $\ZQH$; see (\ref{ZD_sec}). In particular, $\varphi$ maps to $\Da_i$ into $\Oa_i$. Let $(r, j)\in \Oa_i$. Since $\Sa_i$ is a section in $\GaH$, there exists some $s\ge 0$ such that $\tau^{-s} \hspace{-1pt} P_j\in \Da_i$. Thus,
$(s, j)\in \Oa_i$. Since $\Oa_i$ is a section in $\ZQH$, we have $s=r$. So, $\tau^{-r}P_j$ is a preimage of $(r,j)$ in $\Da_i$. The proof of the proposition is completed.

\vspace{3pt}

The following statement is due to Happel in the algebraically closed case; see \cite[(4.5)]{Hap}.

\begin{thm}\label{Der_ARQ}

Let $H$ be a hereditary artin algebra of Dynkin type. Then

\begin{enumerate}[$(1)$]

\vspace{-1.5pt}

\item $\GaH$ embeds in $\Ga_{\hspace{-.8pt}D^{\hspace{.4pt}b\hspace{-.5pt}}(\mmod H)}$ as a convex valued translation subquiver$\hspace{1pt};$

\vspace{.5pt}

\item $\Ga_{\hspace{-.8pt}D^{\hspace{.4pt}b\hspace{-.5pt}}(\mmod H)}\cong \mathbb Z \Da_H$, where $\DaH$ is the projective section in $\Ga_{\hspace{-1pt}{\rm mod}H}.$

\end{enumerate}\end{thm}

\vspace{-4pt}

\noindent {\it Proof.} It is well-known that there exists a full convex embedding of $\mmod H$ in $D^{\hspace{.4pt}b\hspace{-.5pt}}(\mmod H)$, sending a module $M$ to the stalk complex $M[0]$. By Theorem 7.2(1) in \cite{BLP}, this induces a convex valued translation quiver embedding of $\GaH$ in $\Ga_{\hspace{-.8pt}D^{\hspace{.4pt}b\hspace{-.5pt}}(\mmod H)}$. Next, in view of Theorem 7.2 in \cite{BLP}, we may apply Happel's argument in \cite[(4.5)]{Hap} to show that $\Ga_{\hspace{-.8pt}D^{\hspace{.4pt}b\hspace{-.5pt}}({\rm mod}H)}\cong \mathbb{Z} \Da_H$. The proof of the theorem is completed.

\subsection{\sc Hammocks.} Since $\GaH$ is finite and contains no oriented cycle, we may study hammocks in $\GaH$; see (\ref{hammock_defn}).
Given a vertex $(s,k)$ in $\ZQH$, by Proposition \ref{ZD_sec}, $\ZQH$ contains an $(s,k)$-source section $\Da_{s,k}$ generated by the sectional successors of $(s,k)$. And by Lemma \ref{s_sec_unicity}, $\Da_{s,k}$ is a leftmost section in ${\rm Suc}(s,k)$, the full valued translation subquiver of $\ZQH$ generated by the successors of $(s,k)$. By Proposition \ref{stab_af}, we may introduce the following definition.

\begin{defn}\label{hamm_func}

Let $H$ be a hereditary artin algebra with a Dynkin ext-quiver $\QH$. Consider a vertex $k\in \QH$ and the $(0,k)$-source section $\Da$ in $\ZQH$. The {\it extended hammock function} $h_k$ associated with $k\in \QH$ is the unique additive function on ${\rm Suc}(0,k)$ such, for any vertex $(r,i) \in \Da$, that

\begin{enumerate}[$(1)$]

\item $h\hspace{-.5pt}_k(r, i)=1$ if $(r, i)=(0,k);$

\item $h_k(r, i)=d'_1 \cdots d'_t$ if $\Da$ contains a path \vspace{-3.5pt} \hspace{-6pt} $\xymatrixcolsep{14pt}\xymatrix{(0,k) = (r_0, i_0) \ar[r] & (r_1, i_1) \ar[r] & \cdots \ar[r] & (r_t, i_t) = (r, i),\hspace{-2pt} }$ where the arrow $\hspace{-3pt} \xymatrixcolsep{14pt}\xymatrix{(r_{j-1}, i_{j-1})\ar[r] & (r_j, i_j)}\hspace{-3pt}$ has valuation $(d_j, d'_j)$, for $j=1, \ldots, t$.

\end{enumerate}\end{defn}

The following statement says in particular that an extended hammock function is indeed an extension of Brenner's canonical hammock function; see (\ref{hammock_defn}).

\begin{lem}\label{ham_fun_prpt}

Let $H$ be a hereditary artin algebra with a Dynkin ext-quiver $\QH$. Consider the projective modules $P_i$ and  the simple modules $S_i$ in $\GaH$ associated with $i\in \QH$. Let $h_k$ be the extended hammock function associated with some $k\in \QH$. If $(r,i)\in {\rm Suc}(0,k)$ such that $\tau^{-r}\hspace{-1.5pt}P_i\in \GaH,$ then
$$h_k(s,j)=\ell_{S_k}(\tau^{-s}\hspace{-1.5pt}P_j),$$ for any predecessor $(s,j)$ of $(r,i)$ in ${\rm Suc}(0,k).$

\end{lem}

\vspace{-4pt}

\noindent{\it Proof.} \vspace{.5pt} By Proposition \ref{H_alg_Dynkin}, we have a canonical embedding $\varphi: \GaH \to \ZQH: \tau^{-r}\hspace{-1.5pt}P_i \mapsto (r,i)$. Recall that $\varphi$ has a convex image $\Ga$ in $\ZQH$, and maps the $P_k$-source section $\Sa$ in $\GaH$ onto the $(0,k)$-source section $\Da$ in $\ZQH$. \vspace{.5pt} Assume that $(r,i)\in {\rm Suc}(0,k)$ such that $\tau^{-r}\hspace{-1.5pt}P_i\in \GaH$, that is, $(r,i)\in \Ga.$
Since $\Ga$ is convex in $\ZQH$, we see that $\tau^{-s}\hspace{-1.5pt}P_j\in \GaH$, for any predecessor $(s,j)$ of $(r,i)$ in ${\rm Suc}(0,k).$ Therefore, it suffices to show that $h_k(r,i)=\ell_{S_k}(\tau^{-r}\hspace{-1.5pt}P_i).$

By Proposition \ref{uni_sec_path}, the paths in $\ZQH$ from $(0, k)$ to $(r, i)$ have the same length, written as $l_{r,i}$. Since $h_k(0,k)=1=\ell_{S_k}(P_k)$; see (\ref{pi_endo}), we may assume that $l_{r,i}>0$. Suppose first that $(r,i)\in \Da.$ Then, $\Da$ contains a non-trivial sectional path $\hspace{-4pt}\xymatrixcolsep{18pt}\xymatrix{(0,k) = (r_0, i_0) \ar[r] & (r_1, i_1) \ar[r] & \cdots \ar[r] & (r_t, i_t) = (r, i),}$ where $(r_{j-1}, i_{j-1})\to (r_j, i_j)$ has valuation, say $(d_j, d'_j)$, for $j=1, \ldots, t$. \vspace{-3pt} Since $\varphi(\Sa)=\Da$, we have a sectional path $\hspace{-4pt}\xymatrixcolsep{18pt}\xymatrix{P_k=\tau^{-r_0}\!P_{i_0} \ar[r] & \tau^{-r_1}\!P_{i_1} \ar[r] & \cdots \ar[r] & \tau^{-r_t}\!P_{i_t}=\tau^{-r}P_i}$ \vspace{-1pt} in $\Sa$, where the arrow $\tau^{-r_{j-1}}\!P_{i_{j-1}} \to \tau^{-r_j}\!P_{i_j}$ also has valuation $(d_j, d'_j)$, for $j=1, \ldots, t$. In view of Lemma \ref{sp_mul} and Definition \ref{hamm_func}, we see that $\ell_S(\tau^{-r}\hspace{-1pt}P_i)=d_1'\cdots d_t'= h_i(r,i).$

\vspace{.5pt}

Assume now that $(r,i)\notin\Da$. Since $\Da$ is a leftmost section in ${\rm Suc}(0,k);$ see (\ref{s_sec_unicity}), we see that $(r-1, i)\in {\rm Suc}(0,k).$ Let $(s_j, k_j)\to (r, i)$ with valuation $(n_j, n_j')$, $j=1, \ldots, p$, be the arrows in $\ZQH$ ending with $(r,i)$. Since $\Ga$ is convex in $\ZQH$, we see that $(r-1, i), (s_j, k_j)\in {\rm Suc}(0,k) \hspace{.4pt} \cap \Ga$. Since $l_{r-1, i}<l_{s_j, k_j}<l_{r, i}$ for $1\le j\le p$, by the induction hypothesis, $h_k(r-1,i)=\ell_{S_k}(\tau^{1-r}\hspace{-1pt}P_i)$ and $h_k(s_j,k_j)=\ell_{S_k}(\tau^{-s_j}\hspace{-1pt}P_{k_j})$ for $1\le j\le p$. Since $\tau^{-s_j}\hspace{-1.2pt} P_{k_j}\to \tau^{-r}\hspace{-1.2pt} P_i$ and $(s_j, k_j)\to (r,i)$ have the same valuation $(n_j, n'_j)$, for $j=1, \ldots, p$, we deduce from
Proposition \ref{Irr_mor} an almost split sequence \vspace{-5pt}
$$\xymatrix{0\ar[r] & \tau^{1-r}\hspace{-1.2pt}P_i \ar[r] & \oplus_{j=1}^p (\tau^{-s_j}\hspace{-1.5pt} P_{k_j})^{n_j'}\ar[r] & \tau^{-r}\hspace{-1.2pt} P_i\ar[r] & 0}\vspace{-3pt}$$ in $\mmod H$. Since $\ell_{S_k}$ and $h_k$ are additive, we obtain \vspace{-0pt}
$$
\textstyle
\ell_{S_k}(\tau^{-r}\hspace{-2pt} P_i)=\sum_{j=1}^p n_j' \hspace{.5pt}
\ell_{S_k}(\tau^{-s_j}\hspace{-1.2pt}P_{k_j})- \ell_{S_k}(\tau^{1-r}\hspace{-2pt}P_i)=\sum_{j=1}^p n_j' \hspace{.5pt} h_k(s_j, k_j)- h_k(r-1,i)=h_k(r,i).\vspace{-1pt}$$
The proof of the lemma is completed.

\vspace{3pt}

Given a vertex $k$ in $\QH$, we shall write ${\mathcal H}_k$ for the hammock ${\mathcal H}_{S_k}$, where $S_k$ is the simple module in $\GaH$  associated with $k$.

\begin{lem}\label{ham_fun_prpt_2}

Let $H$ be a hereditary artin algebra with a Dynkin ext-quiver $\QH$. Consider  the projective modules $P_i$ and the injective modules $I_i$ in $\GaH$ associated with $i\in \QH$. Let $\mathcal{H}_k$ be the hammock and $h_k$ the extended hammock function associated with some $k\in \QH$. Assume that $(r,i)\in {\rm Suc}(0,k)$ such that $h_k(s, j)\ge 0$ for all proper predecessors $(s, j)$ of $(r,i)$ in ${\rm Suc}(0,k)$.

\begin{enumerate}[$(1)$]

\vspace{-1pt}

\item If $h_k(r,i)>0$, then $\tau^{-r}\hspace{-1.5pt} P_i \in \mathcal{H}_k$.

\vspace{.8pt}

\item If $(r-1, i)\in {\rm Suc}(0,k)$ such that $\tau^{1-r}\hspace{-1.5pt}P_i=I_l$ for some vertex $l\in \QH$, then $h_k(r,i)=-1$ in case $l=k$, and otherwise, $h_k(r,i)=0$.

\end{enumerate}\end{lem}

\vspace{-4pt}

\noindent{\it Proof.} Consider the canonical \vspace{.5pt} embedding $\varphi: \GaH \to \ZQH,$ which has a convex image $\Ga$ in $\ZQH$. Let $\Sa$ be the $P_k$-source section in $\GaH$, and $\Da$ the $(0,k)$-source section in $\ZQH$. Furthermore, let $(s_j, k_j)\to (r, i)$ with valuation $(n_j, n_j')$, $j=1, \ldots, p$, be the arrows in $\ZQH$ ending with $(r,i)$. Given a vertex $(s,j)\in {\rm Suc}(0,k)$, denote by $l_{s,j}$ the length of paths in $\ZQH$ from $(0,k)$ to $(s,j)$. In case $l_{r,i}\le 1$, we have $(r,i) \in \Da$, and consequently,  $\tau^{-r}\hspace{-1.5pt} P_i \in \mathcal{H}_k$; see (\ref{H_alg_Dynkin}) and $(r-1, i)\not\in {\rm Suc}(0,k)$; see (\ref{s_sec_unicity}). Thus, the two statements hold in this case. Consider the case where $l_{r,i}>1$.

Suppose first that $(r-1, i)\in {\rm Suc}(0,k)$ such that $\tau^{1-r}\hspace{-1.5pt}P_i=I_l$, for some vertex $l\in \QH$. Writing $S_l={\rm soc}I_l$, we have a short exact sequence  \vspace{-4pt}
$$\xymatrix{0\ar[r] & S_l \ar[r] & \tau^{1-r}\hspace{-1.5pt}P_i \ar[r]^-g& M\ar[r] & 0,} \vspace{-3pt}$$ where $g$ is minimal left almost split. By the assumption on $(r,i)$ stated in the lemma, $h_k(s_j, k_j)\ge 0$, for all $1\le j\le p$. We may assume that there exists some $0\le t\le p$ such that $(s_j, k_j)\in \Ga$ if and only if $1\le j\le t.$ Given $t<j\le p$, in view of the induction hypothesis, we see that $h_k(s_j, k_j)=0$. Since $\tau^{1-r}P_i\to \tau^{-s_j}\hspace{-1.5pt}P_{k_j}$ and $(r-1, i)\to (s_j, k_j)$ have the same valuation $(n_j', n_j)$ for $1\le j\le t$, \vspace{-1pt} by Proposition \ref{Irr_mor}(1), $M\cong \oplus_{j=1}^t (\tau^{-s_j}\hspace{-1.5pt}P_{k_j})^{n_j'}$. \vspace{-3pt}  Since $h_k$ and $\ell_{S_k}$ are additive, by Lemma \ref{ham_fun_prpt}, we obtain \vspace{-8pt}
$$\textstyle h_k(r,i)= \sum_{j=1}^t n_j' \hspace{.5pt}  h_k(s_j, k_j) - h(1-r, i)  =  \sum_{j=1}^t n_j' \hspace{.5pt}  \ell_{S_k}(\tau^{-s_j}\hspace{-1.5pt}P_{k_j}) - \ell_{S_k}(\tau^{1-r}\hspace{-1.5pt}P_i) =  \ell_{S_k}(M) - \ell_{S_k}(\tau^{1-r}\hspace{-1.5pt}P_i).\vspace{-2pt}$$

If $l=k$, then $\ell_{S_k}(\ell_{S_k}(\tau^{1-r}\hspace{-1.5pt}P_i)=1$; see (\ref{pi_endo}). Thus $\ell_{S_k}(M)=0$, and hence, $h_k(r,i)=-1$. If $l\ne k$, then $\ell_{S_k}(S_l)=0$. Thus $\ell_{S_k}(M) = \ell_{S_k}(\tau^{1-r}\hspace{-1.5pt}P_i)$,  hence, $h_k(r,i)=0$. Statement (2) holds in this case.

Suppose now that $h_k(r, i)>0$. If $(r,i)\in \Da$, then $\tau^{-r}\!P_i\in \Sa \subseteq \mathcal{H}_k$; see (\ref{hammock_prop}). Otherwise, since $\Da$ is a leftmost section of ${\rm Suc}(0,k)$; see (\ref{s_sec_unicity}), $(r-1, i)\in {\rm Suc}(0,k)$ and $h_k(r-1, i)\ge 0$ by the assumption on $(r,i)$. Since $h_k$ is additive and $h_k(r,i)>0$, we may assume that $h_k(s_1, k_1)>0$. By the induction hypothesis, $\tau^{-s_1}\! P_{k_1}\in \mathcal{H}_k$. 
Since $\Ga$ is convex in $\ZQH$, 
we see that $\tau^{1-r}\hspace{-1.5pt}P_i\in \GaH$. Since $h_k(r, i)>0$, as has been shown, $\tau^{1-r}\hspace{-1.5pt}P_i$ is not injective. Then, $\tau^{-r}\hspace{-1.5pt}P_i \in \GaH$ with $h_k(r, i)>0$. By Lemma \ref{ham_fun_prpt}, $\tau^{-r}\hspace{-1.5pt}P_i \in\mathcal{H}_k$. Statement (1) holds in this case. The proof of the lemma is completed.

\vspace{3pt}

The following statement describes all the hammocks in $\GaH$ and tells us an injective module lies in the $\tau$-orbit of  which projective module.

\begin{thm} \label{hamm_neg_section}

Let $H$ be a hereditary artin algebra with a Dynkin ext-quiver $\QH$. Consider  the projective modules $P_i$ and the injective modules $I_i$ in $\GaH$ associated with vertices $i\in \QH$.
Let $\mathcal{H}_k$ be the hammock and $h_k$ the extended hammock function associated with some vertex $k\in\QH$. Then,

\begin{enumerate}[$(1)$]

\vspace{-.5pt}

\item $I_k=\tau^{1-s_k}\hspace{-1.5pt}P_{i_k},$ where $(s_k, i_k)\in {\rm Suc}(0,k)$ \vspace{.5pt} such that $h_k(s_k, i_k)=-1,$ and $h_k(s,j) \ge 0$ for all proper predecessors $(s,j)$ of $(s_k, i_k)$ in ${\rm Suc}(0,k);$

\vspace{.5pt}

\item ${\mathcal H}_k$ \vspace{1pt} is generated by the modules $\tau^{-r}\!P_i\in \GaH,$ where $(r,i)$ lies in the convex hull of $(0,k)$ and $(s_k-1, i_k)$ in $\ZQH$ and $h_k(r,i)>0.$

\end{enumerate}\end{thm}

\vspace{-4pt}

\noindent{\it Proof.} Consider the canonical \vspace{.5pt} embedding $\varphi: \GaH \to \ZQH,$ which has a convex image $\Ga$ in $\ZQH$. By Proposition \ref{hammock_prop}(1) and Lemma \ref{ham_fun_prpt}, we easily deduce Statement (2) from Statement (1). So, we only prove Statement (1). Suppose first that $(s_k, i_k) \in {\rm Suc}(0,k)$ such that $h_k(s_k, i_k)=-1$ and $h_k(s,j) \ge 0$ for all proper predecessors $(s,j)$ of $(s_k, i_k)$ in ${\rm Suc}(0,k).$ By Proposition \ref{H_alg_Dynkin} and Lemma \ref{ham_fun_prpt}, $(s_k, i_k)$ is not in the $(0,k)$-source section $\Da$ in $\ZQH$. By Lemma \ref{s_sec_unicity}(1), $(s_k-1, i_k)\in {\rm Suc}(0,k)$. Since $h_k$ is additive, we deduce from the assumption on $(s_k,i_k)$ that $h(s_k-1, i_k)>0$. By Lemma \ref{ham_fun_prpt_2}(1), $\tau^{1-s_k}\hspace{-.5pt}P_{i_k}\in \mathcal{H}_k.$ If $\tau^{1-s_k}\hspace{-.5pt}P_{i_k}$ is not injective, then $\tau^{-s_k}\hspace{-.5pt}P_{i_k}\in \GaH$, and hence, $h_k(s_k, i_k)\ge 0$; see (\ref{ham_fun_prpt}), a contradiction. Thus, $\tau^{1-s_k}\hspace{-.5pt}P_{i_k}$ is injective, and by Lemma \ref{ham_fun_prpt_2}(2), $\tau^{1-s_k}\hspace{-.5pt}P_{i_k}=I_k$.

Suppose conversely that $I_k=\tau^{-r_k}\hspace{-1.5pt}P_{i_k}$, where $(r_k, i_k)\in \Z\times (\QH)_0$. Since $I_k$ is a succesor of $P_k$ in $\GaH$ by Proposition \ref{hammock_prop}(1), $(r_k, i_k) \in {\rm Suc}(0,k)\cap \Ga$.
Put $s_k=r_k+1$. We claim $h_k(s, j)\ge 0$, for any proper predecessor $(s, j)$ of $(s_k, i_k)$ in ${\rm Suc}(0,k)$. Indeed, if $(s,j)$ is a predecessor of $(r_k, i_k)$ in ${\rm Suc}(0,k)$, then $h_k(s, j)\ge 0$ by Lemma \ref{ham_fun_prpt}. Otherwise, $(s-1,j)$ is a predecessor of $(r_k, i_k)$ in ${\rm Suc}(0,k)$. By Lemma \ref{ham_fun_prpt}, $\tau^{1-s}\hspace{-1.5pt}P_j\in \GaH$. If $\tau^{1-s}\hspace{-1.5pt}P_j$ is not injective, then $\tau^{-s}\hspace{-1.5pt}P_j\in \GaH$, and by Lemma \ref{ham_fun_prpt}, $h_k(s, j)\ge 0$. Suppose that $\tau^{1-s}\hspace{-1.5pt} P_j=I_l$ for some $l\in (\QH)_0.$ Since $(s,j)\ne (s_k, i_k)$, we have $(s-1, j)\ne (r_k, i_k)$. Thus, $l\ne k$, and by Lemma \ref{ham_fun_prpt_2}(2), $h_k(s,j)=0$. This establishes our claim. In particular, by Lemma \ref{ham_fun_prpt_2}(2), $h_k(s_k, i_k)=-1$.
%
The proof of the theorem is completed.


\vspace{2pt}

Theorem \ref{hamm_neg_section} allows us to determine an injective module lies in which $\tau$-orbit of a projective module. For our later purpose, we shall provide an example for each of the types $\mathbb{A}_n$, $\mathbb{D}_n$ and $\mathbb{E}_6$; see (\ref{Dyn_diag}).

\begin{lem} \label{An}

Let $H$ be a hereditary artin algebra with a Dynkin ext-quiver $\QH$. Consider the projective modules $P_i$ and the injective modules $I_i$ in $\GaH$ associated with $i\in \QH$. If $\OQH=\mathbb{A}_n$ with $n\ge 1$, then $I_1=\tau^{-r_{1,n}}\hspace{-1.5pt}P_n$, where $r_{1,n}$ is the number of arrows in the reduced walk in $\QH$ from $1$ to $n$.

\end{lem}

\vspace{-5pt}

\noindent {\it Proof.} Let $\OQH=\mathbb{A}_n$ with $n\ge 1$. Given $i\in (\QH)_0$, write $r_{1,i}=\mathfrak{a}(1,i)$, the number of arrows in the reduced walk in $\QH$ from $1$ to $i$, that is the number of inverse arrows in the reduced walk in $Q^{\rm op}_H$ from $1$ to $i$. In particular, $r_{1,1}=0$. If $n=1$, then $I_1=\tau^{-r_{1,1}}\hspace{-1pt}P_1$. Suppose that $n\ge 2$. By Lemma \ref{sec_paths}(1) and Proposition \ref{ZD_sec}, \vspace{-1pt} the $(0,1)$-source section in $\ZQH$ is $\hspace{-2pt}\xymatrixcolsep{18pt}\xymatrix{(r_{1,1}, 1) \ar[r] & (r_{1,2}, 2) \ar[r] & \cdots \ar[r] & (r_{1,n},  n).}$ \vspace{-1pt} Consider the extended hammock function $h_1$ associated with the vertex $1$. Since $\ZQH$ is trivially valued, by Definition \ref{hamm_func}, $h_1(r_{1,i}, i)=1$ for $i=1, \ldots, n$. Since $h_1$ is additive, we can depict its valuation on the convex hull of $(0,1)$ and $(r_{1,n}+1, n)$ in $\ZQH$ as follows$\,:$\vspace{-1pt}
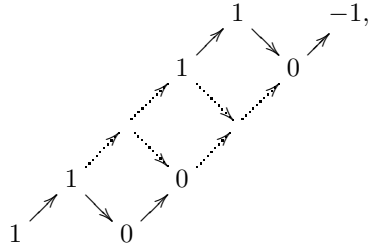
$$\xymatrixcolsep{15pt}\xymatrixrowsep{15pt}\xymatrix@!=6pt{
&&&&1\ar[dr]&&-1,&&\\
&&&1\ar@{.>}[dr]\ar[ur]&&0\ar[ur]&&\\
&&\ar@{.>}[ur]\ar@{.>}[dr]&&\ar@{.>}[ur]&&\\
&1\ar@{.>}[ur]\ar[dr]&&0\ar@{.>}[ur]&&\\
1\ar[ur]&& 0\ar[ur]
} \vspace{-2pt} $$ where the source is $(0,1)$ and the sink is $(r_{1,n}+1, n)$. By Theorem \ref{hamm_neg_section}(1), $I_1=\tau^{-r_{1,n}}\hspace{-1.5pt}P_n$. The proof of the lemma is completed.

\begin{lem}\label{Dn}

Let $H$ be a hereditary artin algebra with a Dynkin ext-quiver $\QH$. Consider the projective modules $P_i$ and the injective modules $I_i$ in $\GaH$ associated with $i\in \QH$. If $\OQH=\mathbb{D}_n$ with $n\ge 4$, then $I_1=\tau^{2-n}\hspace{-1pt} P_1$ in case $n$ is even$\,;$ and $I_1=\tau^{3-n-r_{1,2}}\hspace{-1pt}P_2$ in case $n$ is odd, where $r_{1,2}$ is the number of arrows in the reduced walk in $\QH$ from $1$ to $2$.

\end{lem}

\vspace{-5pt}

\noindent {\it Proof.} Assume that $\OQH=\mathbb{D}_n$ with $n\ge 4$. For $i\in (\QH)_0$, write $r_{1,i}=\mathfrak{a}^+(1,i)$, the number of arrows in the reduced walk in $\QH$, that is the number of inverse arrows in the reduced walk in $Q^{\rm op}_H$, from $1$ to $i$. In particular, $r_{1,1}=0$. By Lemma \ref{sec_paths}(1) and Proposition \ref{ZD_sec}, the $(0,1)$-source section in $\ZQH$ is
\vspace{-2pt}
$$\xymatrixrowsep{15pt}\xymatrixcolsep{16pt}
\xymatrix{
(r_{1,1}, 1)\ar[r] & (r_{1,3}, 3) \ar[r] \ar@<-0.4ex>[d]  & (r_{1,4}, 4) \ar[r] &\cdots \ar[r] & (r_{1,n}, n).\\
&(r_{1,2}, 2)&&
}\vspace{-4pt}$$

Consider the extended hammock function $h_1$ associated with the vertex $1$. Since $\ZQH$ is trivially valued, $h_1(r_{1,i}, i)=1$ for $i=1, 2, \ldots, n$. Consider first the case where $n$ is even. Since $h_1$ is additive, its valuation on the convex hull of $(0,1)$ and $(n-1,1)$ in $\ZQH$ can be depicted as follows$:$ \vspace{-0pt}
$$\xymatrixcolsep{15pt}\xymatrixrowsep{15pt}\xymatrix@!=5pt{
&&&&&&1\ar[dr]&&0\ar[dr]\\
&&&&&1\ar[dr]\ar[ur]&&1\ar[dr]\ar[ur]&&0\ar@{.>}[dr]\\
&&&&1\ar@{.>}[dr]\ar[ur]&&1\ar@{.>}[dr]\ar[ur]&&1\ar@{.>}[dr]\ar[ur]&&\ar@{.>}[dr]\\
&&&\ar@{.>}[ur]\ar@{.>}[dr]&&\ar@{.>}[ur]\ar@{.>}[dr]&&\ar@{.>}[ur]\ar@{.>}[dr]&&
\ar@{.>}[dr]\ar@{.>}[ur]&&0\ar[dr]\\
&&1\ar@{.>}[ur]\ar[dr]&&1\ar[dr]\ar@{.>}[ur]&&1\ar@{.>}[dr]\ar@{.>}[ur]&& \ar@{.>}[ur]\ar@{.>}[dr] &&1\ar[dr]\ar[ur]&&0\ar[dr]\\
&1\ar[r]\ar[ur]\ar[dr] & 1\ar[r] & 1\ar[dr]\ar[ur]\ar[r] & 0\ar[r] & 1\ar[ur]\ar@{.>}[dr]\ar[r] &1 \ar@{.>}[r] & \ar@{.>}[ur]\ar@{.>}[r] \ar@{.>}[dr]& \ar@{.>}[r]& 1\ar[ur] \ar[r] \ar[dr] & 1\ar[r] & 1 \ar[ur]\ar[dr]\ar[r]&0\ar[r]&0 \ar[dr]\\
1\ar[ur]&& 0\ar[ur] && 1\ar[ur]&&\ar@{.>}[ur] && 1\ar[ur] && 0 \ar[ur] && 1\ar[ur]&&-1,
}$$ where the source is $(0,1)$ and the sink is $(n-1, 1)$. By Theorem \ref{hamm_neg_section}(1), $I_1=\tau^{2-n}\hspace{-1pt} P_1$. Suppose next that $n$ is odd. Since $h_1$ additive, the valuation of $h_1$ on the convex hull of $(0,1)$ and $(r_2+n-2, 2)$ in $\ZQH$ can be depicted as follows$\,:$ \vspace{-4pt}
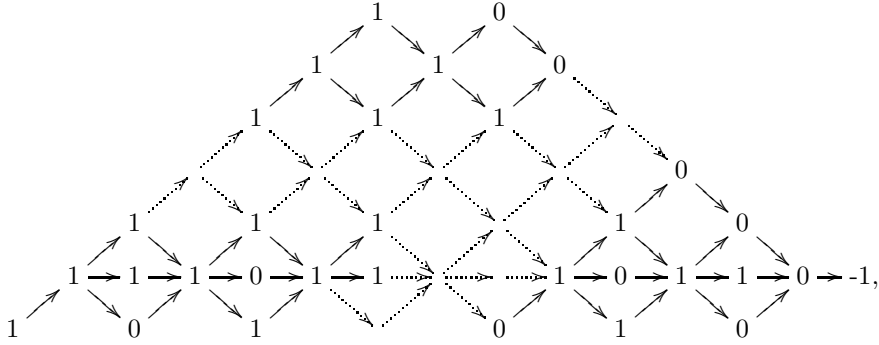
$$\xymatrixcolsep{18pt}\xymatrixrowsep{15pt}\xymatrix@!=5pt{
&&&&&&1\ar[dr]&&0\ar[dr]\\
&&&&&1\ar[dr]\ar[ur]&&1\ar[dr]\ar[ur]&&0\ar@{.>}[dr]\\
&&&&1\ar@{.>}[dr]\ar[ur]&&1\ar@{.>}[dr]\ar[ur]&&1\ar@{.>}[dr]\ar[ur]&& \ar@{.>}[dr]\\
&&&\ar@{.>}[ur]\ar@{.>}[dr]&&\ar@{.>}[ur]\ar@{.>}[dr]&&\ar@{.>}[ur]\ar@{.>}[dr]&&
\ar@{.>}[dr]\ar@{.>}[ur]&&0\ar[dr]\\
&&1\ar@{.>}[ur]\ar[dr]&&1\ar[dr]\ar@{.>}[ur]&&1\ar@{.>}[dr]\ar@{.>}[ur]&& \ar@{.>}[ur] \ar@{.>}[dr] &&1\ar[dr]\ar[ur]&&0\ar[dr]\\
&1\ar[r]\ar[ur]\ar[dr] & 1 \ar[r] & 1\ar[dr]\ar[ur]\ar[r] & 0\ar[r] & 1\ar[ur]\ar@{.>}[dr]\ar[r] &1\ar@{.>}[r] & \ar@{.>}[ur]\ar@{.>}[r] \ar@{.>}[dr]& \ar@{.>}[r]& 1\ar[ur] \ar[r] \ar[dr] & 0 \ar[r] & 1 \ar[ur]\ar[dr]\ar[r]&1\ar[r]& 0 \ar[r] & \mbox{-}1, \\
1\ar[ur]&& 0\ar[ur] && 1\ar[ur]&&\ar@{.>}[ur] && 0\ar[ur] && 1 \ar[ur] && 0\ar[ur]}$$ where the source is $(0,1)$ and the sink is $(r_{1,2}+n-2, 2).$ By Theorem \ref{hamm_neg_section}(1), $I_1=\tau^{-r_{1,2}-n+3}\hspace{-1pt}P_2.$ The proof of the lemma is completed.

\begin{lem} \label{E6}

Let $H$ be a hereditary artin algebra with a Dynkin ext-quiver $\QH$. Consider the projective modules $P_i$ and the injective modules $I_i$ in $\GaH$ associated with $i\in \QH$. If $\OQH=\mathbb{E}_6$, then $I_1=\tau^{-r_{1,6}-3}P_6$, where $r_{1,6}$ is the number of arrows in the reduced walk in $\QH$ from $1$ to $6$.

\end{lem}

\vspace{-4pt}

\noindent{\it Proof.} Assume that $\OQH=\mathbb{E}_6$. Given $i\in (\QH)_0$, let $r_{1,i}$ be the number of arrows in the reduced walk in $\QH$, that is, the number of inverse arrows in the reduced walk in $Q^{\rm op}_H$, from $1$ to $i$. In particular, $r_{1,1}=0$. By Lemma \ref{sec_paths}(1) and Proposition \ref{ZD_sec}, the $(0,1)$-sourced section in $\ZQH$ is \vspace{-2pt}
$$\xymatrixrowsep{16pt}\xymatrixcolsep{16pt}
\xymatrix{
(r_{1,1}, 1)\ar[r] & (r_{1,2}, 2) \ar[r]&(r_{1,3}, 3) \ar[r]\ar@<-0.4ex>[d]&(r_{1,5},5) \ar[r] & (r_{1,6}, 6).\\
&&(r_{1,4}, 4)&&
}$$ \vspace{-8pt}

Consider the extended hammock function $h_1$ associated with the vertex $1$. Since $\ZQH$ is trivially valued, $h_1(r_{1,i}, i)=1$ for $i=1, \ldots, 6$. So, its valuation on the convex hull of $(0,1)$ and $(r_6+4,6)$ in $\ZQH$ can be depicted as follows$:$ \vspace{-4pt}
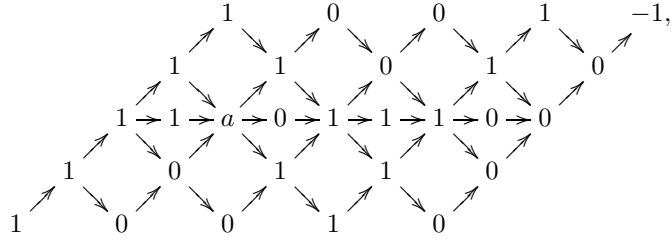
$$\xymatrixcolsep{15pt}\xymatrixrowsep{15pt}\xymatrix@!=5pt{
&&&&1\ar[dr]&&0\ar[dr]&&0\ar[dr]&&1\ar[dr]&&-1, \\
&&&1\ar[dr]\ar[ur]&&1\ar[dr]\ar[ur]&&0\ar[dr]\ar[ur]&&1\ar[dr]\ar[ur]&&0\ar[ur]\\
&&1\ar[r]\ar[dr]\ar[ur]&1\ar[r]&a\ar[ur]\ar[dr]\ar[r]&0\ar[r]&1\ar[dr]\ar[ur]\ar[r]&1\ar[r]&1\ar[dr]\ar[ur]\ar[r]&0\ar[r]&0\ar[ur] \\
&1\ar[dr]\ar[ur] &&0\ar[dr]\ar[ur]&&1\ar[dr]\ar[ur]&&1\ar[dr]\ar[ur]&&0\ar[ur]\\
1\ar[ur]&&0 \ar[ur]&&0\ar[ur]&&1\ar[ur]&&0\ar[ur]
}$$ where the source is $(0,1)$ and the sink is $(r_{1,6}+4, 6).$ By Theorem \ref{hamm_neg_section}(1), $I_1=\tau^{-r_{1,6}-3}P_6$. The proof of the lemma is completed.

\subsection{\sc Pi-permutation and pi-indices} By Proposition \ref{H_alg_Dynkin}, $\GaH$ is finite and convexly embeds in $\ZQH$ in such a way that every module lies in the $\tau$-orbit of a projective module. We need to determine the number of modules and the injective module lying in the $\tau$-orbit of a projective module. For this purpose, we introduce the notions of the pi-permutation and the pi-index function in the following well-known statement; see \cite[(1.8)]{PAu}.

\begin{lem}\label{exq_permut}

Let $H$ be a hereditary artin algebra with a Dynkin ext-quiver $\QH.$ Consider the projective modules and the injective modules $I_i$ in $\GaH$ associated with vertices $i\in \QH$. Then there exists a permutation $\rho$ of $(\QH)_0$, called the \emph{pi-permutation}, and a function
$m: (\QH)_0 \to \mathbb{N}$, called the \emph{pi-index function}, such that $\tau^{-m(i)}\hspace{-.6pt} P_i=I_{\rho(i)}$, for every vertex $i$ in $\QH$.

\end{lem}

The following statement is important in our investigation, where $\mathfrak{a}^+(i,j)$ denotes the number of arrows in the reduced walk in $\QH$ from $i$ to $j$; see (\ref{No_arrows}).

\begin{lem}\label{pi-index}

Let $H$ be a hereditary artin algebra with a Dynkin ext-quiver $\QH$. Consider the pi-permutation $\rho$ and the pi-index function $m$ for $(\QH)_0$. If $i,j$ are vertices in $\OQH$, then $$m(i)-m(j)=\mathfrak{a}^+(\rho(i), \hspace{.3pt} \rho(j))-\mathfrak{a}^+(i,j).
\vspace{.5pt}$$

\end{lem}

\vspace{-4pt}

\noindent{\it Proof.} Consider the projective modules $P_i$ and the injective modules $I_i$ in $\GaH$ associated with $i\in \QH$. By Proposition \ref{H_alg_Dynkin}, we have a canonical embedding $\varphi: \GaH\to \mathbb{Z} Q^{\rm op}_{\hspace{-.8pt}H}$, sending $\tau^{-r}\hspace{-1.5pt}P_i$ to $(r,i).$ Write $r_{p,q}=\mathfrak{a}^+(p,q)$ for $p,q\in (\QH)_0$.
Fix some vertices $i, j\in \QH$. Observe that $r_{i,j}$ is the number of inverse arrows in the reduced walk in $Q^{\rm op}_{\hspace{-.8pt}H}$ from $i$ to $j$. By Lemma \ref{sec_paths}(2), $\ZQH$ contains a sectional path from $(0, i)$ to $(r_{i,j}, j)$, and hence, a sectional path from $(m(j)-r_{i,j}, i)$ to $(m(j), j)$.

On the other hand, since $r_{\rho(i), \hspace{.3pt} \rho(j)}$ is the number of arrows in the reduced walk in $\QH$ from $\rho(i)$ to $\rho(j)$, by Lemma \ref{pic_sec_path}(2), $\GaH$ contains a sectional path from $\tau^{\hspace{.5pt}r_{\rho(i), \hspace{.3pt} \rho(j)}\hspace{-1pt}}I_{\rho(i)}$ to $I_{\rho(j)}$. By definition, $I_{\rho(i)}=\tau^{-m(i)}\hspace{-1.2pt}P_i$ and $I_{\rho(j)}=\tau^{-m(j)}\hspace{-1.2pt}P_j$. Thus, $\GaH$ contains a sectional path from $\tau^{\hspace{.5pt}r_{\rho(i), \hspace{.3pt} \rho(j)}-m(i)\hspace{-1.2pt}} P_i$ to $\tau^{-m(j)}\hspace{-1.2pt}P_j,$ which induces a sectional path in $\mathbb{Z} Q^{\rm op}_{\hspace{-.8pt}H}$ from $(m(i)-r_{\rho(i), \hspace{.3pt} \rho(j)}, i)$ to $(m(j), j)$. By Proposition \ref{ZD_sec}, \vspace{.5pt} both $(m(j)-r_{i,j}, i)$ and $(m(i)-r_{\rho(i), \hspace{.3pt} \rho(j)}, i)$ belong to the $(m(j), j)$-sink section in $\ZQH$, and consequently, $m(i)-r_{\rho(i), \hspace{.3pt} \rho(j)}=m(j)-r_{i,j}$. The proof of the lemma is completed.

\vspace{2pt}

Next, \vspace{.5pt} we shall show that the pi-permutation of $(\QH)_0$ is of order at most two. For doing this, we shall extend it to an automorphism of the ext-graph $\OQH$ in the following statement.

\begin{prop}\label{vq_iso}

Let $H$ be a hereditary artin algebra with a Dynkin ext-quiver $\QH$. Then the pi-permutation $\rho$ of $(\QH)_0$ induces a valued graph automorphism
of $\hspace{2pt}\overline{\hspace{-1.5pt}\QH\hspace{-7pt}} \hspace{7pt},$
called the \emph{pi-automorphism} and denoted again by $\rho.$

\end{prop}

\vspace{-4pt}

\noindent{\it Proof.} Let $\rho$ be the pi-permutation of $(\QH)_0$ and $m$ the pi-index function on $(\QH)_0$. Consider the projective modules $P_i$ and the injective modules $I_i$ in $\GaH$ associated with $i\in (\QH)_0$. Let $\DaH$ be the projective section generated by the $P_i$, and $\SaH$ the injective section generated by the $I_i$, in $\GaH$. We first construct a valued graph isomorphism $\rho^*: \hspace{3pt}\overline{\hspace{-3pt}\Da}_{\hspace{-.5pt}H}\to
\hspace{2pt}\overline{\hspace{-2pt}\Sa}\hspace{-1.5pt}_H,$ sending $P_i$ to $I_{\rho(i)}.$

\vspace{-1pt}

Fix $i,j\in (\QH)_0$. Suppose that $\hspace{3pt}\overline{\hspace{-3pt}\Da}_{\hspace{-.5pt}H}$ contains an edge $\xymatrixcolsep{15pt}\xymatrix{\hspace{-2pt}P_i\ar@{-}[r] & P_j\hspace{-2pt}}\vspace{-1pt}$ with valuation $(d,d')$. \vspace{-3pt}  We claim that $\hspace{2pt}\overline{\hspace{-2pt}\Sa}\hspace{-1.5pt}_H$ contains an edge
$\hspace{-2pt} \xymatrixcolsep{12pt}\xymatrix{I_{\rho(i)} \ar@{-}[r] & I_{\rho(j)}}\hspace{-2pt} \vspace{-2pt} $
with valuation $(d, d')$. First, assume that $\DaH$ contains an arrow $P_i\to P_j$, whose valuation is then $(d,d')$; see (\ref{under_vq_val}). Note that $\tau^{-m(i)}\hspace{-1.5pt}P_i=I_{\rho(i)}$ and $\tau^{-m(j)}\hspace{-1.5pt}P_j=I_{\rho(j)}$. By Lemma \ref{arrow_translate}(2), $\GaH$ contains an arrow $\tau^{-m(j)}\hspace{-1.5pt}P_i\to I_{\rho(j)}$, whose valuation is $(d, d')$; see (\ref{pa_valuation}). So, $m(j)\le m(i)$. \vspace{-1pt} If $m(j)=m(i)$, then $\SaH$ contains an arrow $I_{\rho(i)}\to I_{\rho(j)}$ with valuation $(d, d')$, and hence, $\hspace{2pt}\overline{\hspace{-2pt}\Sa}\hspace{-1.5pt}_H$ contains an edge $\hspace{-2pt} \xymatrixcolsep{12pt}\xymatrix{I_{\rho(i)} \ar@{-}[r] & I_{\rho(j)}}\hspace{-2pt}\vspace{-2pt}$ with valuation $(d, d')$. If $m(j)<m(i)$, then $\GaH$ contains an arrow $I_{\rho(j)} \to \tau^{-m(j)-1}\hspace{-1.5pt}P_i$ with valuation $(d', d)$.
Since $H$ is hereditary, $\tau^{-m(j)-1}\hspace{-1.5pt}P_i$ is injective, and hence, $\tau^{-m(j)-1}\hspace{-1.5pt}P_i=I_{\rho(i)}$. By definition, $m(i)=m(j)+1$, \vspace{-1pt} and hence, $\SaH$ has an arrow $I_{\rho(j)} \to I_{\rho(i)}$ with valuation $(d', d)$. \vspace{-1.5pt} So, $\hspace{2pt}\overline{\hspace{-2pt}\Sa}\hspace{-1.5pt}_H$ contains an edge $\hspace{-2pt} \xymatrixcolsep{12pt}\xymatrix{I_{\rho(j)} \ar@{-}[r] & I_{\rho(i)}}\hspace{-2pt}\vspace{-3pt}$ with valuation $(d', d)$, that is, an edge $\hspace{-2pt} \xymatrixcolsep{12pt}\xymatrix{I_{\rho(i)} \ar@{-}[r] & I_{\rho(j)}}\hspace{-2pt}$ with valuation $(d, d')$.
Next, \vspace{-4pt} assume that $\DaH$ contains an arrow $P_j\to P_i$, whose valuation is then $(d',d)$; see (\ref{under_vq_val}). \vspace{-4.5pt} As has been shown, $\hspace{2pt}\overline{\hspace{-2pt}\Sa}\hspace{-1.5pt}_H$ contains an edge $\hspace{-2pt} \xymatrixcolsep{12pt}\xymatrix{I_{\rho(j)} \ar@{-}[r] & I_{\rho(i)}}\hspace{-2pt}$ with valuation $(d', d)$, that is, an edge $\hspace{-2pt} \xymatrixcolsep{12pt}\xymatrix{I_{\rho(i)} \ar@{-}[r] & I_{\rho(j)}}\hspace{-2pt}$ with valuation $(d, d')$. This establishes our claim.

\vspace{-3pt}

Suppose, on the other hand, that $\hspace{2pt} \overline{\hspace{-2pt}\Sa}\hspace{-1.5pt}_H$ contains an edge $\hspace{-2pt} \xymatrixcolsep{12pt}\xymatrix{I_{\rho(i)} \ar@{-}[r] & I_{\rho(j).}}$ \vspace{-3.5pt} We shall show that $\hspace{2pt}\overline{\hspace{-2pt}\Da}\hspace{-1pt}_H$ contains an edge $\hspace{-2pt} \xymatrixcolsep{12pt}\xymatrix{P_i \ar@{-}[r] & P_j.}$ \vspace{-3.5pt} Note that $P_i=\tau^{\hspace{.5pt}m(i)}\hspace{-.6pt} I_{\rho(i)}$ and $P_j=\tau^{\hspace{.5pt}m(j)}\hspace{-.6pt} I_{\rho(j)}$. Assume first that $\SaH$ contains an arrow $I_{\rho(i)}\to I_{\rho(j)}$. \vspace{.5pt} By Lemma \ref{arrow_translate}(2), $\GaH$ contains an arrow $P_i\to \tau^{m(i)}\! I_{\rho(j)}$. \vspace{-1.5pt} So, $m(i)\le m(j)$. If $m(i)=m(j)$, then $\SaH$ contains an arrow $P_i\to P_j$, and hence, $\hspace{2pt}\overline{\hspace{-2pt}\Da}\hspace{-1pt}_H$ contains an edge $\hspace{-2pt} \xymatrixcolsep{12pt}\xymatrix{P_i \ar@{-}[r] & P_j.}$ \vspace{-2.5pt}
If $m(i)<m(j)$, then $\GaH$ contains an arrow $\tau^{\hspace{.5pt}m(i)+1}\hspace{-1.2pt}I_{\rho(j)}\to P_i$. Since $H$ is hereditary, $\tau^{\hspace{.5pt}m(i)+1}\!I_{\rho(j)}$ is projective. \vspace{-1.5pt} Then $m(j)=m(i)+1$, and hence, $\SaH$ contains an arrow $P_j\to P_i$. So, $\hspace{2pt}\overline{\hspace{-2pt}\Da}\hspace{-1pt}_H$ contains an edge $\hspace{-2pt} \xymatrixcolsep{14pt}\xymatrix{P_i \ar@{-}[r] & P_j.}$ \vspace{-3.5pt} In case $\SaH$ contains an arrow $I_{\rho(j)}\to I_{\rho(i)}$, by a similar argument,
$\hspace{2pt}\overline{\hspace{-2pt}\Da}\hspace{-1pt}_H$ contains an edge $\hspace{-2pt}\xymatrixcolsep{12pt}\xymatrix{P_i \ar@{-}[r] & P_j.}$ \vspace{-1pt} Thus, we have a valued graph isomorphism $\rho^*: \hspace{3pt}\overline{\hspace{-3pt}\Da}_{\hspace{-.5pt}H} \to
\hspace{2pt}\overline{\hspace{-2pt}\Sa}\hspace{-1.5pt}_H,$ sending $P_i$ to $I_{\rho(i)}$.

Finally, \vspace{.5pt} the valued quiver isomorphisms stated in Theorem \ref{pisec} induce a valued graph isomorphism $\OQH^{\hspace{-4pt}\rm op} \to \hspace{2pt}\overline{\hspace{-2pt}\Da\hspace{-.5pt}}_H$ sending $i$ to $P_i$, and a valued graph isomorphism $\hspace{2pt}\overline{\hspace{-2pt}\Sa}\hspace{-1.5pt}_H \to \OQH^{\hspace{-4pt} \rm op}$ sending $I_i$ to $i$. Compo\-sing these two isomorphisms with $\rho^*,$ we obtain a valued graph automorphism $\OQH^{\hspace{-4pt}\rm op} \to \OQH^{\hspace{-4pt}\rm op},$ sending $i$ to $\rho(i)$. This induces a desired valued graph automorphism $\rho: \OQH\to \OQH,$ sending $i$ to $\rho(i)$. The proof of the proposition is completed.

\vspace{3pt}

We collect some easy properties of the pi-automorphism in the following statement.

\begin{lem}\label{rho_property}

Let $H$ be a hereditary artin algebra with a Dynkin ext-quiver $\QH$, and let $\rho$ be the pi-automorphism of $\OQH$.

\begin{enumerate}[$(1)$]

\vspace{-1pt}

\item [$(1)$] If $i$ is a vertex in $\QH$ of weight $3$, then $\rho(i)=i$.

\vspace{.5pt}

\item [$(2)$] If $\hspace{-2pt}\xymatrixcolsep{14pt}\xymatrix{i \ar@{-}[r] & j}\hspace{-2pt}\vspace{-3pt}$ is a non-trivially valued edge in $\OQH$, then $\rho(i)=i$ and $\rho(j)=j$.

\vspace{.5pt}

\item [$(3)$] If $\hspace{-2pt}\xymatrixcolsep{14pt}\xymatrix{i\ar@{-}[r] & j \ar@{-}[r] & k }\hspace{-2pt}\vspace{-1pt}$ is a reduced walk in $\OQH$, where $j$ has only two neighbors $i$ and $k$, such that $\rho(i)=i$ and $\rho(j)=j$, then $\rho(k)=k$.

\end{enumerate}\end{lem}

\vspace{-4.5pt}

\noindent{\it Proof.} (1) Being a Dynkin diagram, $\OQH$ has at most one vertex of weight $3$; see (\ref{Dyn_diag}). And by Lemma \ref{iso-wt}, ${\rm w}(i)={\rm w}(\rho(i))$ for any vertex $i$ in $\QH$. Thus, if ${\rm w}(i)=3$, then $\rho(i)=i$.

(2) Let $\hspace{-2pt}\xymatrixcolsep{14pt}\xymatrix{i \ar@{-}[r] & j}\hspace{-2pt}\vspace{-1pt}$ be an edge in $\OQH$ with a non-trivial valuation $(d_{ij}, d_{ji})$. Since $\OQH$ is a Dynkin digram, this is the only non-trivially valued edge, and moreover, $d_{ij}\neq d_{ji}$; see (\ref{Dyn_diag}). Since $\rho$ is a valuded graph automorphism, $\rho(i), \rho(j) \in \{i,j\}$ and $d_{\rho(i) \rho(j)}=d_{ij}$. If $\rho(i)=j$, then $\rho(j)=i$. So, $d_{ji}=d_{\rho(i) \rho(j)}=d_{ij}$, a contradiction. Thus, $\rho(i)=i$ and $\rho(j)=j$.

(3) Let \vspace{-1pt} $\hspace{-2pt}\xymatrixcolsep{14pt}\xymatrix{i \ar@{-}[r] & j \ar@{-}[r] & k }\hspace{-2pt}\vspace{-3pt}$ be a reduced walk in $\OQH$, where $j$ has only two neighbors $i$ and $k$. If
 $\rho(j)=j$, then $\hspace{-2pt}\xymatrixcolsep{14pt}\xymatrix{j \ar@{-}[r] & \rho(k) }\hspace{-2pt}\vspace{-1pt}$ is an edge in $\OQH$, and hence, $\rho(k)\in \{i,k\}$. If $\rho(i)=i$, then $\rho(k)=k$. The proof of the lemma is completed.

\vspace{3pt}

We are ready to determine the order of the pi-automorphism for every Dynkin type. In particular, the pi-automorphism is the identity for all non-simply laced types.

\begin{thm} \label{order_2} Let $H$ be a hereditary artin algebra of Dynkin type, and let $\rho$ be the pi-automorphism $\rho$ of the ext-graph $\OQH$.

\begin{enumerate}[$(1)$]

\vspace{-1pt}

\item If $\OQH$ is of type $\mathbb{A}_1$, $\mathbb G_2$, $\mathbb F_4$, $\mathbb E_7,$ $\mathbb E_8$, $\mathbb{B}_n (n \hspace{-1.5pt} \ge \hspace{-1.5pt} 2)$, $\mathbb{C}_n (n \hspace{-1.5pt} \ge \hspace{-1.5pt}  3)$ or $\mathbb D_n (n \hspace{-1.5pt} \ge \hspace{-1.5pt} 4 \ even)$, then $\rho$ is the identity.

\vspace{1pt}

\item If $\OQH$ is of type $\mathbb A_n (n\ge 2)$, $\mathbb D_n (n \ge 5 \ odd\hspace{.6pt})$ or $ \mathbb E_6$, then $\rho$ is of order $2$.

\end{enumerate} \end{thm}

\vspace{-4pt}

\noindent{\it Proof.} \vspace{.5pt} We may assume that $\OQH$ is a canonical Dynkin diagram as stated in Definition \ref{Dyn_diag}. The theorem holds trivially in case $\OQH=\mathbb{A}_1$. We shall proceed case by case for all other cases.

\vspace{.5pt}

1) Suppose that $\OQH=\mathbb{G}_2$ or $\mathbb{B}_2$. Then it contains only one edge, which is non-trivially valued. By Lemma \ref{rho_property}(2), $\rho$ is the identity.

\vspace{.5pt}

2) Suppose that $\OQH=\mathbb{F}_4$. By Lemma \ref{rho_property}(2), $\rho(2)=2$ and  $\rho(3)=3$. And by Lemma \ref{rho_property}(3), $\rho(4)=4$. Thus, $\rho(1)=1$. Therefore, $\rho$ is the identity.

\vspace{.5pt}

3) Suppose that $\OQH$ is $\mathbb{E}_7$ or $\mathbb{E}_8$. Since $3$ is the only vertex of weight $3$, by Lemma \ref{rho_property}(1), $\rho(3)=3$. Thus, $\{\rho(2),\rho(4),\rho(5)\}\subseteq \{2,4,5\}$.
Since ${\rm w}(4)=1$, we have $\rho(4)=4$. If $\rho(2)=5$, then $\rho(1)=6$, a contradiction to Lemma \ref{iso-wt} for $\rho(6)=2$. Thus, $\rho(2)=2$, and consequently,  $\rho(1)=1$ and $\rho(5)=5$. Since each of the vertices $5,6,7$ has at most two neighbors, we deduce from \ref{rho_property}(3) that $\rho(i)=i$ for all vertices $i>5.$ So, $\rho$ is the identity.

\vspace{1pt}

4) Suppose that $\OQH=\mathbb{B}_n$ or $\mathbb{C}_n$ with $n\ge 3$. By Lemma \ref{rho_property}(2), $\rho(1)=1$ and $\rho(2)=2$. Since every vertex $i$ with $1<i<n$ has exactly two neighbors, we easily deduce from Lemma \ref{rho_property}(3) that $\rho(i)=i$ for all $3\le i\le n$. So, $\rho$ is the identity.

\vspace{1pt}

5) Suppose that $\OQH=\mathbb{D}_n$ with $n \ge 4$ even. By Lemma \ref{Dn}, $I_1=\tau^{2-n}\hspace{-1pt} P_1$. Using the same argument there, we can show that $I_2=\tau^{2-n}\hspace{-1pt} P_2$. Thus, $\rho(1)=1$ and $\rho(2)=2$. Since $3$ is of weight $3$, by Lemma \ref{rho_property}(1), $\rho(3)=3$. Then, $\rho(4)=4$. Since every vertex $i$ with $4\le i<n$ has exactly two neighbors, by Lemma \ref{rho_property}(3), $\rho(i)=i$ for all $4< i\le n$. So, $\rho$ is the identity.

\vspace{.5pt}

6) Suppose that $\OQH=\mathbb A_n$ with $n\ge 2.$ By Lemma \ref{An}, $I_1=\tau^{-r_{1,n}}P_n$ with $r_{1,n}\ge 0$. So, $\rho(n)=1$. Since every vertex $i$ with $1<i<n$ has two neighbors, $\rho(1)=n$. Since $n-1$ is the only neighbor of $n$, we have $\rho(2)=n-1.$ Similarly, $\rho(i)=n+1-i$ for all $1\le i\le n.$ Hence, $\rho$ is of order $2$.

\vspace{.5pt}

7) Suppose that $\OQH=\mathbb D_n$ with $n \ge 5$ odd. by Lemma \ref{rho_property}(1), $\rho(3)=3$. Since  $\tau^{-s}P_2=I_1$ for some $s\ge 0$; see (\ref{Dn}), $\rho(2)=1$. Since the vertex $4$ has two neighbors, $\rho(2)=1$, and hence, $\rho(4)=4$. Since every vertex $i$ with $3< i<n$ has only two neighbors, by Lemma \ref{rho_property}(3), $\rho(i)=i$ for $5\le i\le n$. So, $\rho$ is of order $2$.

\vspace{.5pt}

8) Suppose $\OQH=\mathbb E_6.$ By Lemma \ref{rho_property}(1), $\rho(3)=3$. Since the vertex $4$ has only one neighbors, $\rho(4)=4.$ By Lemma \ref{E6}, $\tau^{-t}\hspace{-1.2pt}P_6=I_1$ for some $t\ge 0$. So, $\rho(6) = 1$. Then $\rho(1)=6$. Now, $\rho(2)=5$ and $\rho(5)=2$. So $\rho$ is of order $2$. The proof of the theorem is completed.

\subsection{\sc Coxeter order} In case $H$ is of finite representation type, it is well-known that the Coxeter transformation $C_{\hspace{-.5pt}H}$ of $K_0(\mmod H)$ is of finite order; see \cite[(4.1)]{PAu}. We call this finite order the {\it Coxeter order} and write as $|C_{\hspace{-.5pt}H}|$. The following statement reinforces this fact in particular,.


\begin{prop} \label{order_c}

Let $H$ be a hereditary artin algebra with a Dynkin ext-quiver $\QH$. Consider the Auslander-Reiten translation $\tau_{\hspace{-.5pt}_D}$ of $\Ga_{\hspace{-.8pt}D^{\hspace{.4pt}b\hspace{-.5pt}}(\mmod H)},$ the pi-permutation $\rho$ and the pi-index function $m$ for $(\hspace{1.5pt}\overline{\hspace{-.8pt}Q\hspace{-.3pt}}\hspace{-.3pt}_H)_0$. Then, the Coxeter transformation $C_H$ for $H$ is of finite order $|C_H|$ such, for any $i\in \QH,$ that \vspace{-1pt} $$|C_{\hspace{-.5pt}H}|= m(i)+m(\rho(i))+2 \mbox { \, and \,}
\tau_{\hspace{-.5pt}_D}^{-|C_{\hspace{-.5pt}H}|}P_i= P_i[2].$$

\end{prop}

\vspace{-4pt}

\noindent{\it Proof.} First, we claim that $a\!:\hspace{1pt}=m_i+m_{\rho(i)}$ is a constant, for all $i\in (\QH)_0$. Fix $i, j\in (\QH)_0$. By Lemma \ref{pi-index}, $m(i)=m(j)-r_{i,j}+r_{\rho(i), \hspace{.3pt} \rho(j)},$\vspace{.5pt} where $r_{p,q}$ is the the number of arrows in the reduced walk in $\QH$ from $p$ to $q$. Since $\rho^2={\rm id}$ by Theorem \ref{order_2}, $m(\rho(i))=m(\rho(j))-r_{\rho(i),\hspace{.3pt} \rho(j)}+r_{i,j}.$ Thus, $m(i)+m(\rho(i))= m(j)+m(\rho(j)).$ This establishes our claim.

\vspace{.5pt}

Now, since $\tau^{-m(i)}\hspace{-1.3pt} P_i=I_{\rho(i)}$ by definition and $\tau_{\hspace{-.5pt}_D}^{-}\hspace{-1pt} I_{\rho(i)} = P_{\rho(i)}[1]$; see \cite[(4.3)]{Hap} and \cite[(7.2)]{BLP}, we obtain $\tau_{\hspace{-.5pt}_D}^{-(m(i)+1)}P_i=\tau_{\hspace{-.5pt}_D}^{-}(\tau^{-(m(i)}\hspace{-1pt}P_i)=P_{\rho(i)}[1]$.
This yields \vspace{-3pt}
$$\tau_{\hspace{-.5pt}_D}^{-(a+2)} P_i =\tau_{\hspace{-.5pt}_D}^{-(m(\rho(i))+1)}(\tau_{\hspace{-.5pt}_D}^{-(m(i)+1)}P_i)=\tau_{\hspace{-.5pt}_D}^{-(m(\rho(i))+1)}(P_{\rho(i)}[1])=
P_{\rho^2(i)}[1][1]=P_i[2]. \vspace{-1pt} $$

It remains to show that $C_H$ is of order $a+2$. Given $i\in (\QH)_0$, from Proposition \ref{Cox_tau} we see that
$$C_{\hspace{-.5pt}H}^{-(a+2)}\hspace{-1pt}(\undim P_i)= \undim \hspace{.3pt} \tau_{\hspace{-.5pt}_D}^{-(a+2)}\hspace{-1pt} P_i = \undim  P_i[2] = (-1)^2\undim  P_i = \undim  P_i. \vspace{-2pt}$$ Since $\{\undim P_1, \ldots, \undim P_n\}$ is a basis for $K_0(\mmod H)$, $C_{\hspace{-.5pt}H}^{-(a+2)}=\id$. Consider now an integer $t$ with $1\le t\le a+1=m(1)+m(\rho(1))+1$. Suppose first that $t\le m(1)$. Then, $\tau_{\hspace{-.5pt}_D}^{-t}\hspace{-1pt}P_1 =\tau^{-t}\hspace{-1pt}P_1 \in \GaH$ with $\tau^{-t}\hspace{-1.3pt}P_1\not\cong P_1$. Since the modules in $\GaH$ are determined by their composition factors; see \cite[(IX.2.3)]{ARS}, $\undim  \tau^{-t}P_1\ne \undim  P_1$. That is, $C_{\hspace{-.5pt}H}^{-t}(\undim  P_1) \neq \undim  P_1$. \vspace{.5pt} Suppose now that $m(1)+1\le t$. Then $0\le t-(m(1)+1)\le m(\rho(1))$. Therefore,
$$\tau_{\hspace{-.5pt}_D}^{-t} \hspace{-1.3pt} P_1 =
\tau_{\hspace{-.5pt}_D}^{-(t-(m(1)+1))} (\tau_{\hspace{-.5pt}_D}^{-(m(1)+1)} \hspace{-1pt} P_1) = \tau_{\hspace{-.5pt}_D}^{-(t-(m(1)+1))} \hspace{-1pt} (P_{\rho(1)}[1]) = M[1],\vspace{-1pt}$$
where $M=\tau^{-(t-(m(1)+1))} \hspace{-1.3pt} P_{\rho(1)} \in \GaH$. Applying again Proposition \ref{Cox_tau}, we see that \vspace{-1pt} $$C_{\hspace{-.5pt}H}^{-t}(\undim  P_1)=\undim \, \tau_{\hspace{-.5pt}_D}^{-t} \hspace{-1.3pt} P_1= \undim \, M[1]= - \, \undim  M \ne \undim  P_1. \vspace{-2pt} $$ So, $C_{\hspace{-.5pt}H}$ is indeed of order $a+2$. The proof of the proposition is completed.


\vspace{2pt}

The Coxeter orders are given for each canonical Dynkin diagram with a particular orientation in \cite[Pages 289-290]{ARS}. The following statement says that they are independent of the orientation.

\begin{prop} \label{co_table} Let $H$ be a hereditary artin algebra with a Dynkin ext-quiver $\QH$. Then the Coxeter order $|C_{\hspace{-.5pt}H}|$ for $H$ is independent of the orientation of $\QH$ and given by the following table$\,:$

\vspace{-3pt}

$$\begin{tabular}{|c|c|c|c|c|c|c|c|c|c|} \hline \vspace{-11pt} \\
$\overline{\hspace{-1.6pt}Q\hspace{-.5pt}}_H $  & $\mathbb{A}_n$ & $\mathbb{B}_n$ & $\mathbb{C}_n$ & $\mathbb{D}_n$ & $\mathbb{E}_6$ & $\mathbb{E}_7$ &$\mathbb{E}_8$ & $\mathbb{F}_4$ & $\mathbb{G}_2$  \\  \hline
$|C_{\hspace{-.5pt}H}|$  & $n\!+\!1$ & $2n$ & $2n$ & $\!2(n\!-\!1)\!$ & $12$  & $18$  & $30$  & $12$  & $6$ \\
\hline
\end{tabular}$$
\end{prop}

\vspace{1pt}

\noindent{\it Proof.} Assume that $\OQH$ is a canonical Dynkin diagram; see (\ref{Dyn_diag}). 
Consider the pi-permutation $\rho$ and the mi-index function $m$ for $(\QH)_0$. Suppose first that $\OQH=\mathbb{A}_n$ with $n\ge 1.$ In view of Lemma \ref{An}, we see that $\tau^{-r_{1,n}}\hspace{-1.5pt}P_n=I_1$ and $\tau^{-r_{n,1}}\hspace{-1.5pt}P_1=I_n$, where $r_{i,j}$ is the number of arrows in the reduced walk in $\QH$ from $i$ to $j$. Thus, $m(n)=r_{1,n}$ and $m(1)=r_{n.1}$. Since $r_{n,1}$ is also the number of inverse arrows in the reduced walk in $\QH$ from $1$ to $n$, we have $m(1)+m(n)=n-1.$ Since $\rho(1)=n$, from Propsoition \ref{order_c} we deduce that $|C_H|=(n-1)+2=n+1.$

\vspace{1pt}

Suppose now that $\OQH\ne \mathbb{A}_n$. Then, $\OQH$ contains a vertex $s$, which is of weight $3$ or incident to a non-trivially valued edge. By Lemma \ref{rho_property}, $\rho(s)=s$, and by Proposition \ref{order_c}, $|C_H|=2 m(s)+2$. Let $h_s$ be the extended hammock function on ${\rm Suc}(0,s)$. By Theorem \ref{hamm_neg_section}(1), $m(s)$ is such that $h_s(m(s)+1, s)=-1$ and $h_s(s, j) \ge 0$ for all proper predecessors $(s,j)$ of $(m(s)+1, s)$ in ${\rm Suc}(0,s).$

Let $H'$ be a hereditary artin algebra with $\hspace{2.5pt}\overline{\hspace{-1.5pt}Q}_{\hspace{-.9pt}H'\hspace{-1pt}}=\OQH$. Consider the pi-permutation $\rho'$ and the pi-index function $m'$ for $(Q_{\hspace{-.8pt}H'})_0=(\QH)_0$. Let $h'_s$ be the extended hammock function on ${\rm Suc}_{H'}(0,s)$, the full subquiver of $\mathbb{Z}\hspace{.8pt}Q^{\rm op}_{\hspace{-1pt}H'}$ generated by the successors of $(0,s)$. Since  $\rho'(s)=s$, as previously seen, $|C_{H'}|=2 m'(s)+2$, where $m'(s)$ is such that $h'_s(m'(s)+1, s)=-1$ and $h'_s(s, j) \ge 0$ for all proper predecessors $(s,j)$ of $(m'(s)+1, s)$ in ${\rm Suc}_{H'}(0,s).$

By Proposition \ref{ZD_sec} and Lemma \ref{sec_paths}, the vertices of
the $(0,s)$-source section $\Da$ in $\ZQH$ are $(r_i, i)$, where $i\in (\QH)_0$, and $r_i$ is the number of inverse arrows in the reduced walk in $Q^{\rm op}_{\hspace{-1pt}H}$, that is, the number of arrows in the reduced walk in $\QH$, from $s$ to $i$. Since $\Da$ is a leftmost section in ${\rm Suc}(0,s)$; see (\ref{s_sec_unicity}), the vertices of ${\rm Suc}(0,s)$ are $(r_i+p, i)$ with $p\ge 0$ and $i\in (\QH)_0$.
Similarly, the vertices of the $(0,s)$-source section $\Da'$ in $\mathbb{Z}\hspace{.8pt}Q^{\rm op}_{\hspace{-1pt}H'}$ are $(r'_i, i)$ with $i\in (\hspace{2.5pt}\overline{\hspace{-1.5pt}Q}_{\hspace{-.9pt}H'\hspace{-1pt}})_0=(\OQH)_0$ and $r'_i$ the number of arrows in the reduced walk in $Q_{\hspace{-.8pt}H'}$ from $s$ to $i$,
and those of ${\rm Suc}_{H'}(0,s)$ are $(r'_i+p, i)$ with $i\in (\OQH)_0$ and $p\ge 0$.

\vspace{1pt}

Fix $i,j\in \OQH$. \vspace{-2pt} By Lemma \ref{sec_paths}, $\Da$ contains an arrow $(r_i, i)\to (r_j, j)$ if and only if $\OQH$ (that is $\hspace{2.5pt}\overline{\hspace{-1.5pt}Q}_{\hspace{-.9pt}H'\hspace{-1pt}})$ contains an edge $\hspace{-2pt}\xymatrixcolsep{12pt}\xymatrix{i \ar@{-}[r] & j}\hspace{-2pt}$ if and only if $\Da'$ contains an arrow $(r'_i, i)\to (r'_j, j)$;  \vspace{-3pt} and in this case, the arrows $(r_i, i)\to (r_j, j)$ and $(r'_i, i)\to (r'_j, j)$ have the same valuation as $\hspace{-2pt}\xymatrixcolsep{12pt}\xymatrix{i\ar@{-}[r] & j.}$  \vspace{-2pt} So, we have a valued quiver isomorphism $f: \Da\to \Da'$, sending $(r_i, i)$ to $(r'_i, i)$. Clearly, $f$ extends to a valued translation quiver isomorphism $f: {\rm Suc}(0,s)\to {\rm Suc}_{H'}(0,s)$, sending $(r_i+p, i)$ to $(r'_i+p, i)$, for all $p\ge 0$. Since $r_s=r'_s=0$, we have
$f(p, s)=(p,s)$ for all $p\ge 0$.

\vspace{1pt}

By Definition \ref{hamm_func}, $h'_s(r'_i, i)=h'_s(f(r_i, i))=h_s(r_i, i)$ for all $(r_i,i)\in \Da_0$. Since $h_s$ and $h'_s$ are additive, $h'_s(r'_i+p)=h_s(r_i+p, i)$ for all $p\ge 0$ and $i\in (\OQH)_0$. In particular, $h'_s(m(s)+1, s)=h_s(m(s)+1, s)=-1$. If $(r'_i+p, i)$ is a proper predecessor of $(m_i(s)+1, s)$ in $\mathbb{Z}\hspace{.8pt}Q^{\rm op}_{\hspace{-1pt}H'}$, then $(r_i+p, i)$ is a proper predecessor of $(m_i(s)+1, s)$ in $\ZQH$, and hence, $h'_s(r'_i+p, i)=h_s(r_i+p, i)\ge 0$. Therefore, $m'(s)=m(s),$ and consequently, $|C_{H'}|=|C_H|$. Thus, $|C_{\hspace{-.5pt}H}|$ is independent of the orientation of $\QH$. Using the Coxeter orders given in \cite[Pages 289-290]{ARS}, we obtain the table stated in the proposition. The proof of the proposition is completed.

\begin{rem}

As in the first part of the proof of Proposition \ref{co_table}, one can also use hammocks and hammock functions to explicitly compute the Coxeter order case by case.

\end{rem}

\subsection{\sc Shape of the Auslander-Reiten quiver} Applying our previous results, we can describe the precise shape of $\GaH$ in terms of $\QH$ in the following statement.

\begin{thm}\label{ARQ_Dyn}

Let $H$ be a hereditary artin algebra with a canonical Dynkin ext-quiver $\QH$. Let $P_i$ be the projective modules and $I_i$ are the injective modules in $\GaH$ associated with $i\in (\QH)_0$. Then $\tau^{-m(i)}\hspace{-1pt}P_i=I_{\rho(i)},$ and there exists a convex embedding $\GaH\to \ZQH$, sending $\tau^{-r}P_i$ to $(r, i),$ where $\rho$ is the pi-permutation of $(\QH)_0$ and $m$ is the pi-index function on $(\QH)_0$ given case by case as follows.

\begin{enumerate}[$(1)$]

\vspace{-2pt}

\item If $\OQH=\mathbb{A}_n$ with $n\ge 1$, then $\rho(i)=n+1-i$ and $m(i)=\mathfrak{a}^+(1,i)+\frak{a}^-(1, n+1-i),$ 
for $i\in (\QH)_0$.

\vspace{3pt}

\item If $\OQH=\mathbb{E}_6$, then $\rho=(16)(25)$ and $m(i)=5-\mathfrak{a}^+(i,3)+\mathfrak{a}^+(\rho(i),3)$
for $i\in (\QH)_0$.

\vspace{3pt}

\item If $\OQH=\mathbb{D}_n\, (n \ge 5\ odd \,)$, then $\rho=(12)$
 and $m(i)=n-2$ for all $3\le i\le n$. \vspace{.5pt} Moreover, $m(1)=(n-2)-\mathfrak{a}^+(1,3)+\mathfrak{a}^+(2,3)$ 
 and $m(2)=(n-2)+\mathfrak{a}^+(1,3)-\mathfrak{a}^+(2,3).$

\vspace{3pt}

\item If $\OQH=\mathbb{G}_2,$ $\mathbb{F}_4$, $\mathbb{E}_7$, $\mathbb{E}_8,$
$\mathbb{B}_n \,(n\ge 2)$, $\mathbb{C}_n \,(n\ge 3)$ or $\mathbb{D}_n \,(n \ge 4\ even)$, \vspace{1pt} then we have $\rho(i)=i$ and $m(i)=\frac{1}{2} |C_{\hspace{-.5pt}H}|-1$ for all $i\in \QH$, where $|C_H|$ is the Coxeter order for $H$.

\end{enumerate}\end{thm}

\vspace{-4pt}

\noindent{\it Proof.} By Proposition \ref{H_alg_Dynkin}, we have a canonical embedding $\varphi: \GaH\to \ZQH: \tau^{-r}\hspace{-1pt}P_i\mapsto (r,i)$, which has a convex image in $\ZQH$. By Lemma \ref{exq_permut}, $\tau^{-m(i)}P_i=I_{\rho(i)}$ for all $i\in (\QH)_0$.
Write $r_{i,j}=\mathfrak{a}^+(i, j)$, for $i, j\in (\QH)_0$.

(1) Suppose that $\OQH=\mathbb{A}_n$ with $n\ge 1$. As seen in the proof of Theorem \ref{order_2}, $\rho(i)=n+1-i$ for all $\in (\QH)_0$. Moreover, by Lemma \ref{An}, $\tau^{-r_{1,n}}P_n=I_1$. Thus, $m(n)=r_{1,n}$. Given $i\in (\QH)_0$, in view of Lemma \ref{pi-index}, we see that
$m(i)=m(n)-r_{i,n}+r_{n+1-i, 1}=r_{1,n}-r_{i,n}+r_{n+1-i, n}=r_{1,i}+r_{n+1-i, 1}.$

\vspace{.5pt}

(2) Suppose that $\OQH=\mathbb{E}_6$. In view of the proof of Theorem \ref{order_2}, we see that $\rho=(16)(25)$. In particular, $\rho(3)=3$. Since $|C_H|=12;$ see (\ref{co_table}), we deduce from Proposition \ref{order_c} that $m(3)=5.$ Given $i\in \QH$, by Lemma \ref{pi-index}, $m(i)=5-r_{i,3}+r_{\rho(i),3}.$

(3) Suppose that $\OQH=\mathbb{D}_n$ with $n \ge 5$ odd. In view of the proof of Theorem \ref{order_2}, we see that $\rho=(12).$ In particular, $\rho(i)=i$ for all $3\le i\le n.$ Since $|C_H|=2(n-1)$; see (\ref{co_table}), we deduce from Proposition \ref{order_c} that $m(i)=n-2$, for $3\le i\le n$. Moreover, by Lemma \ref{pi-index}, $m(1)=(n-2)-r_{1,3}+r_{2,3}$ and $m(2)=(n-2)+r_{1,3}-r_{2,3}$.

\vspace{.5pt}

(4) Suppose that $\OQH=\mathbb{G}_2,$ $\mathbb{F}_4$, $\mathbb{E}_7$, $\mathbb{E}_8,$  $\mathbb{B}_n (n\ge 2)$, $\mathbb{C}_n (n\ge 3)$ or $\mathbb{D}_n (n \ge 4 \ {\rm even} \hspace{.5pt})$. By Theorem \ref{order_2}(1), we have $\rho(i)=i$, and by Proposition \ref{order_c}, $|C_{\hspace{-.5pt}H}|=2 m(i)+2$, for all $i\in \QH$. The proof of the theorem is completed.




\vspace{2pt}

\begin{eg}

Let $H$ be a hereditary artin algebre with ext-quiver
$$\xymatrixrowsep{16pt}\xymatrixcolsep{18pt}
\xymatrix{
&&&4&&\\
\QH: & 1\ar[r] & 2 \ar[r]&3\ar[u] \ar[r]&5 & \ar[l] 6.
}$$
Since $\QH$ is of type $\mathbb{E}_6$, by Theorem \ref{ARQ_Dyn}(2), we have
$\rho(1)=6;$ $\rho(2)=5;$ $\rho(3)=3$; $\rho(4)=4$; $\rho(5)=2$ and $\rho(6)=1$. Moreover, $m(1)=m(2)=4;$ $m(3)=m(4)=5$; $m(5)=m(6)=6$. Thus, $\GaH$ is of the following shape:
$$\xymatrixcolsep{23pt}\xymatrixrowsep{20pt}
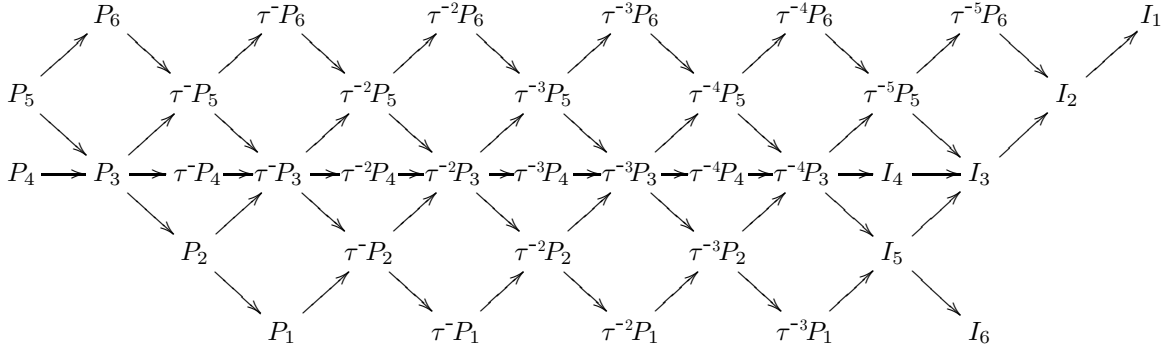
\xymatrix@!=10pt{
&P_6\ar[dr]&&\tau^{\mbox{-}}\hspace{-1.2pt}P_6\ar[dr]&&\tau^{\mbox{-}\scalebox{.6}2}\hspace{-1.2pt}P_6\ar[dr]&&\tau^{\mbox{-}\scalebox{.6}3}\hspace{-1.2pt}P_6\ar[dr]&&\tau^{\mbox{-}\scalebox{.6}4}\hspace{-1.2pt}P_6\ar[dr]&&\tau^{\mbox{-}\scalebox{.6}5}\hspace{-1.2pt}
P_6\ar[dr]&  & \hspace{-2pt} I_1\\
P_5\hspace{-1pt}\ar[dr]\ar[ur]&&\tau^{\mbox{-}}\hspace{-1.2pt}P_5\ar[dr]\ar[ur]&&\tau^{\mbox{-}\scalebox{.6}2}\hspace{-1.2pt}P_5\ar[dr]\ar[ur]&&\tau^{\mbox{-}\scalebox{.6}3}\hspace{-1.2pt}P_5\ar[dr]\ar[ur]&&\tau^{\mbox{-}\scalebox{.6}4}\hspace{-1.2pt}P_5\ar[dr]\ar[ur]&&\tau^{\mbox{-}\scalebox{.6}5}\hspace{-1.2pt}P_5\ar[dr]\ar[ur]&&I_2\ar[ur]&&\\
P_4 \hspace{-1pt} \ar[r]&P_3\ar[dr]\ar[ur]\ar[r]&\tau^{\mbox{-}}\hspace{-1.8pt}P_4 \hspace{-2.5pt} \ar[r]&\hspace{-2.5pt}\tau^{\mbox{-}}\hspace{-1.8pt}P_3\ar[dr]\ar[ur]\ar[r]& \hspace{-2.5pt} \tau^{\mbox{-}\scalebox{.6}2}\hspace{-1.8pt}P_4\hspace{-2.5pt} \ar[r]& \hspace{-2.5pt} \tau^{\mbox{-}\scalebox{.6}2}\hspace{-1.8pt}P_3\ar[dr]\ar[ur]\ar[r] & \hspace{-3pt} \tau^{\mbox{-}\scalebox{.6}3}\hspace{-1.8pt}P_4\hspace{-1.8pt}\ar[r] & \hspace{-2.5pt} \tau^{\mbox{-}\scalebox{.6}3}\hspace{-1.8pt}P_3\hspace{-2pt}\ar[dr]\ar[ur]\ar[r] & \hspace{-2.5pt} \tau^{\mbox{-}\scalebox{.6}4}\hspace{-1.8pt}P_4\hspace{-2pt}\ar[r] & \hspace{-2.5pt}\tau^{\mbox{-}\scalebox{.6}4}\hspace{-1.8pt}P_3\ar[dr]\ar[ur]\ar[r]&I_4\ar[r]&I_3\ar[ur]&&\\
&&P_2\ar[dr]\ar[ur]&&\tau^{\mbox{-}}\hspace{-1.2pt}P_2\ar[dr]\ar[ur]&&\tau^{\mbox{-}\scalebox{.6}2}\hspace{-1.2pt}P_2\ar[dr]\ar[ur]&&\tau^{\mbox{-}\scalebox{.6}3}\hspace{-1.2pt}P_2\ar[dr]\ar[ur]&&I_5\ar[ur]\ar[dr]&&&&\\
&&&P_1\ar[ur]&&\tau^{\mbox{-}}\hspace{-1.2pt}P_1\ar[ur]&&\tau^{\mbox{-}\scalebox{.6}2}\hspace{-1.2pt}P_1\ar[ur]&&\tau^{\mbox{-}\scalebox{.6}3}\hspace{-1.2pt}P_1\ar[ur]&&I_6 }$$

\end{eg}


\vspace{1pt}

\begin{eg} Let $H$ be a hereditary artin algebra with ext-quiver $\xymatrixcolsep{22pt}\xymatrix{Q_H: \hspace{5pt} 1 \ar[r] & 2 \ar[r]^{(1,2)} & 3 & \ar[l] 4.}$ Since $\OQH=\mathbb F_4$, by  Theorem \ref{ARQ_Dyn}(4), $\rho(i)=i$ and $m(i)=5,$ for all $1\le i\le 4$. Thus, $\GaH$ is of the following shape
$$\xymatrixcolsep{24pt}\xymatrixrowsep{25pt}
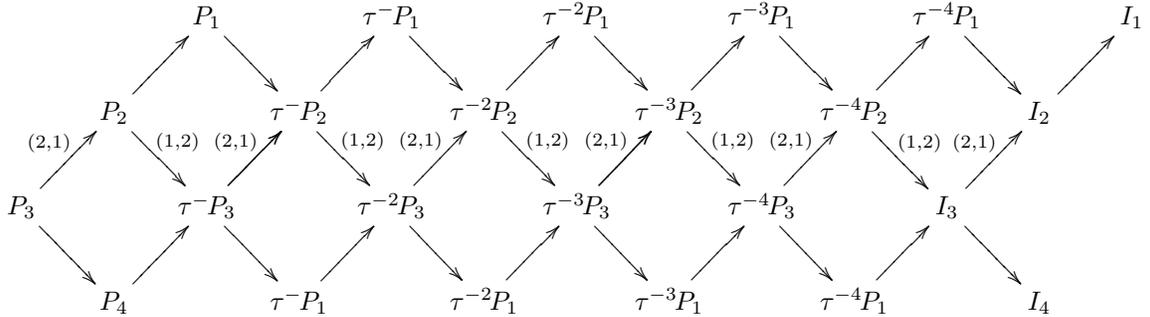
\xymatrix@!=11pt{
&&P_1 \ar[dr] && \tau^{-}\hspace{-1.2pt}P_1 \ar[dr] && \tau^{-2}\hspace{-1.2pt}P_1 \ar[dr] && \tau^{-3}\hspace{-1.2pt} P_1 \ar[dr] && \tau^{-4}\hspace{-1.2pt}P_1 \ar[dr] && I_1 \\
&P_2\ar[dr]^-{\hspace{-4pt}(1,2)} \ar[ur] && \tau^{-}\hspace{-1.2pt}P_2\ar[dr]^-{\hspace{-4pt}(1,2)}\ar[ur] &&\tau^{-2}\hspace{-1.2pt} P_2\ar[dr]^-{\hspace{-4pt}(1,2)}\ar[ur]&& \tau^{-3}\hspace{-1.2pt}P_2\ar[dr]^-{\hspace{-4pt}(1,2)}\ar[ur] &&\tau^{-4}\hspace{-1.2pt} P_2\ar[dr]^-{\hspace{-4pt}(1,2)}\ar[ur]&& I_2\ar[ur] \\
P_3\ar[dr] \ar[ur]^-{(2,1)\hspace{-4pt}} && \tau^{-}\hspace{-1.2pt} P_3 \ar[ur]^-{(2,1)\hspace{-4pt}} \ar[dr] \ar[ur]&&\tau^{-2}\hspace{-1.2pt} P_3\ar[dr] \ar[ur]^-{(2,1)\hspace{-4pt}} && \tau^{-3}\hspace{-1.2pt} P_3 \ar[ur]^-{(2,1)\hspace{-4pt}} \ar[dr] \ar[ur]&&\tau^{-4}\hspace{-1.2pt} P_3\ar[dr] \ar[ur]^-{(2,1)\hspace{-4pt}} &&
I_3\ar[dr] \ar[ur]^-{(2,1)\hspace{-4pt}}  \\
&P_4\ar[ur]&& \tau^{-}\hspace{-1.2pt}P_1\ar[ur]&& \tau^{-2}\hspace{-1.2pt}P_1\ar[ur]&& \tau^{-3}\hspace{-1.2pt}P_1\ar[ur]&& \tau^{-4}\hspace{-1.2pt}P_1\ar[ur]&& I_4 }$$

\end{eg}

\vspace{3pt}

The following statement is interesting in its own right.

\begin{prop}

Let $H$ be a hereditary artin algebra of Dynkin type. If $P$ is a projective module and $I$ is an injective module in $\GaH$, then $I$ is a successor of $P$ in $\GaH$.

\end{prop}

\vspace{-4pt}

\noindent{\it Proof.} Let $P$ be a projective module and $I$ an injective module in $\GaH$. By Proposition \ref{special section}, $\GaH$ contains a $P$-source section $\Da$. Then $I=\tau^s\hspace{-1pt}M$ for some $M\in \Da$ and $s\in \Z$. If $s>0$, then $\tau^{-}\hspace{-1pt}I=\tau^{s-1}M\in \GaH$, which gives a contradiction since $\tau^{-}\hspace{-1pt}I =0$. Thus, $s\le 0$, and consequently, $I$ is a successor of $P$. The proof of the proposition is completed.

\subsection{\sc Applications} Our first application is to compute the number of non-isomorphic indecomposable modules in $\mmod H$. Applying a series of results in \cite{PAu} and \cite{DRi}, we see that this number was given case by case in \cite{DRi} and \cite{Gab}, and also by the number of positive roots of Dynkin diagrams; see \cite{DR, Gab}, which was given in \cite{Cox}. However, using our previous results, we obtain this number in a direct and easy way.

\begin{thm}\label{rep_nb}

Let $H$ be a hereditary artin algebra with a Dynkin ext-quiver $\QH$ and Coxeter order $|C_{\hspace{-.5pt}H}|$. Then, the number of non-isomorphic indecomposable modules in $\mmod H$ is equal to $$\frac{1}{2}n |C_{\hspace{-.5pt}H}|,$$ where $n$ is the number of non-isomorphic simple modules in $\mmod H$.

\end{thm}

\vspace{-4pt}

\noindent{\it Proof.} Consider the pi-permutation $\rho$ and the pi-index function $m$ for $(\QH)_0=\{1, \ldots, n\}$. Let $P_i$ be the projective module and $I_i$ the injective module in $\GaH$ associated with $i\in (\QH)_0$. Since $\tau^{-m(i)}=I_{\rho(i)}$, the $\tau$-orbit of $P_i$ contains exactly $m(i)+1$ modules. Let $t$ be the number of modules in $\GaH$. Since the $P_i$ generate a section in $\GaH$, we see that $t=\sum_{i=1}^n (m(i)+1)$. And since $\rho$ is a permutation, $t=\sum_{i=1}^n (m(\rho(i)+1)$. \vspace{.5pt} Applying Proposition \ref{order_c} yields $2t=\sum_{i=1}^n (m(i)+m(\rho(i))+2)=n |C_{\hspace{-.5pt}H}|.$ The proof of the theorem is completed.


\vspace{3pt}

Now, we shall compute the nilpotency of the radical of $\mmod H$. This has been essentially done by Zacharia using preprojective partitions; see \cite[Section 4]{Zac}. However, we shall provide an alternative approach. Since $\QH$ is a Dynkin diagram, by Propositions \ref{uni_sec_path} and \ref{H_alg_Dynkin}, any two parallel paths in $\GaH$ and $\Ga_{\hspace{-1pt}D^{\hspace{.4pt}b\hspace{-.5pt}}(\mmod H)}$ have the same length. So we may introduce the following definition.

\begin{defn}\label{dist}

Let $H$ be a hereditary artin algebra of Dynkin type. 

\begin{enumerate}[$(1)$]

\vspace{-1pt}

\item If $\GaH$ contains a path from $M$ to $N$, then we define the {\it distance} between $M$ and $N,$ written as ${\rm dist}(M, N),$ to be the length of any path from $M$ to $N$ in $\GaH$.

\item If $\Ga_{\hspace{-1pt}D^{\hspace{.4pt}b\hspace{-.5pt}}(\mmod H)}$ contains a path from $M^\cdt$ to $N^\cdt$, then we define the {\it distance} between $M^\cdt$ and $N^\cdt,$ written as ${\rm dist}(M^\cdt, N^\cdt\hspace{-2pt})$, to be the length of any path  from $M^\cdt$ to $N^\cdt$ in $\Ga_{\hspace{-1pt}D^{\hspace{.4pt}b\hspace{-.5pt}}(\mmod H)}$.

\end{enumerate} \end{defn}

The following statement allows us to easily determine the depth of any map in $\mmod H$.

\begin{lem} \label{dpt_dist}

Let $H$ be a hereditary artin algebra of Dynkin type. If $f: M\to N$ is a non-zero map with $M, N\in \GaH$, then ${\rm dp}(f)={\rm dist}(M, N).$

\end{lem}

\vspace{-4pt}

\noindent{\it Proof.} Let $f: M\to N$ be a non-zero map with $M, N\in \GaH$. Since $H$ is representation-finite, $f\not\in \rad^s(M, N)$ for some $s> 0$; see \cite[(V.7.6)]{ARS}. Thus, ${\rm dp}(f)=t<\infty$. If $t=0$, then $M= N$ and ${\rm dist}(M,N)=0$. Otherwise, by Lemma \ref{t_irreducibles}, $\GaH$ contains a path from $M$ to $N$ of length $t$. Hence, ${\rm dist}(M,N)=t.$ The proof of the lemma is completed.

\vspace{3pt}

We are ready to provide a new proof for Zacharia's theorem as follows; see \cite[(4.11)]{Zac}.

\begin{thm}\label{Mod_nil}

Let $H$ be a hereditary artin algebra with a Dynkin ext-quiver $\QH$ and  Coxeter order $|C_{\hspace{-.5pt}H}|$. Then, the radical of $\mmod H$ is nilpotent of nilpotency $|C_{\hspace{-.5pt}H}|-1$.

\end{thm}

\vspace{-4pt}

\noindent{\it Proof.} Given $i\in (\QH)_0=\{1, \ldots, n\},$ consider the simple module $S_i\in \GaH$ with a projective cover $\pi_i: P_i\to S_i$ and an injective envelope $\iota_i: S_i\to I_i$, where $P_i, I_i\in \GaH$.
Applying Theorem 2.7 in \cite{CL} followed by Lemma \ref{dpt_dist}, we see that the nilpotency of $\rad(\mmod H)$ is equal to \vspace{-2pt}
$$\sup\{\ndp(\iota_1 \pi_1), \ldots, \ndp(\iota_n \pi_n)\}+1=\sup\{{\rm dist}(P_1, I_1), \ldots, {\rm dist}P_n,I_n)\} +1.\vspace{-2pt}$$

Thus, it suffices to show that $\text{dist}(P_i, I_i)=|C_{\hspace{-.5pt}H}|-2$, for any $i\in (\QH)_0$. In fact,
since $\GaH$ is a convex valued translation subquiver of $\Ga_{\hspace{-1pt}D^b(\mmod H)}$; see (\ref{Der_ARQ}), we have $\text{dist}(M[0], N[0])=\text{dist}(M, N),$ for all $M, N\in \GaH$. And since $\tau_{\hspace{-.5pt}_D}^-I_i[0]=P_i[1]$; see \cite[(7.2)]{BLP}, we deduce \vspace{-2pt} $$\text{dist}(P_i[0],P_i[1])=\text{dist}(P_i[0],I_i[0])+\text{dist}(I_i[0],P_i[1])=\text{dist}(P_i, I_i)+2.\vspace{-2pt}$$
As a consequence, we obtain \vspace{-2pt} $$\text{dist}(P_i[0],P_i[2])= \text{dist}(P_i[0],P_i[1]) + \text{dist}(P_i[1],P_i[2])=2\hspace{.5pt}\text{dist}(P_i[0],P_i[1])=2(\text{dist}(P_i, I_i)+2).\vspace{-2pt}$$
On the other hand, by Proposition \ref{order_c},
$\text{dist}(P_i[0],P_i[2]) =\text{dist}(P_i[0], \tau_{\hspace{-.5pt}_D}^{-|C_{\hspace{-.5pt}H}|}P_i[0]) = 2\hspace{.5pt}|C_{\hspace{-.5pt}H}|.$ That is, $|C_{\hspace{-.5pt}H}|=\text{dist}(P_i, I_i)+2$. The proof of the theorem is completed.

\vspace{3pt}

Next, we shall study the radical of the bounded derived category $D^{\hspace{.5pt}b\hspace{-.6pt}}(\mmod H)$. Being a triangulated category, $D^{\hspace{.5pt}b\hspace{-.6pt}}(\mmod H)$ coincides with its projectively stable category and its injectively stable category; see \cite[(2.4)]{LiN}. Since $D^{\hspace{.5pt}b\hspace{-.6pt}}(\mmod H)$ has almost split triangles, the Auslander-Reiten translation $\tau_{\hspace{-1pt}_D}$ of $\Ga_{D^{\hspace{.5pt}b\hspace{-.6pt}}(\mmod H)}$ induces an auto-equivalence $\tau_{\hspace{-1pt}_D}: D^{\hspace{.5pt}b\hspace{-.6pt}}(\mmod H)\to D^{\hspace{.5pt}b\hspace{-.6pt}}(\mmod H)$; see \cite[(4.10)]{LiN}.

\begin{thm}\label{Der_nil}
Let $H$ be a hereditary algebra of Dynkin type with Coxeter order $|C_{\hspace{-.5pt}H}|$. Then the radical of $D^{\hspace{.5pt}b\hspace{-.6pt}}(\mmod H)$ is nilpotent of nilpotency $|C_{\hspace{-.5pt}H}|-1$.
\end{thm}

\noindent{\it Proof.} Write $\mathscr{D}=D^{\hspace{.5pt}b\hspace{-.6pt}}(\mmod H)$. Since $\mmod H$ is a convex subcategory in $\mathscr{D}$, by Theorem \ref{Mod_nil}, $\rad^{|C_H|-2}(\mathscr{D})\ne 0.$ Suppose that $\rad^r(\mathscr{D})\ne 0$ with $r\ge 0$. Then, $\rad^r_{\mathscr{D}}(M, N[s])\ne 0$, for some $M, N\in \GaH$ and $s\in \Z$. Write $M=\tau^{-t}P$ for some $t\ge 0$, where $P$ is a projective module in $\GaH$. Applying the equivalence $\tau_{\hspace{-1pt}_D}$, \vspace{1pt} we obtain $\rad^r_{\mathscr{D}}(P[0], (\tau_{\hspace{-1pt}_D}N[0])[s])\ne 0$. Since $P$ is projective, $(\tau_{\hspace{-1pt}_D}N[0])[s]=L[0]$ for some $L\in \GaH$. \vspace{1pt} Then, $\rad^r(P, L)\ne 0$, and by Theorem \ref{Mod_nil}, $r\le |C_H|-2.$ \vspace{.5pt} So, $\rad(\mathscr{D})$ is nilpotent of nilpotency $|C_{\hspace{-.5pt}H}|-1$. The proof of the theorem is completed.

\vspace{3pt}

We conclude with the cluster category associated with $H$. Let $\mathscr{D}^{\hspace{.5pt}b\hspace{-.7pt}}(H)$ be a skeleton of $D^{\hspace{.4pt}b\hspace{-.5pt}}(\mmod H)$, containing the complexes in $\Ga_{\hspace{-1.5pt}D^{\hspace{.4pt}b\hspace{-.5pt}}(\mmod H)}$. Then, $\mathscr{D}^{\hspace{.5pt}b\hspace{-.7pt}}(H)$ is a
Hom-finite Krull-Schmidt $R$-category, which has almost split sequences. Note that the Auslander-Reiten quiver $\Ga_{\hspace{-.8pt}\mathscr{D}^{\hspace{.4pt}b\hspace{-.5pt}}(H)}$ of $\mathscr{D}^{\hspace{.4pt}b\hspace{-.5pt}}(H)$ coincides with $\Ga_{\hspace{-1.5pt}D^{\hspace{.5pt}b\hspace{-.7pt}}(\mmod H)}$, and its
Auslander-Reiten translation $\tau_{\hspace{-.6pt}_\mathscr{D}}$ extends to an automorphism of $\mathscr{D}^{\hspace{.5pt}b\hspace{-.7pt}}(H)$. Thus, we have an automorphism $F=\tau_{\hspace{-.6pt}_\mathscr{D}}^{-1}\circ [1]$ of  $\mathscr{D}^{\hspace{.5pt}b\hspace{-.7pt}}(H)$ such, for any indecomposable objects $M^\pdt, N^\pdt\in \mathscr{D}^{\hspace{.5pt}b\hspace{-.7pt}}(H)$, that $F^{\hspace{.5pt}p\hspace{-.6pt}}(M^\pdt)\not\cong M^\pdt$ for all $p\ne 0$ and $\Hom_{\mathscr{D}^{\hspace{.4pt}b\hspace{-.5pt}}(H)}(M^\pdt, F^{\hspace{.5pt}p\hspace{-.6pt}}(N^\pdt))=0$ for all but finitely many $p\in \Z.$ Thus, the action of the group $\cF$ generated by $F$ on $\mathscr{D}^{\hspace{.5pt}b\hspace{-.7pt}}(H)$ is free and locally bounded. As did in \cite{BMRRT}, we define the {\it cluster category} associated with $H$ to be the orbit category $$\mathscr{D}^{\hspace{.4pt}b\hspace{-.5pt}}(H)/\mathcal{F}=:\mathscr{C}_H,$$
which is Hom-finite and Krull-Schmidt; see \ref{orbit_cat_prop}.

\begin{thm}\label{clus_nil}

Let $H$ be a hereditary artin algebra of Dynkin type with Coxeter order $|C_H|$, and let $\mathscr{C}_H$ be the cluster category associated with $H$.

\begin{enumerate}[$(1)$]

\vspace{-1.5pt}

\item The radical of $\mathscr{C}_H$ is nilpotent of nilpotency $|C_{\hspace{-.5pt}H}|-1$.

\vspace{.5pt}

\item The number of non-isomorphich indecomposable objects in $\mathscr{C}_H$ is equal to \vspace{-1pt} $$\frac{1}{2}n(|C_H|+2),\vspace{-1pt}$$ where $n$ is the number of non-isomorphic simple modules in $\mmod H$.

\end{enumerate} \end{thm}

\vspace{-4pt}

\noindent{\it Proof.} First, Statement (1) follows immediately from Proposition \ref{orbit_cat_prop}(1) and Theorem \ref{Der_nil}. Let  $P_1, \ldots, P_n$ be the projective modules and  $I_1, \ldots, I_n$ the injective modules  in $\GaH$ with ${\rm top} P_i={\rm soc} I_i$. Recall that $F=\tau_{\hspace{-.6pt}_\mathscr{D}}^{-1}\circ [1]$, where $\tau_{\hspace{-.6pt}_\mathscr{D}}$ is the Auslander-Reiten translation for $\mathscr{D}^{\hspace{.5pt}b\hspace{-.7pt}}(H)$. Given $p\in \Z$, write $\mathcal{S}^{(p)}=\{ M[p] \mid M\in \GaH\} \cup \{P_1[p+1], \ldots, P_n[p+1]\}$ and $\mathcal{T}^{(p)}=F^p(\mathcal{S}^{(0)})$. Since $\tau_{\hspace{-.6pt}_\mathscr{D}}^{-}I_i=P_i[1]$, we obtain $F(I_i)=P_i[2]$ for $1\le i\le n$. Thus, $\mathcal{T}^{(1)}=F(\mathcal{S}^{(0)}).$ In general,
$\mathcal{T}^{(p)}=F^p(\mathcal{S}^{(0)})$ for all $p\in \Z$. Since $\mathcal{F}=\{F^p \mid p\in \Z\}$, we see that $\mathcal{S}^{(0)}$ is a complete set of representatives of the $\mathcal{F}$-orbits of non-isomorphic indecomposable objects in $\mathscr{D}^{\hspace{.5pt}b\hspace{-.7pt}}(H)$. Now, Statement (2) follows from  Proposition \ref{orbit_cat_prop} and Theorem \ref{rep_nb}. The proof of the theorem is completed.

\begin{rem}

The number given in Theorem \ref{clus_nil}(2) coincides with the number cluster variables in the cluster algebra associated with $\QH$; see \cite[(5.9.1)]{FWZ}.

\end{rem}


\end{document}